\documentclass[12pt,twoside,english]{article}
\usepackage{ifthen}                % if then else Abfragen
\newboolean{boldeqnitalic}
\setboolean{boldeqnitalic}{true} 

\usepackage{caption}
\usepackage[intlimits] {amsmath}   % muss for defjs# stehen 
\usepackage{defjs2}                % eigene Definitionen    
\usepackage{epsfig}
\usepackage{psfig}
\usepackage{graphicx}
\usepackage{psfrag}
\usepackage[usenames]{color}
\usepackage{colortbl}
\usepackage{ifthen}
\usepackage{fancyhdr}
\usepackage{rotating}
\usepackage{multirow}
\usepackage{booktabs}              % Dicke Tabellenlinien
\usepackage{amssymb}
\usepackage{amsbsy}
\usepackage{bm}                    % fette Schrift in Formeln
\usepackage{babel}
\usepackage{theorem}
\usepackage{upgreek}  % gerade (roman) griechische Buchstaben z.B. \upalpha
\usepackage{pifont}   % zusaetzliche Schriften, z.B. eingekreiste Zahlen
\usepackage{mathrsfs}
\usepackage{datetime}
\usepackage{tikz}
\usetikzlibrary{matrix}
\usepackage[all]{xy}
\usepackage{rotating}
\usepackage{subfigure}
\definecolor{dunkelgrau}{rgb}{0.8,0.8,0.8}
\definecolor{hellgrau}{rgb}{0.90,0.90,0.90} %... slightly darker 
\newcommand\defi{\mathrel{\overset{\makebox[0pt]{\mbox{\normalfont\scriptsize\sffamily def}}}{=}}}
\usepackage[authoryear]{natbib}   
\fancypagestyle{plain}{%
  \fancyhead{}%
  \fancyfoot[c]{\sffamily\thepage}%
}
\makeatletter                           % Leere Fuellseiten
\def\cleardoublepage{\clearpage\if@twoside \ifodd\c@page\else
  \hbox{}
  \vspace*{\fill}
  \thispagestyle{empty}
  \newpage
  \if@twocolumn\hbox{}\newpage\fi\fi\fi}
\makeatother
\begin{document}
\unitlength1.0cm
\frenchspacing

\thispagestyle{empty}
 
\vspace{-3mm}
\begin{center}
  {\bf \large Convergence and Error Analysis of FE-HMM/FE$^2$}\\[2mm]
  {\bf \large for Energetically Consistent Micro-Coupling Conditions} \\[2mm]
  {\bf \large in Linear Elastic Solids} \\[2mm]
\end{center}

\vspace{4mm}
\ce{Andreas Fischer, Bernhard Eidel$^{\ast}$}

\vspace{4mm}

\ce{\small Heisenberg-Group, Institute of Mechanics, Department Mechanical Engineering} 
\ce{\small University Siegen, 57068 Siegen, Paul-Bonatz-Str. 9-11, Germany} 
\ce{\small $^{\ast}$e-mail: bernhard.eidel@uni-siegen.de, phone: +49 271 740 2224, fax: +49 271 740 2436} 
\vspace{2mm}

%\bigskip
%
%\begin{center}
%{\bf \large Highlights}
%
%\bigskip
%
%{\footnotesize
%\begin{minipage}{15.5cm} 
%
%\vspace*{-2mm}
%
%\begin{itemize}
% \item Micro and macro convergence analysis for periodic, constant traction and linear displacement coupling conditions for linear and quadratic shape functions. \\[-6mm]
% \item Comparison of two methods for the constant traction coupling condition. \\[-6mm]
% \item {\color{black}Superconvergent, recovery-type a posteriori error estimation on the macroscale accurately captures the macro discretization error.}       \\[-6mm]
% \item Assessment of optimal uniform micro-/macro mesh refinement strategy: full order for minimal effort.    \\[-6mm]
%\end{itemize}
%
%\end{minipage}
%}
%\end{center}

\bigskip

\begin{center}
{\bf \large Abstract}

\bigskip

{\footnotesize
\begin{minipage}{14.5cm} 
\noindent
A cornerstone of numerical homogenization is the equivalence of the microscopic and the macroscopic energy densities, which is referred to as Hill-Mandel condition. Among these coupling conditions, the cases of periodic, linear displacement and constant traction conditions are most prominent in engineering applications. While the stiffness hierarchy of these coupling conditions is a theoretically established and numerically verified result, very little is known about the numerical errors and convergence properties for each of them in various norms. The present work addresses these aspects both on the macroscale and the microscale for linear as well as quadratic finite element shape functions.
The analysis addresses aspects of (i) regularity and how its loss affects the convergence behavior on both scales compared with the a priori estimates, of (ii) error propagation from micro to macro and of (iii) optimal micro-macro mesh refinement strategy. For constant traction conditions two different approaches are compared. The performance of a {\color{black}recovery}-type error estimation based on superconvergence is assessed. All results of the present work are valid for both the Finite Element Heterogeneous Multiscale Method FE-HMM and for FE$^2$.  
\end{minipage}
}
\end{center}

{\bf Keywords:}
Computational homogenization; Macro-to-micro modeling; Convergence analysis; Error estimation; Finite element methods \hfill 
%{\color{black} vers.\,\today \, at \currenttime}\\
%\hspace*{6.3cm} 
%\\[2mm]
%vers.\,\today \, at \currenttime\\
%{\color{red}{\bf Agenda:} micro-, macro-: one word or two words? $\mathbb{A}^0$ or what -- look to the first version.} 

%----------------------------------------------------------------------------------------------------------------

\section{Introduction}
\label{sec:intro}

The overall aim of computational homogenization is to compute effective properties of microheterogeneous materials.
This can be done in an a priori fashion in that effective properties are the result of pre-computations, which enables the identification of parameters in a constitutive law. In strongly nonlinear regimes as for inelastic material behavior, effective properties are rather calculated on the fly in direct micro-macro transitions. The first approach can be seen as a sequential or hierarchical multiscale method, the second variant as a concurrent multiscale method. 
In either case the methods aim at an trade-off of accuracy with efficiency by capturing the real microstructure along with a sampling of it in volumes of confined size. If the sampling regions are statistically representative, they are referred to as representative volume element (RVE). While an RVE is uncritical to identify for periodic microstructures, the proper choice of the RVE for non-periodic microstructures and random heterogeneous materials is still an item of ongoing research, see e.g. \cite{Ostoja-Starzewski2006}, \cite{Doskar-Zeman-etal-2018} and references therein. 

Among concurrent two-scale methods with direct micro macro transitions the so-called FE$^2$ has been advanced in different directions of non-linear solid mechanics and used in a multitude of engineering applications, \cite{MichelMoulinecSuquet1999}, \cite{Miehe-etal-1999a}, \cite{Miehe-etal-1999b}, \cite{Fish-etal-damage1999}, \cite{FeyelChaboche2000}, \cite{Kouznetsova-etal2001}, \cite{Kouznetsova-etal2002}, \cite{Kanit-etal-2003}, \cite{Peric-etal-2010}, \cite{GeersKouznetsovaBrekelmans2010a}, \cite{GeersKouznetsovaBrekelmans2010b}, \cite{LarssonRunesson-etal-2011}, \cite{Schröder2014}, \cite{Saeb-Steinmann-Javili-2016}, \cite{JaviliSaebSteinmann2017}. More recently, the FE$^2$ framework has been extended to transient computational homogenization \cite{Pham-etal2013} and to the elastodynamics of metamaterials and of phononic crystals \cite{Sridhar-etal2018}. 
 
In spite of these advancements, there is a remarkable lack of knowledge about the mathematical properties of FE$^2$ in a fully or semi-discrete setting as a two-scale finite element method. Here, the so-called Finite Element Heterogeneous Multiscale Method FE-HMM has made substantial contributions providing unified error estimates that comprise the macro error, the total micro error, and the modeling error \cite{Assyr2005}, \cite{AssyrSchwab2005}, \cite{E-Ming-Zhang-2005}, \cite{Ohlberger2005}, \cite{Assyr2009}, \cite{Assyr-etal2012}. This advances the understanding and knowledge of FE-HMM, although the obtained results are currently restricted to linear problems; in solid mechanics to the purely linear setting of linear elastic material behavior along with geometrical linearity, \cite{Assyr2006}, \cite{JeckerAbdulle2016}. Beyond the theoretical relevance of a priori error estimates they are of practical relevance, since they prescribe, how a uniform micro-macro refinement strategy shall be carried out in order to achieve optimal convergence for minimal computational efforts. 

FE$^2$ and FE-HMM have been developped independently from each other and on almost parallel avenues without joint links or crossroads, FE$^2$ in mechanics, FE-HMM in mathematics as an off-spring of the very general Heterogeneous Mulitscale Method HMM, \cite{E-Engquist-2003}, \cite{E-Engquist-Huang-2003}. Quite recently it was shown that despite minor differences in the numerical setup\footnote{While the micro-macro stiffness transfer in FE$^2$ refers to the homogenized tangent moduli, FE-HMM refers to the microstiffness matrix along with a transfer operator.} the two methods are equivalent and, as a consequence, the a priori estimates of FE-HMM equally apply for FE$^2$, \cite{EidelFischer2018}.
 
A theoretical sound and commonly accepted link between the scales is Hill's postulate of energy equivalence between micro and macro energy densities, \cite{Hill1963}, \cite{Hill1972}. It is applied in both FE$^2$ and FE-HMM, for the latter method without reference to Hill's work, \cite{E-Engquist-2003}. Several micro boundary conditions (BC) fulfill the postulate, among them (i) the linear displacement BC also referred to as kinematically uniform BC (KUBC) or Dirichlet BC (DBC), the (ii) constant traction BC (TBC), also called static uniform BC (SUBC) or Neumann BC, and (iii) the periodic BC (PBC). These three micro-coupling conditions, which are frequently called canonical in view of their practical relevance can be ordered according to their stiffness in that PBC has its lower bound by Neumann BC and its upper bound by Dirichlet BC.
  
The present work aims to advance the understanding of FE$^2$ and FE-HMM in their numerical characteristics. For that aim we address the following aspects for the particular case of linear elasticity in a geometrical linear frame, since this setting allows for the direct comparison with the a priori estimates of FE-HMM: 
\begin{enumerate}
 \item The numerical error and its convergence is analyzed for the following set of energetically consistent BCs, for (i) KUBC/Dirichlet, for (ii) PBC and for (iii) SUBC/Neumann. 
       \\[1mm]
       While the hierarchy of stiffnesses for conditions (i)--(iii) is a theoretically established result that was verified in numerous numerical simulations, \cite{Suquet-1987}, \cite{Miehe-2003}, \cite{Kanit-etal-2003}, \cite{Peric-etal-2010}, the convergence for the different coupling conditions is largely unexplored; an exception is \cite{Yue-E2007} for the scalar-valued field problem of transport/conductivity. One of the guiding questions is whether there are significant differences in errors and convergence between the micro coupling conditions, and if so, whether they are generally valid similar to the mentioned stiffness hierarchy. Moreover the analysis compares the measured convergence orders against the nominal a priori estimates in different norms. In this context, a discrimination between micro error convergence on the microscale and on the macroscale is relevant, where the latter implies a micro-to-macro error propagation and a somewhat unusual convergence estimate. In either case the regularity of the boundary value problems both on the macro- and the microscale --and its loss due to singularities-- is of importance.      \\
       Beyond the analysis of different coupling conditions, we compare two different numerical approaches for constant traction conditions; the approach based on a simple mass-type diagonal perturbation of the stiffness matrix introduced by \cite{Miehe-Koch-2002}\footnote{similarly used at finite strains in \cite{Miehe-2003}.} with the more recent approach of \cite{JaviliSaebSteinmann2017}. 
       
 \item Error estimation based on the Superconvergent Patch Recovery (SPR) and its validation by comparison with the exact error.
       \\[1mm]
       The true numerical error can be calculated quite accurately by comparison of the apprximate solution with a reference solution on extremely fine grids. In engineering practice however, suchlike overkill solutions along with error calculation in a postprocessing step are not feasible. For that reason error estimation provides an efficient way to analyze accuracy on-the-fly given that the error estimation is validated.
       
 \item For the above analyses of error computation and estimation a set of benchmark problems is considered. They span the range from highly regular boundary value problems (BVP) up to singularity-dominated cases for both the macro as well as the microscale in $n_{dim}=2$.
\end{enumerate}

%---------------------------------------------------------------------------------------------------------

%---------------------------------------------------------------------------------------------------------
\section{The finite element heterogeneous multiscale method FE-HMM in a nutshell}
\label{sec:FEHMM}
%---------------------------------------------------------------------------------------------------------

To put things into perspective and for ready reference this section outlines an FE-HMM formulation for linear elasticity cf. \cite{EidelFischer2018}.

\subsection{Model problem of linear elasticity} 
\label{subsec:ModelProblemLinearElasticity}

We consider a body $\mathcal{B}$, a bounded subset of $\mathbb{R}^d$, $n_{dim}=2,3$, with boundary 
$\partial \mathcal{B} = \partial \mathcal{B}_D \cup \partial \mathcal{B}_N$ where the Dirichlet boundary $\partial\mathcal{B}_D$ and the Neumann boundary $\partial\mathcal{B}_N$ are disjoint sets.
The closure of the body $\mathcal{B}$ is denoted by $\overline{\mathcal{B}}$.
The body, which exhibits an inhomogeneous microstructure, is subject to body forces $\bm f$ and surface tractions $\bar{\bm t}$ and in static equilibrium.
 
\subsubsection{The microproblem}
\label{subsec:StrongFormLinearElasticity}

The displacement $\bm u^{\epsilon} = (u_1^{\epsilon}, \ldots, u_{n_{dim}}^{\epsilon})$ of the body is given by the solution of the system 
\begin{equation}
\label{eq:StrongFormMicro-2}
\renewcommand{\arraystretch}{1.6}
%\left.
   \begin{array}{rcl}
   {\color{black}- \, \sigma^{\epsilon}_{ij,j}} 
   &=& f_i   \qquad \mbox{in} \quad \mathcal{B}  \\
   u_i^{\epsilon} 
   &=& \bar{u}_i \qquad \mbox{on} \quad \partial \mathcal{B}_{D} \\
    \sigma^{\epsilon}_{ij} \, n_j 
   &=& \bar{t}_i   \qquad \, \, \mbox{on} \quad \partial \mathcal{B}_{N} \\
\end{array}
%\right\} \; 
\end{equation}
 The constitutive law is assumed to be linear elastic $\sigma^{\epsilon}_{ij}=\mathbb{A}^{\epsilon}_{ijlm}\, {\color{black}\varepsilon_{lm}}$ where $\mathbb{A}^{\epsilon}_{ijlm}$ is the fourth order elasticity tensor and $\varepsilon_{ij}$ the infinitesimal strain tensor with $\varepsilon_{ij}(\bm u^{\epsilon}) = 1/2 \left(u_{i,j}^{\epsilon} + u_{j,i}^{\epsilon}\right)$ or more compact, $\bm \varepsilon (\bm u^{\epsilon}) = \bm L \, \bm u^{\epsilon}$ with the linear differential operator $\bm L$.
 Superscript $\epsilon$ throughout indicates the dependency of suchlike marked quantities on the heterogeneity of the elastic material. 

 In \eqref{eq:StrongFormMicro-2}$_{3}$, $\bm n=(n_1, \ldots, n_{dim})$ is the unit outward normal to $\partial \mathcal{B}$.  

Multiplying the strong form \eqref{eq:StrongFormMicro-2} by a test function $\bm v \in \mathcal{V}$ and the application of the Green formula yield the variational form:

Find $\bm u^{\epsilon}$ such that 
\begin{equation}
\label{eq:weak-form-for-u-epsilon}
 B_{\epsilon} (\bm u^{\epsilon}, \bm v) \defi 
  \int_{\mathcal{B}} {\color{black}\mathbb{\sigma}^{\epsilon}(\bm u^{\epsilon})} : \bm \varepsilon(\bm v) \, dV 
           = \int_{\mathcal{B}} \bm f \cdot \bm v \, dV \, + \,  \int_{\partial \mathcal{B}_N} \bar{\bm t} \cdot \bm v \,dA 
  \defi \bm F(\bm v)\, ,
\end{equation}
which must hold for all $\bm v \in \mathcal{V}$, where $\mathcal{V}$ is the space of admissible (virtual) displacements that fulfill homogeneous Dirichlet BC
\begin{equation}
\label{eq:HilbertSpaceV}
  \mathcal{V}=\{\bm v; \bm v \in H^1(\mathcal{B})^{n_{dim}}, \bm v|_{\partial \mathcal{B}_D} = \bm 0 \} \, .
\end{equation}

%---------------------------------------------------------------------------------------------------------
\subsubsection{The macroproblem}
\label{subsec:Variational-FE-HMM-macro}

The strong form of the macroscopic/homogenized boundary value problem (BVP) reads 
\begin{equation}
\label{eq:Homogenized-Strong-Form}
\renewcommand{\arraystretch}{1.6}
%\left.
\begin{array}{rcl}
 {\color{black}- \, \sigma^{0}_{ij,j}} &=& \langle f_i \rangle  \qquad \, \, \, \mbox{in} \quad \mathcal{B}  \\
%   \label{eq:PDE-homogenized} \\
  u_i^{0} &=& \langle \bar{u}_i\rangle_{\Gamma} \qquad \mbox{on} \quad \partial \mathcal{B}_{D}    \\
%    \label{eq:BCs-Dirichlet-homogenized} \\
  \sigma^{0}_{ij} \, n_j &=& \langle \bar{t}_i\rangle_{\Gamma}   \qquad \, \, \mbox{on} \quad \partial \mathcal{B}_{N} \\
%  \label{eq:BCs-Neumann-homogenized}
\end{array}
%\right\} \; 
\end{equation}
for a derivation see \cite{EidelFischer2018}. The macroscopic displacement is denoted by $u_{i}^0$ and $\mathbb{A}^{0}$ is the homogenized elasticity tensor. In \eqref{eq:Homogenized-Strong-Form}$_{1,3}$ $\sigma^{0}_{ij}$ is the macroscopic stress obtained by a volume average over the microdomain. 

The values for the Dirichlet as well as Neumann BC in \eqref{eq:Homogenized-Strong-Form}$_{2,3}$ are obtained by surface averages of corresponding BC in \eqref{eq:StrongFormMicro-2}$_{2,3}$, for details we refer to \cite{EidelFischer2018}. Similarly, $\langle f_i \rangle$ is the volume average of body forces in \eqref{eq:StrongFormMicro-2}$_{1}$.

The solution of the homogenized problem is obtained from the variational form 
\begin{equation}
 B_0 (\bm u^0, \bm v) =  \int_{\mathcal{B}} {\color{black} \bm \sigma^0(\bm u^0)}: \bm \varepsilon(\bm v) \, dV   
                      =  \int_{\mathcal{B}} \bm f \cdot \bm v \, dV \, + \,  \int_{\partial \mathcal{B}_N} \bar{\bm t} \cdot \bm v \,  dA
                         \qquad \forall \, \bm v \in \mathcal V \, ,
 \label{eq:VariationalFormHomogenizedProblem}                         
\end{equation}
which follows from multiplying the strong form \eqref{eq:Homogenized-Strong-Form} by test functions $\bm v$ along with the application of Green's formula. For notational convenience we skip in \eqref{eq:VariationalFormHomogenizedProblem} and in the rest of the paper the averaging symbols $\langle \bullet \rangle$, $\langle \bullet \rangle_{\Gamma}$ for $\bm f$, $\bar{\bm u}$ and $\bar{\bm t}$ but keep in mind that these quantities follow from volume and surface averages, respectively. 
 
%.........................................................................................................
 
We consider the piecewise linear continuous FEM in macro- and microspace, respectively. 

We define a macro finite element space as 
\begin{equation}
 \mathcal{S}^p_{\partial \mathcal{B}_D}(\mathcal{B}, {\mathcal T}_H) = \left\{ \bm u^H \in H^1(\mathcal{B})^d; \bm u^H|_{\partial \mathcal{B}_D} = \bar{\bm u}; \bm u^H|_{K} \in {\mathcal{P}}^{p}(K)^{n_{dim}}, \, \forall \, K \in {\cal T}_{H} \right\} \, ,
 \label{eq:MacroFESpace}
\end{equation}
where ${\mathcal P}^{p}$ is the space of (here: linear, $p=1$, or quadratic, $p=2$) polynomials on the element $K$, 
${\mathcal T}_H$ the (quasi-uniform) triangulation of $\mathcal{B} \, \subset \, \mathbb{R}^{n_{dim}}$. Superscript $H$ denotes the characteristic element size, with $H \gg \epsilon$ for efficiency. The space $\mathcal{S}^{p}_{\partial \mathcal{B}_D}$ is a subspace of $\mathcal{V}$ defined in \eqref{eq:HilbertSpaceV}.
 
For the solution of \eqref{eq:StrongFormMicro-2} in the macrodomain
we use the two-scale FEM framework of the FE-HMM as originally proposed in \cite{E-Engquist-2003}
and analyzed for elliptic PDEs in \cite{E-Ming-Zhang-2005}, and, with the focus on linear elasticity, in 
\cite{Assyr2006}.

The macrosolution of the FE-HMM is given by the following variational form:

Find $\bm u^H \in \mathcal{S}_{\mathcal{B}_D}(\mathcal{B}, \mathcal{T}_H)$ such that
\begin{equation}
\label{eq:VariationalFormulationHMM}
 B_H (\bm u^H, \bm v^H) = \int_{\mathcal{B}} \bm f \cdot \bm v^H  \, dV \, + \,  \int_{\partial \mathcal{B}_N} \bar{\bm t} \cdot \bm v^H \, dA
                          \qquad \forall \bm v^H \in \mathcal S_{\partial \mathcal{B}_D} (\mathcal{B}, \mathcal{T}_H) \, ,
\end{equation}
which reads as a standard finite element formulation.  

\subsection{The modified macro bilinear form of FE-HMM}
\label{subsec:Modified-Macro-Bilinear-Form}

If the homogenized constitutive 
tensor $\mathbb{A}^{0}(\bm x)$ is explicitly known, the bilinear form $ B_H (\bm u^H, \bm v^H)$ can be calculated
using standard numerical quadrature according to \eqref{eq:ModifiedBilinearForm-1}, where $\bm x_{K_l}$ and $\omega_{K_l}$ are the quadrature points and 
quadrature weights, respectively
\begin{eqnarray}
 B_H (\bm u^H, \bm v^H) &=& \sum_{K\in \mathcal T_H} \sum_{l=1}^{N_{qp}} \omega_{K_l} \,
                            {\color{black}\left[ {\color{black} \bm \sigma^0}(\bm u^H(\bm x_{K_l})) : \bm \varepsilon(\bm v^H(\bm x_{K_l})) \right]}
                            \label{eq:ModifiedBilinearForm-1} \\
                        &\approx& \sum_{K\in \mathcal T_H} \sum_{l=1}^{N_{qp}} \omega_{K_l} \, \left[ \dfrac{1}{{\color{black}|K_l|}} 
                            {\color{black} \int_{K_l} {\color{black}\bm \sigma^{\epsilon}}(\bm u^h_{K_l}) : \bm \varepsilon(\bm v^h_{K_l}) \, dV} \right] \, .
                            \label{eq:ModifiedBilinearForm-2}                                
\end{eqnarray}

Since $\mathbb{A}^{0}(\bm x)$ is typically not known for heterogeneous materials, the ansatz of FE-HMM is to approximate 
the virtual work expression at point $\bm x_{K_l}$ in the semidiscrete form \eqref{eq:ModifiedBilinearForm-1} by another bilinear form using the known microheterogeneous elasticity tensor $\mathbb{A}^{\epsilon}$, see \eqref{eq:ModifiedBilinearForm-2}. According to this approximation, the solution $\bm u_{K_l}^h$ is obtained on microsampling domains $K_{\delta_l}=\bm x_{K_l} + \delta \, [-1/2, +1/2]^{n_{dim}}$, $\delta \geq \epsilon$, which are each centered at the quadrature points $\bm x_{K_l}$ of $K$, $l=1, \ldots, N_{qp}$. For a visualization see Fig.~\ref{fig:MicMac-Problem-Meshing-etc}. These microsampling domains with volume $|K_{\delta_l}|$ provide the additive contribution to the stiffness matrix of the macro finite element. In order to avoid too heavy notation we will replace $K_{\delta_l}$ by $K_{l}$.  
\begin{Figure}[htbp]
   \begin{minipage}{16.0cm}  
   \includegraphics[width=15.5cm, angle=0]{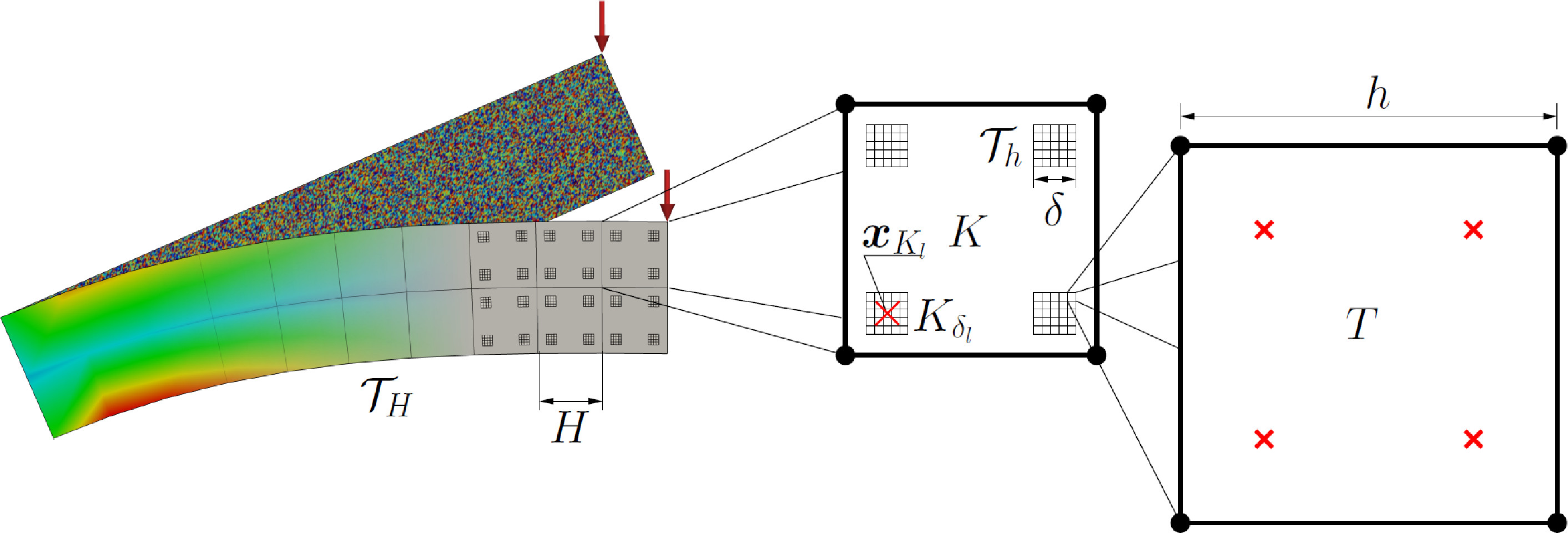}  
\\[2mm]
{\hspace*{3.5cm} (a) \hspace*{4.8cm} (b) \hspace*{3.2cm} (c)}
   \end{minipage}
\caption{FE-HMM as a two-scale finite element method: (a) Macroscopic BVP with macrotriangulation $\mathcal{T}_H$, (b) one macro finite element $K$ of size $H$ with microdomains/RVEs $K_{\delta_l}$ of triangulation $\mathcal{T}_h$, centered at the macro quadrature points $\bm x_{K_l}$, (c) micro finite element $T$ of size $h$ with standard quadrature points. 
\label{fig:MicMac-Problem-Meshing-etc}}
\end{Figure}

The approximation of \eqref{eq:ModifiedBilinearForm-1} by \eqref{eq:ModifiedBilinearForm-2} indicates that FE-HMM crucially relies on a modified quadrature rule and fulfills the equality of the macroenergy density with the microenergy density, thus in agreement with Hill's postulate, \cite{Hill1963}, \cite{Hill1972}.
 
\subsection{Variational formulation of the microproblem}
\label{subsec:Variational-FE-HMM-micro}

It can be shown that the FE-HMM microproblem resembles the discrete version of the cell problem of asymptotic expansion, if it is formulated for each microdomain $K_l$ in $K$ with $l=1, \ldots, N_{qp}$, $K \in \mathcal{T}_H$ like this:
\\[2mm]
Find $\bm u^h_{K_l}$ such that the conditions for macro-micro coupling and for the micro bilinear form \eqref{eq:micro-problem-vers1} are fulfilled:
\begin{equation}
\label{eq:micro-problem-vers1}
\renewcommand{\arraystretch}{1.6}
\left.\begin{array}{rcl}
 \left(\bm u^h_{K_l} -  \bm u^{H}_{lin, K_l} \right) &\in& \mathcal{S}^q (K_l, \mathcal{T}_h)  \\[2mm] 
 B_{K_l}(\bm u^h_{K_l}, \bm w^h_{K_l})  &:=& 
  \displaystyle{\int_{K_l}} {\color{black}\bm \sigma^{\epsilon}}(\bm u^h_{K_l}) : \bm \varepsilon(\bm w^h_{K_l}) \, dV = 0  \\ 
 & &  \forall  \, \bm w^h_{K_l} \in \mathcal{S}^q (K_l, \mathcal{T}_h)   \, ,  
\end{array} \quad \right\} \; 
\end{equation} 
 
where the micro finite element space $\mathcal{S}^q_{}(K_l, \mathcal{T}_h)$ is defined by
\begin{equation}
   \mathcal{S}^q (K_l, \mathcal{T}_h) = \{ \bm w^h \in \mathcal{W}(K_l); \bm w^h \in (\mathcal{P}^q(T))^{n_{dim}}, \, T\in \,\mathcal{T}_h \}\, .
   \label{eq:Periodic-micro-FE-space}
\end{equation}
In \eqref{eq:Periodic-micro-FE-space} $\mathcal{T}_h$ is a quasi-uniform discretization of the sampling domain $K_l$ with mesh size $h \ll \varepsilon$ resolving the finescale and $\mathcal{P}^q$ is the space of polynomials on the element $T$. In the present work we consider linear and quadratic shape functions, $q=1,2$. The particular choice of the Sobolev space $\mathcal{W}(K_l)$ sets the boundary conditions for the micro problems, cf. \cite{Assyr2009}, Sec. 3.2. Among the coupling conditions that fulfill Hill's postulate we consider (i) periodic BC (PBCs), (ii) kinematically uniform displacement conditions (KUBC), and (iii) constant traction conditions (TBC).
 
The linearization of $\bm u^H$ in \eqref{eq:micro-problem-vers1}$_1$ is carried out at the quadrature point $\bm x_{K_l}$  
\begin{equation}
 \bm u^{H}_{lin, K_l} = \bm u^H (\bm x_{K_l}) + (\bm x - \bm x_{K_l}) \cdot \nabla \bm u^H(\bm x_{K_l}) \, .
 \label{eq:linearization_uH}
\end{equation}
It ensures a homogeneous deformation on the microdomain and resembles therein the unit cell problem of asymptotic homogenization
(the FE-HMM perspective) and thus is in the frame of strain-driven first order computational homogenization. 
 
%-----------------------------------------------------------------------------------------------------------------------------------------------------------------------------

For the solution of \eqref{eq:micro-problem-vers1} a basis $\{N_I^H\}_{I=1}^{M_{mac}}$ for the macro finite element space $\mathcal{S}^p_0(\mathcal{B}, \mathcal{T}_H)$ is employed in order to represent the macrosolution $\bm u^H$ of \eqref{eq:VariationalFormulationHMM}. Similarly, a basis $\{N_i^h\}_{i=1}^{M_{mic}}$ of the micro finite element space $\mathcal{S}^q_0(K_l, \mathcal{T}_h)$, \eqref{eq:Periodic-micro-FE-space}, is introduced in order to represent the solution $\bm u^h$ of a microproblem. $M_{mac}$ denotes the number of nodes of the macrodomain, and $M_{mic}$ denotes the number of nodes of each microdomain. 
Hence, the macro- and the microsolution follow the representation
\begin{equation}
\bm u^H = \sum_{I=1}^{M_{mac}} N_I^H \, \bm d_I^H\, , \qquad \bm u^h = \sum_{i=1}^{M_{mic}} N_i^h \, \bm d_i^h\, ,
\label{eq:micro-displacement-vector-in-the-fe-basis-shortened}
\end{equation}
where $\bm d_I^H$ is the displacement vector of macronode $I$, and $\bm d_i^h$ is the displacement vector for micronode $i$.

%-----------------------------------------------------------------------------------------------------------------------------------------------------------------------------

\subsection{Macrostiffness calculation} 
\label{subsec:bottom-up-macrostiffness-calculation}

The macro bilinear form $B^e_H(\bm u^H, \bm v^H)$ is the virtual internal work for a macro finite element. The corresponding bilinear form in terms of the shape functions $B^e_H(\bm N_I^H, \bm N_J^H)$ results in the macro element stiffness matrix contribution $\bm k^{e,mac}_{IJ}$ for macronodes $I, J$, a $d \times d$ matrix. It holds  
\begin{equation}
\bm k^{e,mac}_{IJ} = B_H^{e} (\bm N_I^H, \bm N_J^H) = \sum_{l=1}^{N_{qp}} \dfrac{\omega_{K_l}}{|K_l|} 
                                        \int_{K_l} (\bm L \bm u^{h(I)}_{K_l})^T \mathbb{A}^{\epsilon} (\bm x) \, \bm L \bm u^{h(J)}_{K_l} \, dV \, .
                                        \label{eq:ModifiedBilinearForm-2-for-varphiH}    
\end{equation}
In \eqref{eq:ModifiedBilinearForm-2-for-varphiH} $\bm u^{h(I)}_{K_l}$ is the counterpart of $\bm u_{K_l}^h$ in \eqref{eq:micro-problem-vers1}. It is the dimensionless solution of the microproblem on $K_l$, which is driven by the shape function $N_I^H$ at macronode $I$. In the following, we add $x_i, i=1, \ldots, n_{dim}$ to account for the vector-valued field problem of dimension $n_{dim}$. Consequently, $\bm u^{h(I,x_i)}_{K_l}$ is the microsolution driven by a macroelement unit-displacement state $\bm u^{H(I,x_i)}_{lin, K_l}$ at node $I$ in $x_i$-direction.

For stiffness calculation, problem \eqref{eq:micro-problem-vers1} is reformulated in that $\bm u^{h(I,x_i)}_{K_l}$ replaces
$\bm u^{h}_{K_l}$.

For the coupling of ${\bm u}^{H(I,x_i)}_{lin, K_l}$ with $\bm u^{h(I,x_i)}_{K_l}$ the two fields are expanded into the same basis $\{N_i^h\}_{i=1}^{M_{mic}}$ of $\mathcal{S}^q(K_l, \mathcal{T}_h)$,

\begin{equation}
   {\bm u}^{H(I,x_i)}_{lin, K_l} = \sum_{m=1}^{M_{mic}} \, N^h_{m, K_l} {\bm d}^{H(I,x_i)}_{m} \, , 
   \qquad 
   \bm u^{h(I,x_i)}_{K_l} = \sum_{m=1}^{M_{mic}} N^h_{m, K_l} \, \bm d^{h(I,x_i)}_{m} \, .
   \label{eq:varphi-h-by-psi-shortened}  
\end{equation}

%---------------------------------------------------------------------------------------------------------
The solution of the microproblems for the minimizers $\bm d^{h(I,x_i)}$ is presented in Sec.~\ref{subsec:Solution-of-microproblems}. The macroelement stiffness matrix according to \eqref{eq:ModifiedBilinearForm-2-for-varphiH} yields after some algebra
\begin{eqnarray}          
\bm k^{e,mac}_{IJ} 
         &=& B^e_H\left[\bm N_I^H, \bm N_J^H\right]  \nonumber\\
         &=& \label{eq:k-mac-element-5}
           \sum_{l=1}^{N_{qp}} \dfrac{\omega_{K_l}}{|K_l|} 
           \, \left( \bm d^{h(I)} \right)^T \, \bm K^{mic}_{K_l} \,  \bm d^{h(J)}  \, ,   
\end{eqnarray}          
where $\bm d^{h(I)}=\left(\,\bm d^{h(I,x_1)} | \bm d^{h(I,x_2)} | \bm d^{h(I,x_3)} \, \right)$ for $n_{dim}=3$. A detailed derivation of 
\eqref{eq:k-mac-element-5} is presented in the appendix, Sec.~\ref{sect:Miscellaneous}.

The assembly of $\bm k^{e,mac}_{IJ}$ results in $\bm k^{e,mac}$ and implies a column-wise assembly of $\bm d^{h(I)}$ for $I=1,\ldots, N_{node}$ that results in the transformation matrix $\bm T_{K_l}$

\begin{eqnarray}
 \bm k^{e,mac}_{K} &=& \sum_{l=1}^{N_{qp}} \dfrac{\omega_{K_l}}{|K_l|} 
               \, \, \bm T^{T}_{K_l}\, \bm K^{mic}_{K_l} \, \bm T_{K_l} \label{eq:k-mac-element-6}   \\
 \mbox{with} \quad \bm T_{K_l} &=& \bigg[ \Big[ \big[ \bm d^{h(I,x_i)} \big]_{i=1,\ldots,n_{dim}} \Big]_{I=1, \ldots, N_{node}} \bigg] \label{eq:k-mac-element-7} \,.        
\end{eqnarray}
 
The matrix dimensions
\begin{eqnarray}
  \bm T^ T_{K_l} \in \mathbb{R}^{(N_{node} \cdot n_{dim}) \times (M_{mic} \cdot n_{dim})} ,  
  &\bm K^{mic}_{K_l}  \in \mathbb{R}^{(M_{mic}\cdot n_{dim}) \times (M_{mic}\cdot n_{dim})},&   
  \bm T_{K_l} \in \mathbb{R}^{(M_{mic} \cdot n_{dim}) \times (N_{node} \cdot n_{dim})} \nonumber  \\ 
  &\bm k^{e,mac}_K \in \mathbb{R}^{(N_{node} \cdot n_{dim})\times (N_{node} \cdot n_{dim})}& \nonumber 
\end{eqnarray}
underpin that $\bm T_{K_l}$ is not only the agency of a micro-to-macro stiffness transfer but also a compression operator that transforms $\bm K^{mic}_{K}$ into $\bm k^{e,mac}_{K}$.  
 
\begin{Figure}[htbp]
   \begin{minipage}{16.0cm}  
    \begin{center} 
    	\psfrag{x1}{$x_1$}
    	\psfrag{x2}{$x_2$}
    	\psfrag{\&}[c][Bl]{\&}
          \includegraphics[height=3.2cm, angle=0]{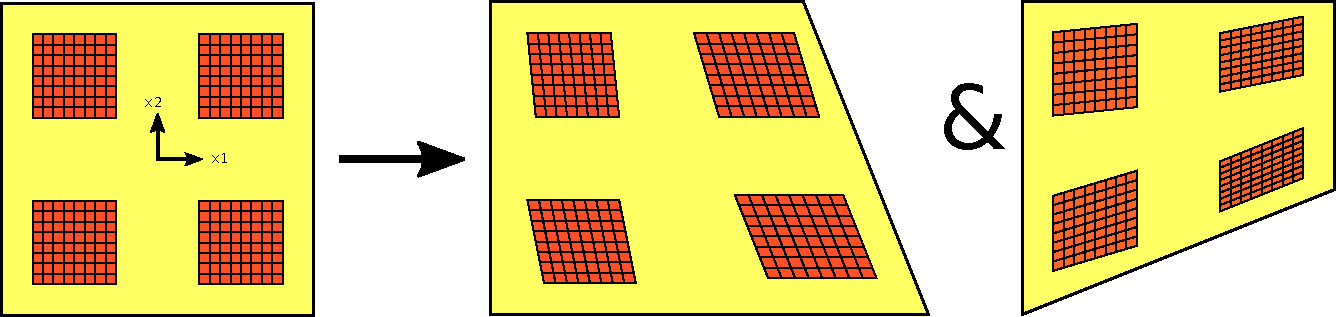}
    \end{center}
   \end{minipage} 
\caption{\color{black} Unit displacement states in $x_i$-directions, $i=1, 2$ applied to the lower right macro element node $I$ and the uniform deformations of the microdomains/RVEs.
 \label{fig:Determination-of-varphi-H}}
\end{Figure}

In the context of stiffness computation, a macro element shape function represents a unit displacement state for macro node  $I, I=1, \ldots, N_{node}$ in each direction of space $x_i\, |\, i=1, \ldots, n_{dim}$.
They drive the microproblem in terms of the corresponding nodal values ${\bm d}^{H(I,x_i)}_{m}, m=1, \ldots, M_{mic}$ in each microdomain to evaluate the macroelement stiffness $\bm k^{e,mac}_{IJ}$. Each unit displacement state in $x_i$-direction induces in ${\bm d}^{H(I,x_i)}$ nonzero components only in $x_i$, for example
${\bm d}^{H(I,x_i)}|_{i=2} = \left[ 0, {d}^{H(I,x_2)}_{1,x_2}, 0, \hdots, 0, {d}^{H(I,x_2)}_{M_{mic},x_2}, 0\right]^T$.
  
Figure \ref{fig:Determination-of-varphi-H} visualizes two unit displacement states applied to the lower right macronode in two directions along with the corresponding uniform RVE-deformations that follow from linearization according to \eqref{eq:linearization_uH}.  
  
%-----------------------------------------------------------------------------------------------------------------------------------------------------------------------------
\newcommand\myeq{\mathrel{\overset{\makebox[0pt]{\mbox{\normalfont\sffamily !}}}{=}}}

\subsection{Solution of the microproblems} 
\label{subsec:Solution-of-microproblems}

The microstiffness matrix is obtained by Gauss-Legendre numerical quadrature on the element level, the assembly of the element stiffness matrices results in the total stiffness matrix for an RVE.
 
With the RVE microstiffness matrix in hand, the microproblem can be solved. Here the method of Lagrange multipliers is chosen such that a saddlepoint problem must be solved. The total energy for a macro unit displacement state reads
\begin{eqnarray}
 \mathcal{L} (\bm d^{h(I,x_i)}, \bm \lambda^{(I,x_i)}) &=& \dfrac{1}{2} \left(\bm d^{h(I,x_i)}\right)^T \bm{K}_{K_l}^{mic} \, \bm d^{h(I,x_i)}
 + \bm \lambda^{(I,x_i)\,T} \, \bm G \, \left(\bm d^{h(I,x_i)} - \overline{\bm d}^{H(I,x_i)}\right) 
 \label{eq:Lagrange-functional} \\
 & & \mbox{for }I=1, \ldots, N_{node}, \, \, \mbox{and} \, \, x_i\,|\,i=1, \ldots, n_{dim}\, ,\nonumber
\end{eqnarray}
where $\bm G$ contains the kinematical coupling constraints. Details of various coupling conditions and their implementation in the frame of Lagrange multipliers are described in Sec.~\ref{sec:Implementation-CouplingConditions}. In either case the Lagrange multipliers represent external forces which enforce the micro coupling condition. 

The vector of Lagrange multipliers $\bm \lambda^{(I,x_i)} \in \mathbb R^{(1+L)\cdot n_{dim}}$, where $L$ depends on the type of microcoupling, reads for $n_{dim}=3$ as
\begin{equation}
\bm \lambda^{(I,x_i)}=\{\lambda^{(I,x_i)}_{0,x_1}, \lambda^{(I,x_i)}_{0,x_2}, \lambda^{(I,x_i)}_{0,x_3}, \lambda^{(I,x_i)}_{1,x_1}, \lambda^{(I,x_i)}_{1,x_2}, \lambda^{(I,x_i)}_{1,x_3}, 
\ldots ,\lambda^{(I,x_i)}_{L,x_1}, \lambda^{(I,x_i)}_{L,x_2}, \lambda^{(I,x_i)}_{L,x_3}\}^T \, . 
\end{equation}
The first variation of $\mathcal{L}$ with respect to $\bm d^{h(I,x_i)}$ and $\bm \lambda^{(I,x_i)}$ results in the  stationarity conditions
\begin{equation}
\left[ \begin{array}{cc}
 \bm K^{mic}_{K_l} & \bm G^T \\
 \bm G  & \bm 0
\end{array} \right] 
\left[ \begin{array}{c}
\bm d^{h(I,x_i)} \\
\bm \lambda^{(I,x_i)}\\
\end{array}\right]
=
\left[ \begin{array}{c}
\bm 0   \\
\bm G \, {\bm d}^{H(I,x_i)}    \\
\end{array}\right]  \, \, \mbox{for} \, \, I=1, \ldots, N_{node}, \, \, i=1, \ldots, n_{dim} \, .
\label{Solve4alpha-lambda}
\end{equation}
The solution of \eqref{Solve4alpha-lambda} subject to $N_{node} \cdot n_{dim}$ right hand sides can be carried out efficiently since the coefficient matrix in \eqref{Solve4alpha-lambda} is constant. The solution vectors are augmented to full matrices, hence, $\bm d^{h(I,x_i)} \rightarrow \bm T$, $\bm \lambda^{(I,x_i)} \rightarrow \bm \Lambda$, ${\bm d}^{H(I,x_i)} \rightarrow \bm d^H$.   

The solution of \eqref{Solve4alpha-lambda} serves the purpose to compute the transformation matrix $\bm T_{K_l}$ according to \eqref{eq:k-mac-element-7}.  After the consecutive solution of the global macroproblem for $\bm u^H$, the microproblems have to be solved. Then, \eqref{Solve4alpha-lambda} is driven by the true macroscopic displacement vector $\bm d^H$, which results in the true microdisplacements $\bm d^{h}$.

\section{The coupling conditions}
\label{sec:Implementation-CouplingConditions}

This section gives a brief account of the implementation of the coupling conditions (Dirichlet, periodic and Neumann) in a Lagrange-Multiplier framework. Doing so, the particular format of the constraint matrix $\bm G$ is detailed.
  
 \subsection{Dirichlet coupling}
 \label{subsec:Implementation-DirichletCoupling}
 
 The simplest coupling condition is the Dirichlet coupling, where linear displacements following from the macro displacement field are applied to the boundaries of a microdomain
 
 \begin{equation}
 \bm u^{h(I,x_i)}_{K_l} = \bm u^{H(I,x_i)}_{lin,K_l} \quad \text{on} \ \partial K_l \, .
 \end{equation}
 
 In this case the constraint matrix $\bm G$ contains $L \cdot n_{dim}$ rows with $L$ the number of boundary nodes and $n_{dim}$ the number of degrees of freedom per node. Each row contains a 1 pointing at a degree of freedom of a boundary node and 0 elsewhere. By doing so \eqref{Solve4alpha-lambda} directly couples the nodal micro displacements $\bm d^{h(I,x_i)}$ on the RVE boundary to the nodal micro displacements following from the macro displacement field ${\bm d}^{H(I,x_i)}$. 
 
 The expression on the right hand side of the system of equations $\bm G \, {\bm d}^{H(I,x_i)}$ can be derived by inserting the micro coordinates of the boundary nodes in the reference configuration into the linearized macro shape functions.  
 
 \subsection{Periodic coupling}
 \label{subsec:periodic-coupling}
 
 Periodic coupling conditions imply periodic displacements and anti-periodic tractions on the boundaries of the microdomain. It holds
 \begin{align}
 \left( \bm u^{h,(I,x_i)}_{K_l} - \bm u^{H(I,x_i)}_{lin,K_l} \right)^+ &= \left( \bm u^{h,(I,x_i)}_{K_l} - \bm u^{H(I,x_i)}_{lin,K_l} \right)^-   \\
 \bm t^+ &= - \bm t^- \, .
 \end{align}
 
 The boundary of the microdomain is here split up in a part $\partial K_l^+$ and a part $\partial K_l^-$ such that $\partial K_l^+ \cup \partial K_l^- = \partial K_l$ having opposite outward normal vectors $\bm n^+ = - \bm n^-$.
 
 The constraint matrix $\bm G$ then contains $L \cdot n_{dim}$ rows with $L$ the number of non-redundant periodic couples and $n_{dim}$ the number of degrees of freedom per node. The single rows of the constraint matrix contain a 1 pointing at a degree of freedom of a node on $\partial K_l^+$ and a -1 pointing at the corresponding degree of freedom of the node on $\partial K_l^-$, the other entries of $\bm G$ vanish.
 
 Periodic displacement fluctuations only eliminate the rotational rigid body motions, the rigid body translations however are not discarded by PBC. For that reason they must be eliminated by an additional condition  
 \begin{equation}
  \int_{K_l} \left(\bm u^{h(I,x_i)}_{K_l} -  {\bm u}^{H(I,x_i)}_{lin, K_l} \right) \, dV = \bm c \, .
  \label{eq:normalization-condition-discrete-version} 
\end{equation}
Since \eqref{eq:normalization-condition-discrete-version} is a normalization condition for the periodic fluctuations, the particular choice of the constant is inconsequential for the microsolution. Here we choose $\bm c = \bm 0$.
 
By multiplying the transpose of $\bm G$ with the Lagrange multipliers, the anti-periodic tractions for each couple of periodic nodes are realized.
 
For $n_{dim}=3$ matrix $\bm G$ exhibits the format 
\begin{equation}
\bm G = 
\left[ \begin{array}{ccccccc}
 b_1 & 0   & 0   & \, \ldots \ldots \, & b_{M_{mic}} & 0 & 0   \\
 0   & b_1 & 0   & \, \ldots \ldots \, & 0 & b_{M_{mic}} & 0   \\
 0   & 0   & b_1 & \, \ldots \ldots \, & 0 & 0 & b_{M_{mic}}   \\[2mm] 
 \multicolumn{7}{c}{\overline{\bm G}}                          \\
\end{array} \right] \, .
\label{eq:D-matrix-3d}
\end{equation}
In $\bm G$, the first $d$ rows contain the normalization condition for the fluctuations, to which the first $d$ Lagrange multipliers are associated. The coefficients $b_i$, $i=1,\ldots, M_{mic}$ in \eqref{eq:D-matrix-3d} follow from the coupling condition \eqref{eq:normalization-condition-discrete-version}. It holds  
\begin{eqnarray}
   \int_{K_{\delta}} \bm u^{h(I,x_i)} \, d V &=& \sum_{T \in \mathcal{T}_h} \sum_{m=1}^{M_{mic}} \bm d_m^{h(I,x_i)}
                       \underbrace{\int_{T} \bm N^h_m \, dV}_{=:b_m} 
%                      =  \sum_{T \in \mathcal{T}_h} \sum_{m=1}^{n_{node}} \alpha_m \, b_m
                       =  \sum_{T \in \mathcal{T}_h} \bm d^{h(I,x_i)} \cdot \bm b \\
\int_{K_{\delta}} {\bm u}^{H(I,x_i)}_{lin} \, dV &=& \sum_{T \in \mathcal{T}_h} \sum_{m=1}^{M_{mic}} \bm d_m^{H(I,x_i)}
                             \int_{T} \bm N^h_m \, dV 
%                      =  \sum_{T \in \mathcal{T}_h} \sum_{m=1}^{n_{node}} \alpha_m \, b_m
                       =  \sum_{T \in \mathcal{T}_h} {\bm d}^{H(I,x_i)}  \cdot \bm b \\                       
%\myeq 0\, \\
\mbox{where} \quad  b_m &=& \int_{T} \bm N^h_m \, dV \, .%=\dfrac{1}{n_{node}} \, |T| 
\end{eqnarray}

\subsection{Neumann coupling} 
\label{subsec:Implementation-NeumannCoupling}
 
The Neumann coupling condition of constant tractions on the element boundary reads as
 
\begin{equation}
 \bm{t}(\bm x) \ = \ \boldsymbol{\sigma}_{K_l}(\bm u^{H(I,x_i)}) \, \bm n (\bm x) \quad \text{on} \ \partial K_l \, .
 \label{eq:neumann_coupling_traction}
\end{equation}
 
Since stress of the macro Gauss point $\boldsymbol{\sigma}_{K_l}(\bm u^{H(I,x_i)})$ for macro unit displacement states is not known, the condition is reformulated to a weak constraint in terms of a macroscopic strain $\boldsymbol{\varepsilon}_{K_l}(\bm u^{H(I,x_i)})$ cf. \cite{Miehe-Koch-2002} 
 
\begin{equation}
 \dfrac{1}{\vert K_l \vert} \int_{\partial K_l}^{} \text{sym} [\bm{u^{h(I,x_i)}}(\bm x) \otimes \bm n] \,dA \ = \ \boldsymbol{\varepsilon}_{K_l}(\bm u^{H(I,x_i)}) \, .
\end{equation}
 
Introducing a discrete nodal normal vector 
 
\begin{equation}
 \bm n_q := \frac{1}{2}[\bm x_{q+1}- \bm x_{q-1}] \times \bm e_3
\end{equation}
 
with $\bm x_{q-1}$ and $\bm x_{q+1}$ being the neighbor nodes of node $q$ on the boundary, we get 
 
\begin{equation}
 \dfrac{1}{\vert K_l \vert} \sum_{q=1}^{L} \text{sym} [\bm d^{h(I,x_i)}_q \otimes \bm n_q] \ = \  \boldsymbol{\varepsilon}_{K_l}(\bm u^{H(I,x_i)}) \, , 
 \label{eq:neumann_cc_int}
\end{equation}
 
where $L$ is the number of nodes on the boundary of the microdomain. The neighboring nodes $\bm x_{q-1}$ and $\bm x_{q+1}$ have to be oriented so that $\bm n_q$ is an outward normal vector. 
 
Expression \eqref{eq:neumann_cc_int} can be recast into a matrix representation
 
\begin{equation}
 \sum_{q=1}^{M} \bm{G}_q \bm{d}^{h(I,x_i)}_q \ = \  \boldsymbol{\varepsilon}_{K_l}(\bm u^{H(I,x_i)}) \, ,
\end{equation}
 
where $\bm{G}_q$ is depending on the normal vector $\bm n_q$ and reads as 
 
\begin{equation}
 \bm G_q \ := \ \dfrac{1}{\vert K_l \vert} \left[ \begin{array}{c c}
 2n_1 & 0 \\
 0 & 2n_2 \\
 n_2 & n_1
 \end{array} \right]_q \, ,
\end{equation}
 
for $n_{dim}=2$. The constraint matrix $\bm G$ follows from assembling the single $\bm G_q$ matrices
 
 \begin{equation}
 \bm{G} \bm{d}^{h(I,x_i)} = \boldsymbol{\varepsilon}_{K_l}(\bm u^{H(I,x_i)}) \quad \text{on} \quad \partial K_l \, .
 \label{eq:neumann_Gd_epsilon}
 \end{equation}
 
 The term on the right hand side of \eqref{Solve4alpha-lambda} can also be replaced by the strains $\boldsymbol{\varepsilon}_{K_l}(\bm u^{H(I,x_i)})$.
 
Again, the constant traction BC alone does not eliminate the rigid body motions of the RVE and the corresponding zero eigenvalues of the stiffness matrix. The rigid body motions must be eliminated by additional kinematical constraints (3 for $n_{dim}=2$ and 6 for $n_{dim}=3$), which can be realized e.g. by a so-called semi-Dirichlet coupling introduced in  \cite{JaviliSaebSteinmann2017}.   

%\begin{Figure}[htbp]
% \begin{minipage}{16.5cm}  
%      \centering
%      \subfigure[]{\psfrag{A}[c][tc]{A}
% 	\psfrag{B}[c][bl]{B}
% 	\psfrag{tAx}[t][c]{$t^\text{A}_{x_1}$}
% 	\psfrag{tay}[c][Br]{$t^\text{A}_{x_2}$}
% 	\psfrag{tB}[cl][bl]{$t^\text{B}_{x_2}$}
% 	\psfrag{bo}[lc][lb]{$\partial K_{l, \, 0}$}
% 	\psfrag{b}[lc][rb]{$\partial K_l, \bm \zeta \ne \bm 0$}
% 	\psfrag{bto}[r][r]{$\partial K_l, \bm \zeta = \bm 0$}
% 	\includegraphics[width=0.42\linewidth]{semi_dirichlet}}
%        \hspace*{11mm}
%        \subfigure[]{\psfrag{A}[c][tc]{A}
% 	\psfrag{B}[l][l]{B}
% 	\psfrag{bo}[cl][Bl]{$\partial K_{l, \, 0}$}
% 	\psfrag{bto}[r][r]{$\partial K_l, \bm \zeta = \bm 0$}
% 	\psfrag{eta}{$\eta$}
% 	\includegraphics[width=0.38\linewidth]{semi_dirichlet_eta}}
%  \end{minipage}  
%        \caption{{\bf Semi-Dirichlet BC for Neumann coupling:}(a) Elimination of rigid body motions in $n_{dim}=2$ by adding three Dirichlet BC at single nodes, which renders the system statically determined. Undeformed microdomain $\partial K_{l, \, 0}$ fixed at node A and node B and deformed microdomain $\partial K_l$ with and without the additional displacement constraints, (b) variation of $\eta$ to realize the constant traction condition without spurious forces in A and B.}
% 	\label{fig:semi_dirichlet_total}
%\end{Figure}

\begin{Figure}[htbp] 
 \begin{minipage}{16.5cm}   
      \centering
 	\includegraphics[width=0.82\linewidth]{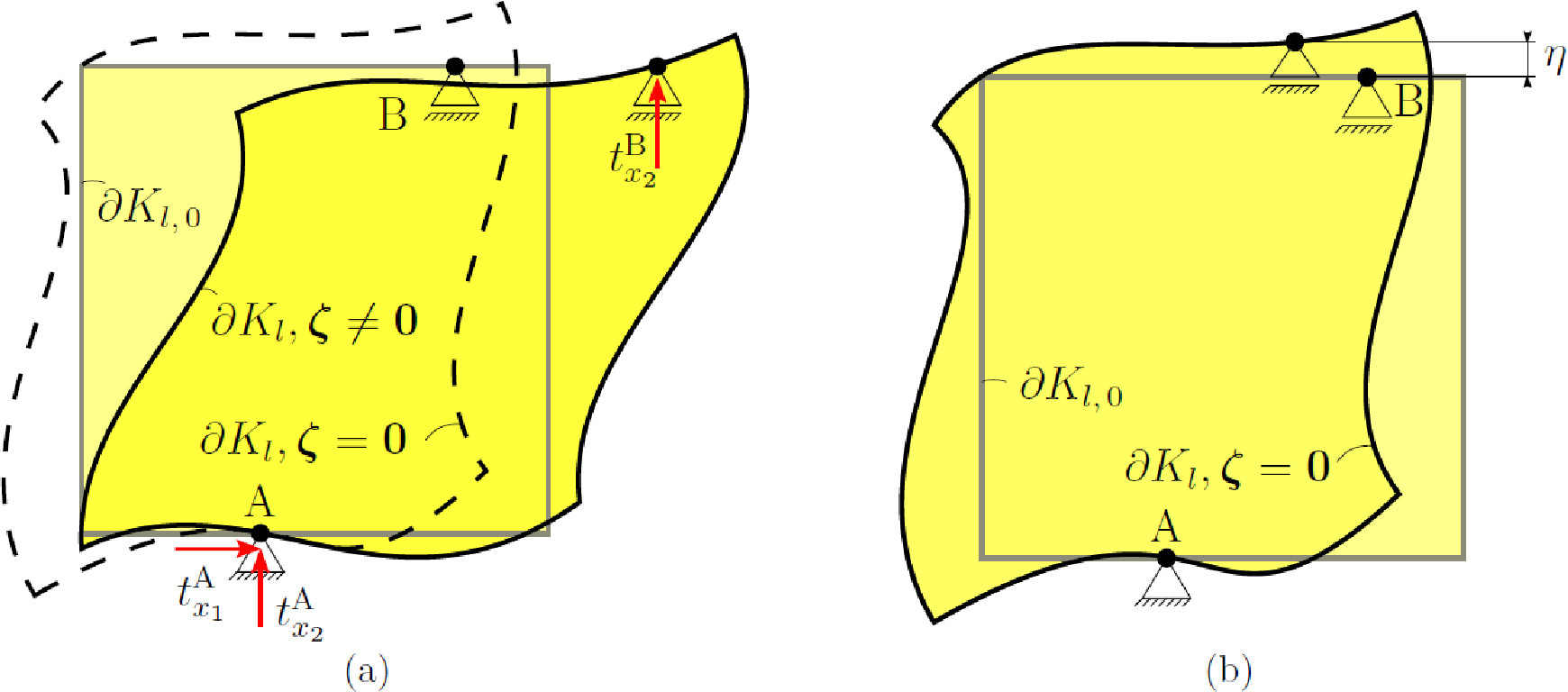} 
  \end{minipage}   
        \caption{{\bf Semi-Dirichlet BC for Neumann coupling:}(a) Elimination of rigid body motions in $n_{dim}=2$ by adding three Dirichlet BC at single nodes, which renders the system statically determined. Undeformed microdomain $\partial K_{l, \, 0}$ fixed at node A and node B and deformed microdomain $\partial K_l$ with and without the additional displacement constraints, (b) variation of $\eta$ to realize the constant traction condition without spurious forces in A and B.} 
 	\label{fig:semi_dirichlet_total} 
\end{Figure}

In Fig.~\ref{fig:semi_dirichlet_total} (a) the additional Dirichlet BC in points A and B and their influence on the {\color{black}reaction forces} on the boundary are shown. The {\color{black}reaction forces} in A and B can be given in terms of the stress $\boldsymbol{\sigma}_{K_l}$ in the corresponding macro quadrature point 
 
 \begin{align}
 \bm t^A &= \boldsymbol{\sigma}_{K_l} \cdot \bm n^A + \boldsymbol{\zeta}^A \, , \\
 \bm t^B &= \boldsymbol{\sigma}_{K_l} \cdot \bm n^B + \boldsymbol{\zeta}^B \, .
 \end{align}
 
 To enforce the Dirichlet BC in A and B, an additional force is needed which influences the {\color{black}reaction forces} in A and B. In order to satisfy \eqref{eq:neumann_coupling_traction}, the Dirichlet BC have to be chosen such that $\boldsymbol{\zeta}^A = \boldsymbol{\zeta}^B = \bm 0$.
 
 To do so the Dirichlet condition in point B, where (here) the node is merely fixed in $x_2$ direction, is modified by moving point B from its former position about $\eta$ in $x_2$ direction. This has to be done until \eqref{eq:neumann_Gd_epsilon} and $\boldsymbol{\zeta}^A = \boldsymbol{\zeta}^B = \bm 0$ are satisfied. Writing down these conditions in a residual vector $ \bm R $ leads to
 
 \begin{equation}
 \bm R (\boldsymbol{\sigma}_{K_l}, \eta) = \left[ \bm G \bm d^{h(I,x_i)} - \boldsymbol{\varepsilon}_{K_l}^{(I,x_i)}, \boldsymbol{\zeta}^{(B(I,x_i))} \right] \stackrel{!}{=} \bm 0 \, , 
 \label{eq:neumann_coupling_error_vector}
 \end{equation}
 
which is solved using the Newton-Raphson scheme. Therefore \eqref{eq:neumann_coupling_error_vector} has to be linearized which results in 
 
 \begin{equation}
 \bm R (\boldsymbol{\sigma}_{K_l , \, i+1}, \eta_{i+1}) = \bm R (\boldsymbol{\sigma}_{K_l , \, i}, \eta_{i}) + \left. \dfrac{\partial \bm R}{\partial \boldsymbol{\sigma}_{K_l}} \right|_i : \Delta \boldsymbol{\sigma}_{K_l, \, i} + \left. \dfrac{\partial \bm R}{\partial \eta} \right|_i  \, \Delta \eta_i \stackrel{!}{=} \bm 0 \, .
 \end{equation}
 
The Lagrange multipliers which follow from the Neumann coupling condition contain the macroscopic stresses in the corresponding quadrature point. The Lagrange multipliers which follow from the semi-Dirichlet coupling contain the additional forces required to enforce the semi-Dirichlet constraints, which must vanish.

In contrast to the iterative, hence expensive solution using the semi-Dirichlet coupling where the micro system of equations has to be solved at least twice, the approach of \cite{Miehe-Koch-2002} enforces regularity of the microproblem by adding a perturbation to the entries on the diagonal of the micro stiffness matrix. Section \ref{sec:NumericalExamples} will provide a quantitative comparison of the methods. 

\subsection{Numerical implications of different coupling conditions}
 
The above described coupling conditions all lead to the same system of equations \eqref{Solve4alpha-lambda}

which has to be solved for the micro displacements. The numerical effort depends on the size of the system of equations and, in the case of Neumann coupling with semi-Dirichlet coupling, additionally on the number of required iterations. 
 
While the micro stiffness matrix $\bm K^{mic}_{K_{l}}$ exhibits the same format for all of the described coupling conditions, the constraints matrix $\bm G$ does not. For Dirichlet coupling $\bm G$ has $L \cdot n_{dim}$ rows with $L$ the number of boundary nodes and $n_{dim}$ the number of degrees of freedom per node. In case of periodic coupling $\bm G$ has $L \cdot n_{dim}$ rows with $L$ the number of non-redundant periodic couples. For Neumann coupling the number of rows equals the number of strain components $\bm \varepsilon$ and the number of additional semi-Dirichlet coupling conditions.
 
Especially for fine micro discretizations with many boundary nodes the system of equations for Neumann coupling will be smaller than for Dirichlet and periodic coupling. 
 
For Dirichlet and periodic coupling the system of equations has to be solved only once, for Neumann coupling realized by the semi-Dirichlet approach the set of equations has to be solved in each of the iterations.

{\color{black}

It should be mentioned that at least for Dirichlet and periodic coupling conditions the method of Lagrange multipliers can be realized by the more efficient direct use of the macroscopic displacement field.

\begin{Table}[htbp]
	\begin{minipage}{16.5cm}  
		\centering
		\renewcommand{\arraystretch}{1.2} 
		\begin{tabular}{c c c}
			\hline
			 Coupling & Direct solution method & Lagrange multiplier method\\
			\hline 
			Dirichlet  & $2 \cdot (N-2)^2$ & $2 \cdot \left(N^2 + 4(N-1) \right)$ \\
			\hline  
			Periodic & $2 \cdot \left((N-2)^2 + (N-1) + (N-2) \right)$ & $2 \cdot \left(N^2 + N + (N-1)\right)$    \\
			\hline                 
		\end{tabular}
		\newline 
	\end{minipage}
	\caption{{{\color{black}{\bf Comparison of the degrees of freedom in 2D:} Direct solution method versus Lagrange multiplier method for Dirichlet and periodic coupling.}}
	\label{tab:Dof_direct_solving_vs_Lagrange}}
\end{Table}

For a uniform micro mesh in 2D with $N$ nodes per edge Tab.~\ref{tab:Dof_direct_solving_vs_Lagrange} displays the number of degrees of freedom for both methods in the cases of Dirichlet and periodic coupling conditions. Especially for Dirichlet coupling the direct solution of the microproblem is a convenient option. The number of degrees of freedom for the direct implementation of Dirichlet coupling is reduced to those of the $(N-2)^2$ nodes in the interior of the microdomain, whereas for the Lagrange multiplier method not only the degrees of freedom of the $N^2$ nodes have to be considered, but additionally the degrees of freedom to impose the coupling conditions on the $4 (N-1)$ boundary nodes. Furthermore, Dirichlet coupling conditions can be easily realized by means of static condensation.
}

%---------------------------------------------------------------------------------------------------------
\section{A priori error estimates and a posteriori error estimation}
%---------------------------------------------------------------------------------------------------------
\label{sec:Apriori-and-Aposteriori-ErrorEstimates}

This section~\ref{sec:Apriori-and-Aposteriori-ErrorEstimates} provides the unified a priori estimates covering the macro error, the micro error and the modeling error. Moreover, the {\color{black}recovery}-type error estimator of Zienkiewicz-Zhu based on superconvergent stress and strain is introduced and contrasted to error computation based on reference solutions.

\subsection{A priori estimates} 
\label{subsec:A-priori-estimates}

FE-HMM as a particular instance of the most general Heterogeneous Multiscale Method HMM \cite{E-Engquist-2003}, \cite{E-Engquist-Huang-2003}, \cite{E-Engquist-Li-Ren-VandenEijnden2007}, \cite{Assyr-etal2012} has its foundation in mathematical homogenization by asymptotic expansion, \cite{Bensoussan-Lions-Papanicolau-BOOK-1976}, \cite{Sanchez-Palencia-BOOK-1980}, \cite{Allaire1992}, \cite{Cioranescu-Donato-BOOK-1999}. 

A priori estimates for various types of partial differential equations (PDEs) have been derived for FE-HMM by virtue of its foundation in mathematical homogenization; for the elliptic case we refer to \cite{E-Ming-Zhang-2005}, \cite{Ohlberger2005}, for the elliptic case of linear elasticity in a geometrical linear setting to \cite{Assyr2006}, \cite{Assyr2009}. 
A posteriori error analysis along with upper and lower bounds of a residual-based error estimator have been presented in \cite{AssyrNonnenmacher2011}, for an adaptive strategy governed by quantities of interest we refer to \cite{AssyrNonnenmacher2013}.

The total FE-HMM error can be decomposed into three parts
\begin{equation}
   || \bm u^0 - \bm u^H || \, \leq \, \underbrace{|| \bm u^0 - \bm u^{0,H} ||}_{\displaystyle e_{mac}} 
                             \, + \, \underbrace{|| \bm u^{0,H} - \widetilde{\bm u}^H ||}_{\displaystyle e_{mod}}
                             \, + \, \underbrace{|| \widetilde{\bm u}^H - \bm u^H ||}_{\displaystyle e_{mic}} \, ,
   \label{eq:Error-decomposition-mac-mod-mic}
\end{equation}
where $e_{mac}$, $e_{mod}$, $e_{mic}$ are the macro error, the modeling error, and the micro error.

Here, $\bm u^0$ is the solution of the homogenized problem \eqref{eq:Homogenized-Strong-Form}, $\bm u^H$ is the FE-HMM solution, $\bm u^{0,H}$ is the standard (single-scale) FEM solution of problem \eqref{eq:VariationalFormHomogenizedProblem} that is obtained through exact $\mathbb{A}^{0}$; and $\widetilde{\bm u}^H$ is the FE-HMM solution obtained through exact microfunctions (in $W(K_l)$).

For sufficiently regular problems the following a priori estimates hold in the $L^2$-norm, the $H^1$-norm and the energy-norm (definition of these norms in Appendix \ref{subsec:Definition-norms}):  
\begin{eqnarray}
   || \bm u^0 - \bm u^H ||_{L^2(\mathcal{B})} &\leq& C\left( H^{p+1} + \left(\dfrac{h}{\epsilon}\right)^{2q} \right) + e_{mod}    \, ,
    \label{eq:Total-Error-estimate-L2} \\
   || \bm u^0 - \bm u^H ||_{H^1(\mathcal{B})} &\leq& C\left( H^p + \left(\dfrac{h}{\epsilon}\right)^{2q} \right) + e_{mod}   \, ,
   \label{eq:Total-Error-estimate-H1} \\
   || \bm u^0 - \bm u^H ||_{A(\mathcal{B})} &\leq& C\left( H^p + \left(\dfrac{h}{\epsilon}\right)^{2q} \right) + e_{mod}   \, .
   \label{eq:Total-Error-estimate-Energy}  
\end{eqnarray}
 
For $e_{mod}$ in \eqref{eq:Total-Error-estimate-L2}--\eqref{eq:Total-Error-estimate-Energy} it holds 

\begin{equation}
\label{eq:ModelingError}
 e_{mod} = \left\{ \begin{array}{ll}
         0 & \mbox{for periodic coupling with} \,\, \delta/\epsilon \in \mathbb{N} \\
         {\color{black}C\,\dfrac{\epsilon}{\delta}}  & \mbox{for Dirichlet coupling with} \,\, \delta > \epsilon.\end{array} \right. \,    
\end{equation} 
{\color{black} given that the hypotheses hold, that the elasticity tensor $\mathbb{A}^{\epsilon}$ is periodic on the RVE and, that the micro solution is sufficiently smooth, \cite{JeckerAbdulle2016}.}

The modeling error for Dirichlet coupling in \eqref{eq:Total-Error-estimate-L2}--\eqref{eq:Total-Error-estimate-Energy} is due to boundary layers \cite{E-Ming-Zhang-2005} (Thm. 1.2), \cite{Assyr2009}. So even for $H \rightarrow 0$ and $h \rightarrow 0$ there is a residual error. 
\\[4mm]
{\bf Remark 1} 
\\[2mm]
(i) Order $2q$ of the micro error in the $L^2$-norm according to \eqref{eq:Total-Error-estimate-L2} seems to contradict standard FEM results. Even more, its order each in the  $H^1$- and energy- norms scales with $2q$ in the same order as in the $L^2$-norm, a phenomenon which is referred to as \emph{superconvergence} in the context of standard (single-scale) finite element methods. For the latter however, superconvergence is not inconditional, since that kind of superconvergence is not only restricted to particular element sites but also to the rectangular shape of them \cite{Barlow1976}. Notice that the latter superconvergence can be used for the construction of a {\color{black}recovery}-type error estimator based on the so-called Superconvergent Patch Recovergy (SPR) introduced by \cite{SPR}, \cite{SPR2}, a concept which is adopted in Sec.~\ref{subsec:superconvergence} of the present work.

(ii) The alleged inconsistency of the micro convergence error order is resolved by the fact that \eqref{eq:Total-Error-estimate-L2}--\eqref{eq:Total-Error-estimate-Energy} describe the micro error as propagated to the macroscale; it is measured by macro quantities, i.e. by $\bm u^H$ in the $L^2$-norm, and by macroscopic stress and strain in the energy-norm. In contrast to this propagated micro error on the macroscale, the micro error on the microscale, which is measured by micro quantities, scales in the order of $\mathcal O((h/\epsilon)^{q+1})$ in the $L^2$-norm and of $\mathcal O((h/\epsilon)^{q})$ in the $H^1$- and energy-norm thus being consistent with estimates of standard finite elements. 
 
(iii) For its composition covering both the macro error as well as the micro error, the estimates \eqref{eq:Total-Error-estimate-L2}--\eqref{eq:Total-Error-estimate-Energy} enable strategies to achieve the optimal convergence order for minimal computational costs in uniform micro-macro discretizations; they answer the practical question on how to improve in two-scale finite element frameworks the accuracy by $H-/h$-refinements on both the macro- and the microscale most efficiently.

\begin{Table}[htbp]
\begin{minipage}{16.5cm}  
\centering
\renewcommand{\arraystretch}{1.2} 
\begin{tabular}{c|cc}
               \hline 
  macro-,micro-FEM  &  $L^2$-norm                    & $H^1$-/energy-norm \\
               \hline  
  $P^p$, $P^q$      & $N_{mic} = (N_{mac})^{p+1/2q}$ & $N_{mic} = (N_{mac})^{p/2q}$    \\
               \hline                 
\end{tabular}
\newline 
   \end{minipage}
\caption{{\color{black}{\bf} Optimal uniform micro-macro refinement strategies: full order for minimal effort. {\color{black} $N_{mic}$ denotes the number of unknowns on the microscale, $N_{mac}$ on the macroscale.}
\label{tab:Best-mic-mac-RefinementStrategies}}}
\end{Table}

Table \ref{tab:Best-mic-mac-RefinementStrategies} displays the optimal uniform micro-macro refinement strategies for the error in the  $L^2$-norm and the $H^1$-/energy-norm. Of course, the strategy's dependency on the polynomial order of macro shape functions $p$ and $q$ on the microscale crucially relies on sufficient regularity of the corresponding BVPs. 

{\color{black}(iv) For numerical convergence analyses of the macro error it is enough to compute the total error at various macrodiscretizations $H$ keeping the micro errors constant by employing a constant micro discretization $h$. Consequently, the reference solution for error computation is $\bm u^{{\color{black}H},\text{ref}}(H \rightarrow 0, h=const.)$. 
Convergence analyses of the micro error are carried out analogously; in this case the reference solution is $\bm u^{{\color{black}H},\text{ref}}(H=const., h \rightarrow 0)$.}   

{\color{black}(v) To our knowledge no estimate for the modeling error along with Neumann BC is available in mathematical  literature.}

\subsection{Error computation}
\label{subsec:error-computation} 

For {\color{black} the special case of} a micro error {\color{black} convergence} analysis on the microscale (as e.g. on a selected microdomain) any macroscopic influence must be switched off. However, since different micro discretizations imply numerical differences in the stiffness approximation, they influence macroscopic displacements, which themselves influence through the postprocessing the microscopic quantities. Consequently, for the micro error analysis on the microscale the macrosolution is kept fixed and only the postprocessing is executed and enters the micro error analysis. 

{\color{black}Since the estimates \eqref{eq:Total-Error-estimate-L2}--\eqref{eq:Total-Error-estimate-Energy} are carried out on the macroscale and measured in macroscopic quantities, error analyses for their numerical verification equally operate on the macroscale.} 
 
The integrals for error calculation in the norms \eqref{eq:L2-norm}--\eqref{eq:energy-norm} are approximated by numerical integration of Gauss-Legendre. The computations are carried out on macro element level of the discretization for the reference solution. For the error in the $L^2$-norm it follows

\begin{equation}
\Vert \bm u^{{\color{black}H},\text{ref}} - \bm u^H \Vert _{L^2(\Omega)} \ = \ \left[ \sum_{K \in \mathcal{T}_{\text{ref}}}^{} \left( \sum_{i=1}^{ngp} \omega_i \left( \bm{u}^{{\color{black}H},\text{ref}}(\bm x_i^{\text{ref}}) - \bm u^H(\bm x_i^{\text{ref}}) \right)^2 \text{det} \bm J \right) \right]^{1/2} \, .
\label{eq:error_calc_num}
\end{equation}

For evaluating \eqref{eq:error_calc_num} the displacements of both the standard FE-HMM solution $\bm u^H$ and the reference solution $\bm u^{{\color{black}H},\text{ref}}$ have to be known in the quadrature points of the reference solution $\bm x_i^{\text{ref}}$. In the simplest case --when only the micro error {\color{black} convergence} is analyzed-- both solutions are computed on the same macro discretization {\color{black}$\mathcal{T}_{\text{ref}}=\mathcal{T}_{H}$} and the elements and their quadrature points therefore coincide, {\color{black}cf. Remark 1 (iv)}.

If the macro error or the total error is investigated, the reference solution has a finer macro triangulation than the single FE-HMM solutions. In this case the results of the FE-HMM solution are projected onto the finer grid of the reference solution.

%\begin{Figure}[htbp]
%	\centering
%	\psfrag{t1}{$\mathcal{T}_{\text{ref}}$}
%	\psfrag{t2}{\textcolor{red}{$\mathcal{T}_H$}}
%	\psfrag{xi}{$x_i^{\text{ref}}$}
%	\includegraphics[width=0.4\linewidth]{projection}
%	\\[6mm]
%	\caption{Projection from a coarse macro mesh (red) onto the quadrature points of an element (green) in the fine mesh of the reference solution (black) for linear shape functions.}
%	\label{fig:error_calc_projection}
%\end{Figure}

\begin{Figure}[htbp]
	\centering
%	\psfrag{t1}{$\mathcal{T}_{\text{ref}}$}
%	\psfrag{t2}{\textcolor{red}{$\mathcal{T}_H$}}
%	\psfrag{xi}{$x_i^{\text{ref}}$}
	\includegraphics[width=0.40\linewidth]{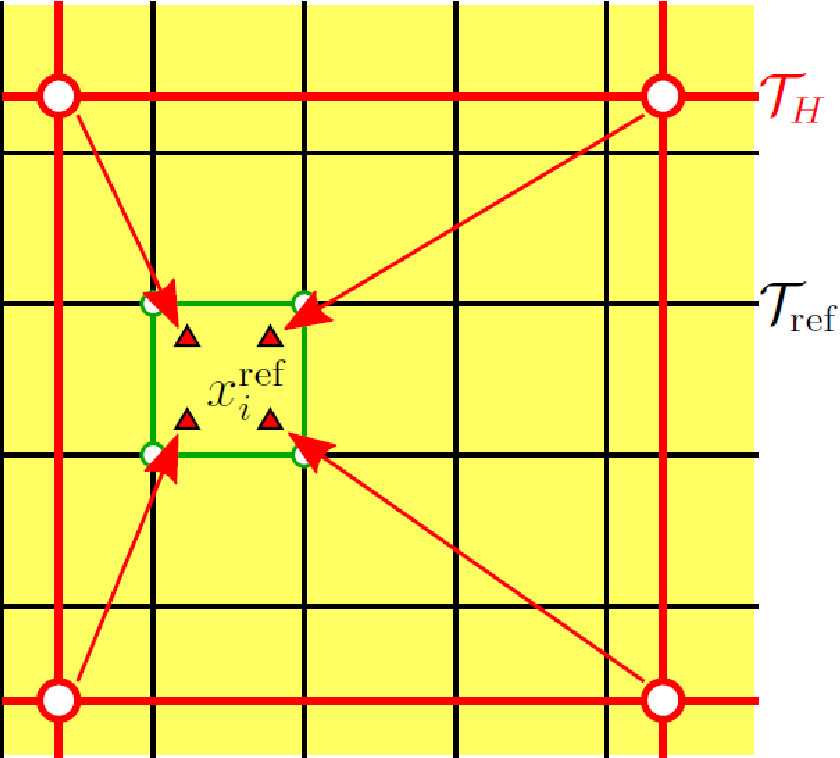}
%	\\[6mm]
	\caption{Projection from a coarse macro mesh (red) onto the quadrature points of an element (green) in the fine mesh of the reference solution (black) for linear shape functions.}
	\label{fig:error_calc_projection}
\end{Figure}%\end{Figure}

Figure \ref{fig:error_calc_projection} schematically displays the projection from a rather coarse macro triangulation $\mathcal{T}_H$ onto the finer reference triangulation $\mathcal{T}_{\text{ref}}$ for one element of the reference mesh. Therein, the quantities of the coarse mesh are projected onto the quadrature points of the reference solution $x_i^{\text{ref}}$ such that the error of the quantities of interest can be calculated, e.g. for the displacement error in the $L^2$-norm according to \eqref{eq:error_calc_num}.

\bigskip

\color{black}
If the absolute figures of the errors are of interest, the total discretization error, its macro and micro parts can be efficiently computed by merely two reference solutions as for example in Tab.~\ref{tab:rationale-for-error-decomposition-i}.
\begin{Table}[htbp]
{\color{black}  
\begin{minipage}{16.5cm}  
\centering
\renewcommand{\arraystretch}{1.2} 
\begin{tabular}{lll}
\hline
step & type of error                          & reference solution $\bm u^{{\color{black}H},\text{ref}}$ \\
\hline
1.) & total error at $H$ and $h$: ${e}_{tot}$ & $\bm u^{H}(H\rightarrow 0, h\rightarrow 0)$ \\
2.) & macro error: ${e}_{mac}$            & $\bm u^{H}(H, h\rightarrow 0)$. \\
3.) & resultant micro error: ${e}_{mic} = {e}_{tot} - {e}_{mic}$   & $[\bm u^{H}(H \rightarrow 0, h)]$   \\
\hline
\end{tabular} 
\end{minipage}
}
\caption{{\color{black} Rationale for the decomposition of the estimated error into its macro and micro parts.}}
\label{tab:rationale-for-error-decomposition-i} 
\end{Table}
 
The computation of reference solutions in 1.) and 2.) in Tab.~\ref{tab:rationale-for-error-decomposition-i} can efficiently be carried out in single-scale finite element simulations on the macroscale using the homogenized elasticity tensor $\mathbb A^{0,h}(h\rightarrow 0)$ determined in a  preprocessing step.

\color{black}
 
%---------------------------------------------------------------------------------------------------------
\subsection{Error estimation based on the Superconvergent Patch Recovery (SPR)}
%---------------------------------------------------------------------------------------------------------
\label{subsec:superconvergence}

In (engineering) practice, error computation as described in \ref{subsec:error-computation} is prohibitive. Instead, the total error is estimated, which is carried out on the particular discretization in use. For that purpose the present work uses the {\color{black}recovery}-type error estimation of Zienkiewicz and Zhu, which exploits superconvergence of stress and strain. In \cite{SPR}, \cite{SPR2} a procedure for the transfer of the superconvergence property from superconvergent, inner element points to element nodes referred to as ''superconvergent patch recovery'' (SPR) was proposed. Based on these recovered superconvergent nodal values the same authors constructed an error estimator that guided adaptive mesh refinement. 
 
\subsubsection{The SPR {\color{black} on the macroscale}}
 \label{subsubsec:SPR-for-MacroFEM}

For ready reference, the rationale of the SPR is briefly re-iterated for linear and quadratic shape functions, where we restrict to the $n_{dim}=2$ case for convenience. Strain and stress are calculated at superconvergent element sites that is for $p$=$1$ in the center of a rectangular element, for $p$=$2$ in the 2$\times$2 points of Gauss-Legendre quadrature. These values are transferred by a least-square procedure to the finite element node in the direct neighborhood, for a visualization see Fig.~\ref{fig:SPR_Elements}. Elements having such a node in common are referred to as the patch in the superconvergent recovery procedure.   

%\begin{Figure}[htbp]
%\centering
%\includegraphics[height=4.0cm]{patch_linear} \hspace*{10mm}
%\includegraphics[height=4.0cm]{patch_quad}
%\caption{Recovery of nodal stresses from stresses of surrounding superconvergent points (marked by a $\triangle$) for 4-node elements (left) and 9-node elements (right). The nodal stresses of the red marked nodes can be calculated using the shown patches.}
%\label{fig:SPR_Elements}
%\end{Figure}

\begin{Figure}[htbp]
\centering
\includegraphics[height=4.0cm]{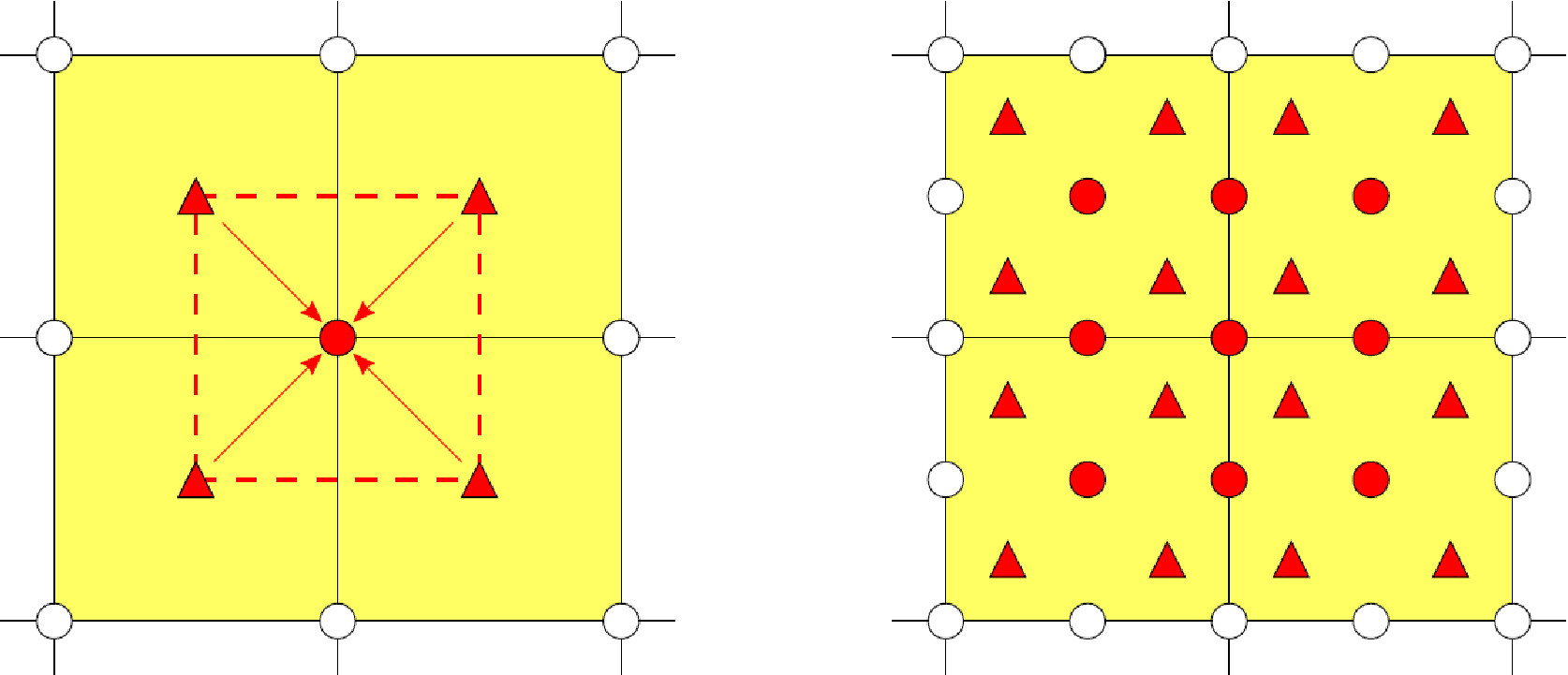} 
\caption{Recovery of nodal stresses from stresses of surrounding superconvergent points (marked by a $\triangle$) for 4-node elements (left) and 9-node elements (right). The nodal stresses of the red marked nodes can be calculated using the shown patches.}
\label{fig:SPR_Elements}
\end{Figure}

Stresses on the patch are prescribed component-wise by 
\begin{equation}
\sigma_p^\star = \ \bm P \,\bm a  
\label{eq_ch7_48}
\end{equation}
with, for the case of linear shape functions,  
\begin{equation}
\bm P \ = \ \left[ 1, \ x, \ y, \ xy \right] \quad \text{and} \quad \bm a \ = \ \left[ a_1, \ a_2, \ a_3, \ a_4 \right] \, .
\end{equation}
Vector $\bm P$ contains polynomial terms of bilinear shape functions for $n_{dim}=2$, no matter whether it is a 4-node or 9-node quadrilateral, since the patch around a finite element node consists of four superconvergent points in either case. For the determination of the unknown vector $\bm a$ the function
\begin{align}
F(\bm a) \ &= \ \sum_{i=1}^{n} \, (\sigma_h(x_i,\, y_i) - \sigma_p^\star (x_i, \, y_i))^2 \nonumber \\
&= \ \sum_{i=1}^{n} \, (\sigma_h (x_i, \, y_i) - \bm P (x_i, \, y_i)\mathbf{a} )^2
\end{align}
has to be minimized. Therein, $(x_i, \, y_i)$ are the coordinates of the superconvergent points, $n$ is the number of superconvergent points of the total patch and $\sigma_h (x_i, \, y_i)$ are the stresses in these superconvergent points. Minimization of $F(\bm a)$ implies that $\bm a$ fulfills the condition
\begin{equation}
\sum_{i=1}^{n} \, \bm P^T (x_i, \, y_i) \, \bm P (x_i, \, y_i) \, \bm a \ = \ \sum_{i=1}^{n} \, \bm P^T (x_i, \, y_i) \, \sigma_h (x_i, \, y_i) \, ,
\end{equation}
which can be solved for $\bm a$ 
\begin{equation}
\bm a \ = \ \bm A^{-1} \, \bm b
\end{equation}
with 
\begin{equation}
\bm A \ = \ \sum_{i=1}^{n} \, \bm P^T (x_i, \, y_i) \, \bm P(x_i, \, y_i) \quad \text{and} \quad \bm b \ = \ \sum_{i=1}^{n} \, \bm P^T (x_i, \, y_i) \, \sigma_h (x_i, \, y_i) \, .
\end{equation}
Stresses in the central node of the patch can be recovered by inserting its nodal coordinates $(x_N, \, y_N)$ into the $\bm P$-vector in \eqref{eq_ch7_48}.

%\begin{Figure}[htbp]
%	\centering
%	\includegraphics[height=4.0cm]{patch_boundary_linear} \hspace*{10mm}
%	\includegraphics[height=4.0cm]{patch_corner_linear}
%	\caption{Recovery of nodal stresses using patch recovery for boundary and corner nodes.}
%	\label{fig:SPR_boundary_nodes}
%\end{Figure}

\begin{Figure}[htbp]
	\centering
	\includegraphics[height=4.0cm]{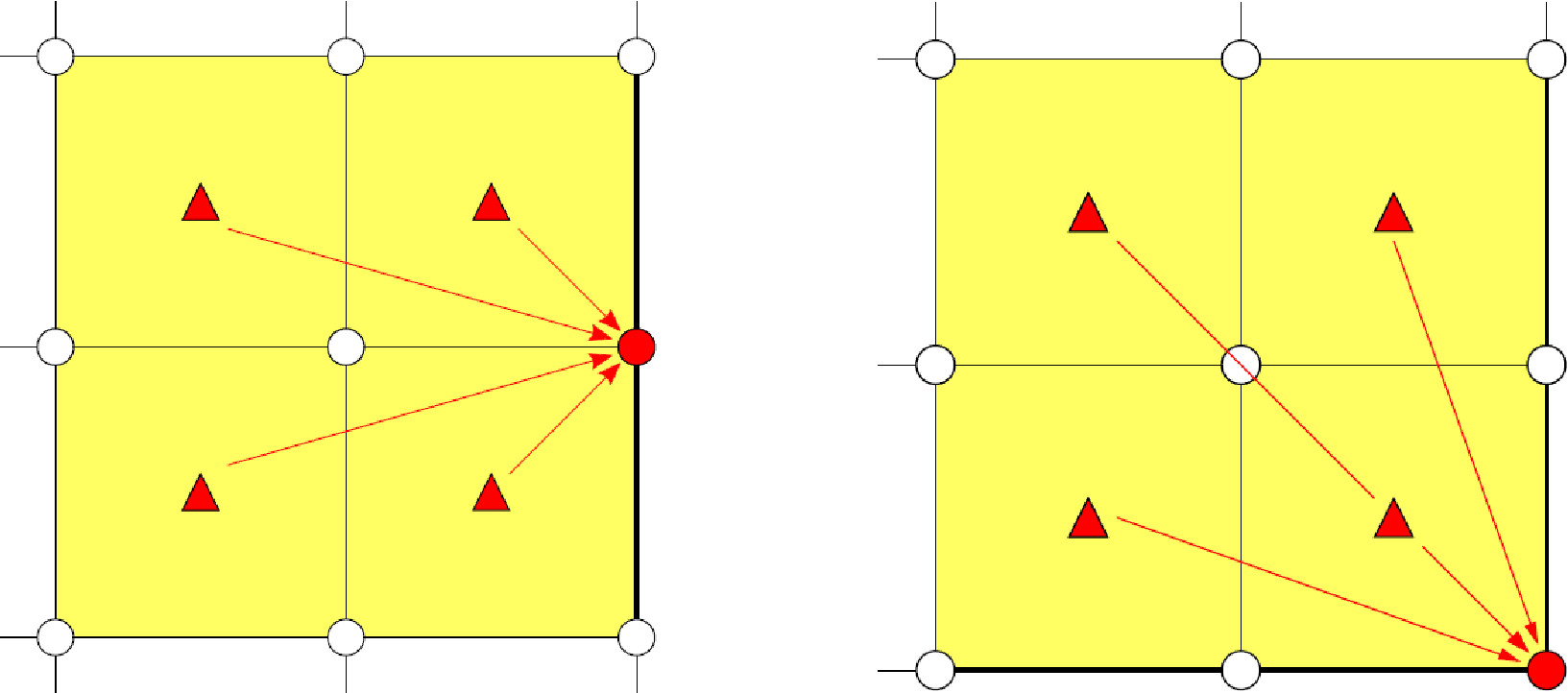}  
	\caption{Recovery of nodal stresses using patch recovery for boundary and corner nodes.}
	\label{fig:SPR_boundary_nodes}
\end{Figure}

Figure \ref{fig:SPR_boundary_nodes} shows the patches for boundary nodes lying either on edges or at corners. Corner nodes adjoin only one element which is insufficient for the calculation method described above. A similar situation arises for any node on the boundary which is adjoined to two elements. Here the patches have to be complemented by further elements.

Another peculiarity exists for patches of 9-node elements. For all of the red marked nodes in Fig. \ref{fig:SPR_Elements} --except of the central-one-- there are two or even more patches available to compute the nodal stresses. In this case the nodal values are calculated by simply averaging the results from the single patches.

\subsubsection{Error estimator and effectivity index}
\label{subsubsec:error_estimator}

As described above, the error estimator is built on superconvergent stress $\boldsymbol{\sigma}^\star$ and strain $\boldsymbol{\varepsilon}^\star$. Of course, the procedure is not applicable for an error estimate in the $L^2$-norm of displacements, since for the existing continuity of displacements the {\color{black}recovery}-type error estimator cannot be constructed. The estimated error in the energy-norm reads as

\begin{eqnarray}
|| \bar{\bm e} ||_{A(\Omega)} = || \bm u^\star - \bm u^H ||_{A(\Omega)} &=& \sqrt{\int_{\Omega} \left( \boldsymbol{\sigma}^\star - \boldsymbol{\sigma}^H \right) \colon \left( \boldsymbol{\varepsilon}^\star - \boldsymbol{\varepsilon}^H \right) \, dV} \, , \\
&\approx& \ \left[ \sum_{K \in \mathcal{T}_{H}}^{} \left( \sum_{i=1}^{ngp} \omega_i \left( \boldsymbol{\sigma}^\star - \boldsymbol{\sigma}^H \right)(\bm x_i^H) \colon \left( \boldsymbol{\varepsilon}^\star - \boldsymbol{\varepsilon}^H \right)(\bm x_i^H) \ \text{det} \bm J \right) \right]^{1/2} \, .
\label{eq:error_estimator}
\end{eqnarray}
  
Compared to the error computation based on a reference solution

\begin{equation}
|| \bm e ||_{A(\Omega)} = \Vert \bm u^{0} - \bm u^H \Vert _{A(\Omega)} \ \approx \ \left[ \sum_{K \in \mathcal{T}_{\text{ref}}}^{} \left( \sum_{i=1}^{ngp} \omega_i \left( \boldsymbol{\sigma}^{\text{ref}} - \boldsymbol{\sigma}^H \right)(\bm x_i^{\text{ref}}) \colon \left( \boldsymbol{\varepsilon}^{\text{ref}} - \boldsymbol{\varepsilon}^H \right)(\bm x_i^\text{ref}) \ \text{det} \bm J \right) \right]^{1/2}
\label{eq:error_estimator_vs_reference}
\end{equation}

the numerical effort of the error estimation is clearly much smaller, since the integration of the error is carried out on the corresponding macro mesh with triangulation $\mathcal{T}_H$ instead of the reference mesh with triangulation $\mathcal{T}_{\text{ref}}$. More important, error estimation can be carried out on-the-fly, no additional solution on a different mesh is required.

The quality of the error estimator is typically assessed by the so-called effectivity index $\theta$ which is defined as the ratio of the estimated error $\bar{ \bm e}$ to the true error $\bm e$ 

\begin{equation}
\theta = \dfrac{\Vert \bar{ \bm e} \Vert}{\Vert \bm e \Vert} \, .
\end{equation}

For consistency the effectivity index must tend to unity as the exact error tends to zero which can easily be shown if the error of stress and strain is considered. For the case of stresses entering the error analysis, the effectivity index follows to 

\begin{equation}
\label{eq:effectivity-index}
\theta = \dfrac{\Vert \bar{ \bm e}_{\sigma} \Vert}{\Vert \bm e_\sigma \Vert} = \dfrac{\Vert \bm \sigma^\star - \bm \sigma^H \Vert}{\Vert \bm \sigma^0 - \bm \sigma^H \Vert} = \dfrac{\Vert \bm \sigma^\star - \bm \sigma^0 + \bm \sigma^0 - \bm \sigma^H \Vert}{\Vert \bm \sigma^0 - \bm \sigma^H \Vert} \, .
\end{equation}
 
The numerator in \eqref{eq:effectivity-index} contains the error of standard stresses with respect to superconvergent stresses. A distinction of cases provides an upper and a lower bound for $\theta$

\begin{equation}
\label{eq:Theta-with-Bounds}
%\left( 
1 - \dfrac{\Vert \bm \sigma^\star - \bm \sigma^0 \Vert}{\Vert \bm \sigma^0 - \bm \sigma^H \Vert} 
%\right) 
\le \theta \le 
% \left( 
1 + \dfrac{\Vert \bm \sigma^\star - \bm \sigma^0 \Vert}{\Vert \bm \sigma^0 - \bm \sigma^H \Vert} 
%\right) 
\, .
\end{equation} 

Since the error of superconvergent quantities is expected to converge in higher order than the error of standard quantities, both bounds tend to unity as the error tends to zero.
\\[4mm]
{\color{black}  
{\bf Remark 2} 
\\[2mm]Imagine the case of error estimation on the macroscale for two different (macro-, micro-) discretizations, e.g. $(H,h)$ and $(H,h/2)$; the figures of the error estimates are expected to differ. Does this difference indicate that the error estimator on the macroscale includes the microdiscretization error? If not, why not and what else is indicated by the difference? 
\\[1mm] 
Here, the error estimator operates on the macroscale and exclusively estimates the macro discretization error at macro element size $H$ along with a given micro constitutive law\footnote{Recall, that the reference solution for the computation of the true macro error is $\bm u^{H}(H \rightarrow 0, h)$.}. The micro constitutive law is given in terms of its type and its material parameters, the latter depend on microdiscretization $h$. For the present case of linear elasticity the approximation $\mathbb{A}^{0,h}$ and its coefficients converge for sufficiently regular problems in the order $\mathcal{O}(h^{2q})$ to $\mathbb{A}^{0}$ for $h\rightarrow 0$ \cite{JeckerAbdulle2016}, \cite{EidelFischer2018}. On the microscale the deviation of $\mathbb{A}^{0,h/2}$ to $\mathbb{A}^{0,h}$ indicates a discretization error as the deviation of $\bm u^{h/2}$ to $\bm u^{h}$ does, which is hence accessible to an error estimator working on the microscale. In the error estimation working on the macroscale however, the deviation of $\mathbb{A}^{0,h/2}$ to $\mathbb{A}^{0,h}$ is not a discretization error but indicates a modeling-type error in terms of different constitutive laws, in the present setting in terms of different model parameters for the same \emph{type} of constitutive law.  
}

%---------------------------------------------------------------------------------------------------------
\section{Numerical examples}   
\label{sec:NumericalExamples} 
%---------------------------------------------------------------------------------------------------------
 
In this section a thorough convergence and error analysis is carried out for the three coupling conditions employing linear and quadratic shape functions on both the macro- and the microscale.   
 
First, the three micro-coupling conditions are compared in the microscale setting of (i) a matrix-inclusion problem, (ii) a chessboard-type microstructure, and (iii) a sine wave distribution of material stiffness, where the strength and quality of the stiffness contrast between different phases and its impact on the convergence properties is a key aspect of investigation. 

The convergence analysis measures the micro error both on the microscale (i.e. on one microdomain) and on the macroscale (as the total micro error that is propagated to the macroscale). Estimates for sufficiently regular problems are provided in Sec.~\ref{subsec:A-priori-estimates}. For the microerror on the macroscale the order $2q$ is expected in all three norms, see \eqref{eq:Total-Error-estimate-L2}--\eqref{eq:Total-Error-estimate-Energy}. For the micro error as measured on the microscale, order $q+1$ is expected in the $L^2$-norm, and $q$ in the $H^1$- and energy-norm.  

Similarly, the regularity of the macro-BVP is examined through the convergence of the macro error. Here, the results of a clamped, square plate with low regularity due to notch effects at the clamped boundary are contrasted to a tapered cantilever of proven high regularity.  

Moreover, the above examples, which have plane strain conditions and loading by volume forces in common, serve the purpose to compare the estimated error with the true error and to verify the optimal uniform micro-macro mesh refinement strategies of Tab.~\ref{tab:Best-mic-mac-RefinementStrategies}.
 
\subsection{Micro convergence analysis}

The macro problem common to all micro problems is a square cantilever subject to a volume forces of $\bm f = [0, -10]^T$ $[F/L^2]$. The coupling conditions which will be analyzed are Dirichlet, Neumann and periodic coupling. 

\subsubsection{Matrix-inclusion problem}
\label{subsubsec:mat-incl}

\begin{Figure}[htbp]
	\centering
	\includegraphics[width=0.32\linewidth]{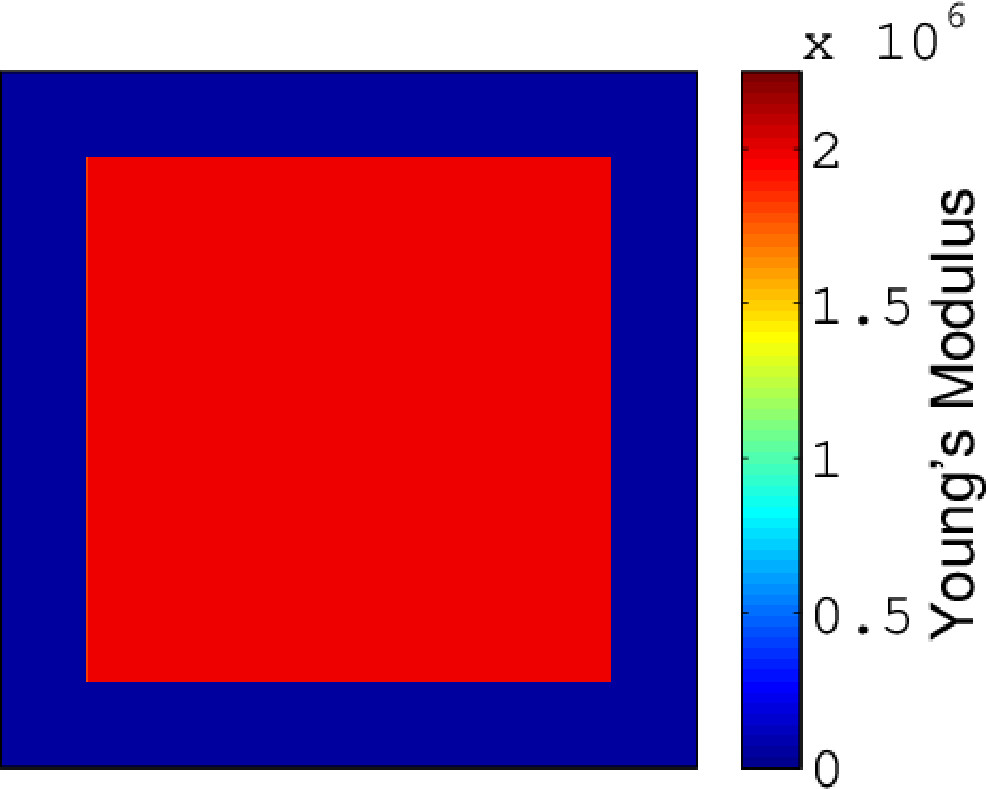}
	\caption{\textbf{Matrix-inclusion microstructure.} Distribution of Young's modulus on the micro domain.}
	\label{fig:matincl}
\end{Figure}

In the first numerical example we consider the microstructure of a stiff inclusion in a soft matrix, which is displayed in Fig. \ref{fig:matincl}. The Young's moduli of the inclusion $E_i = 200\,000$ $[F/L^2]$ and the matrix phase $E_m = 40\,000$ $[F/L^2]$ exhibit a contrast of $E_i/E_m = 50$, for the Poisson's ratio it holds $\nu = 0.2$. The volume ratio of the inclusion phase is $V_i/V_{tot} = 9/16$. {\color{black}The square RVE exhibits side length $\epsilon=0.005$, which is maintained for all examples in the present work.}

\begin{Figure}[htbp]
	\centering
	\includegraphics[width=0.32\linewidth]{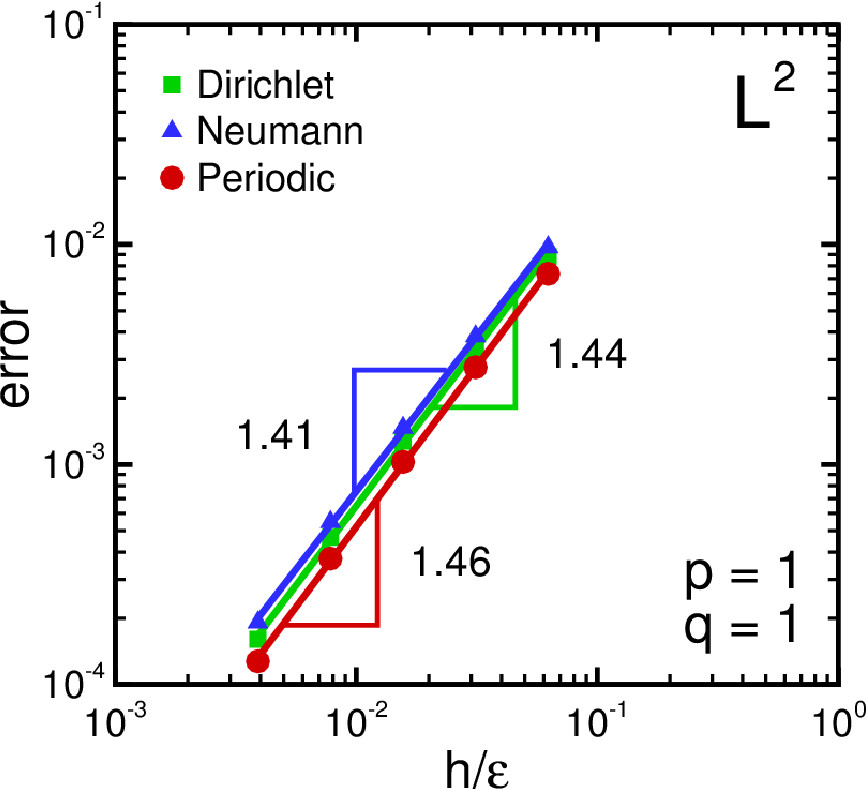} \hfill
	\includegraphics[width=0.32\linewidth]{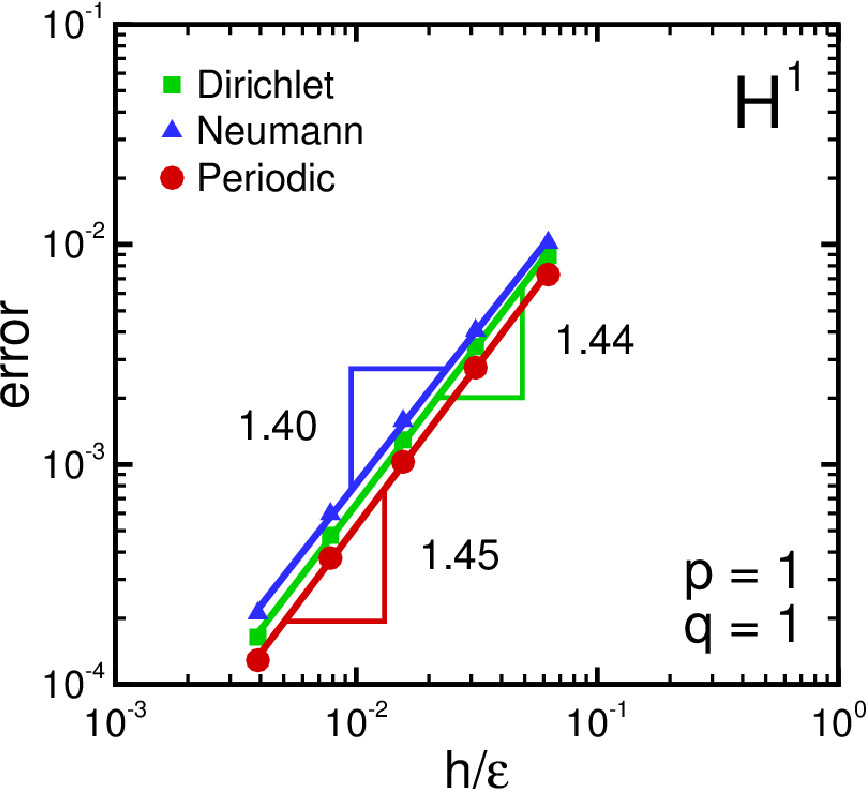} \hfill
	\includegraphics[width=0.32\linewidth]{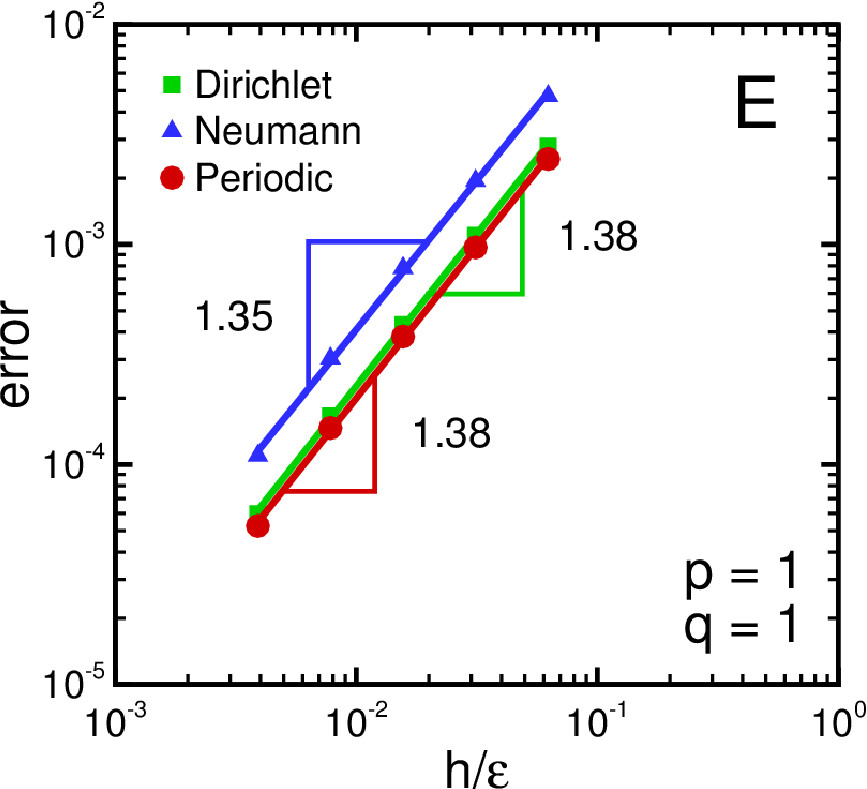} \\
    \includegraphics[width=0.32\linewidth]{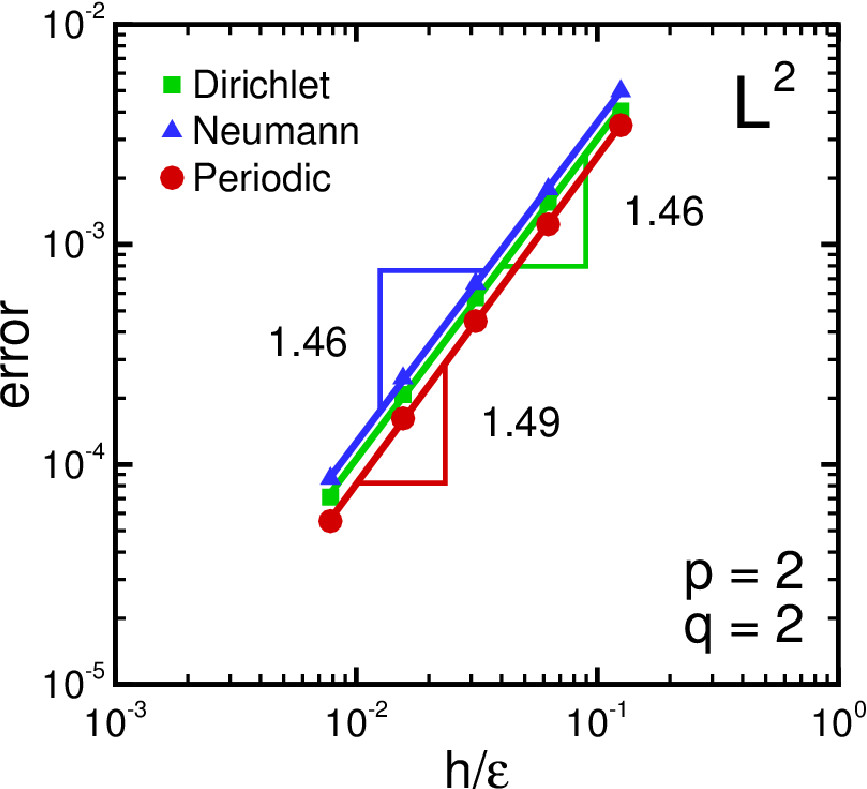} \hfill
	\includegraphics[width=0.32\linewidth]{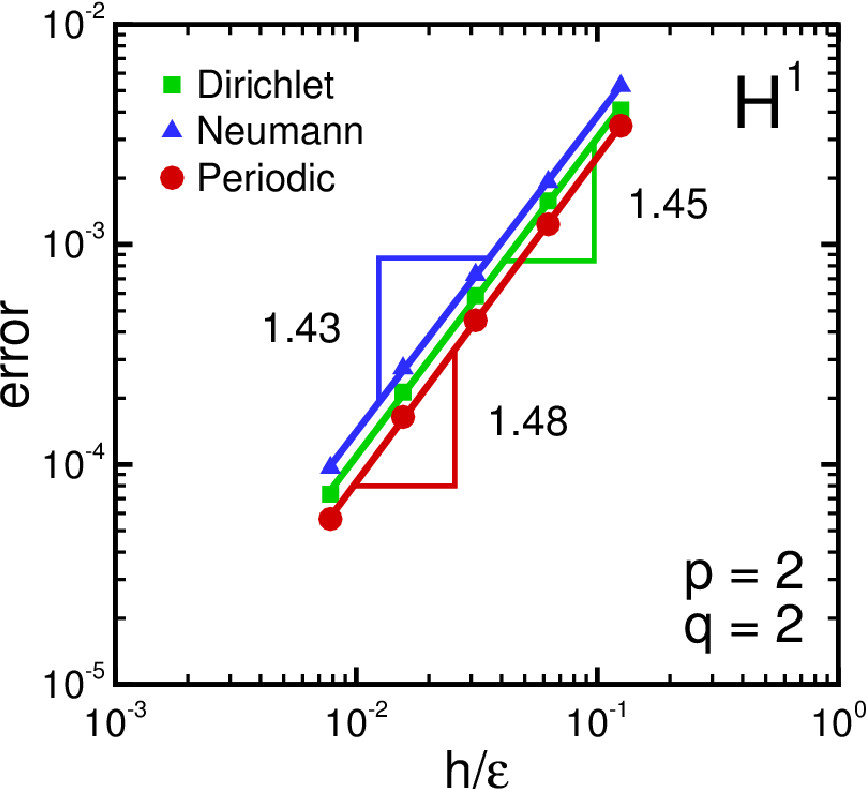} \hfill
	\includegraphics[width=0.32\linewidth]{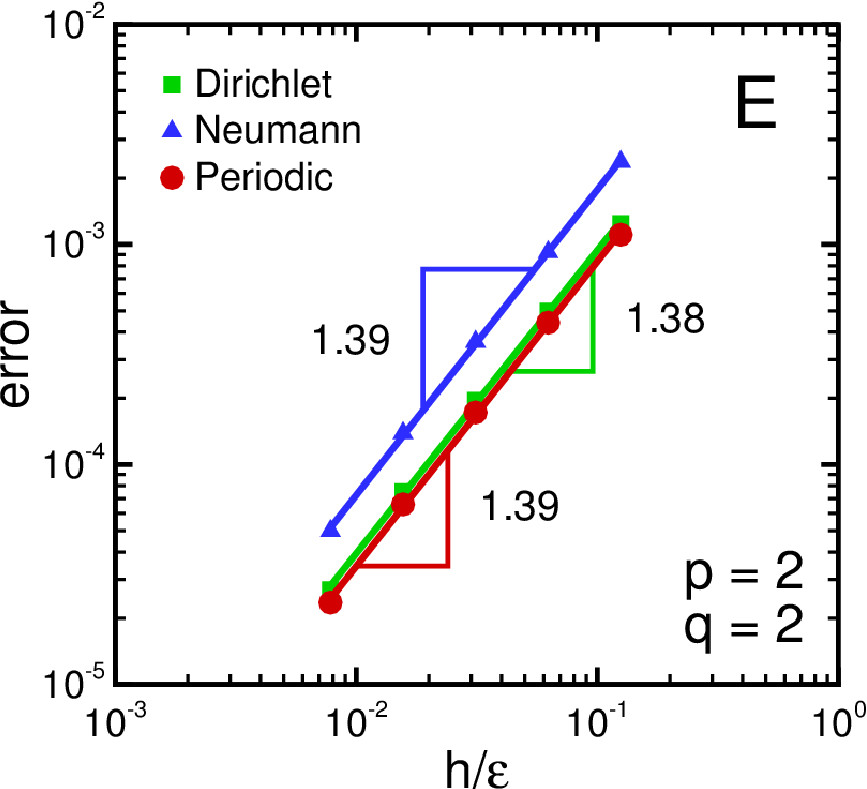}
	\caption{\textbf{Micro error convergence on the macroscale for matrix-inclusion problem:} (first row) linear shape functions $p$=$q$=$1$, (second row) quadratic shape functions $p$=$q$=$2$, from left to right: $L^2$-, $H^1$-, energy-norm.}
	\label{fig:diagramm_matincl_q1}
\end{Figure}
  
The simulation results for linear shape functions are displayed in Fig. \ref{fig:diagramm_matincl_q1} (first row). The different coupling conditions show minor deviations from each other in the convergence order. The values of the calculated errors are in good agreement between all coupling conditions, only the error for Neumann coupling in the energy-norm is slightly larger. Notice that we use here and in the following relative errors, i.e. $|| \bm e ||_{(\Omega)}/|| \bm u ||_{(\Omega)}$.

The observed order reduction from theoretical order $2q=2$ for $q=1$ to approximately 1.4 in all three norms is due to the stiffness-jump at the inclusion-matrix interface along with the high contrast of the Young's moduli of the two phases, the corresponding notch effect lowers the regularity of the microproblem.

The diagrams in Fig. \ref{fig:diagramm_matincl_q1} (second row) display the simulation results for quadratic shape functions, $p$=$q$=$2$. In all of the above described aspects we observe even quantitatively almost the same behavior as for $p$=$q$=$1$. Hence, it is the singularity of the problem which overrules the theoretical convergence order, i.e. quadratic shape functions do not cure the problem of low regularity.
 
\begin{Figure}[htbp]
	\centering
	\includegraphics[width=0.32\linewidth]{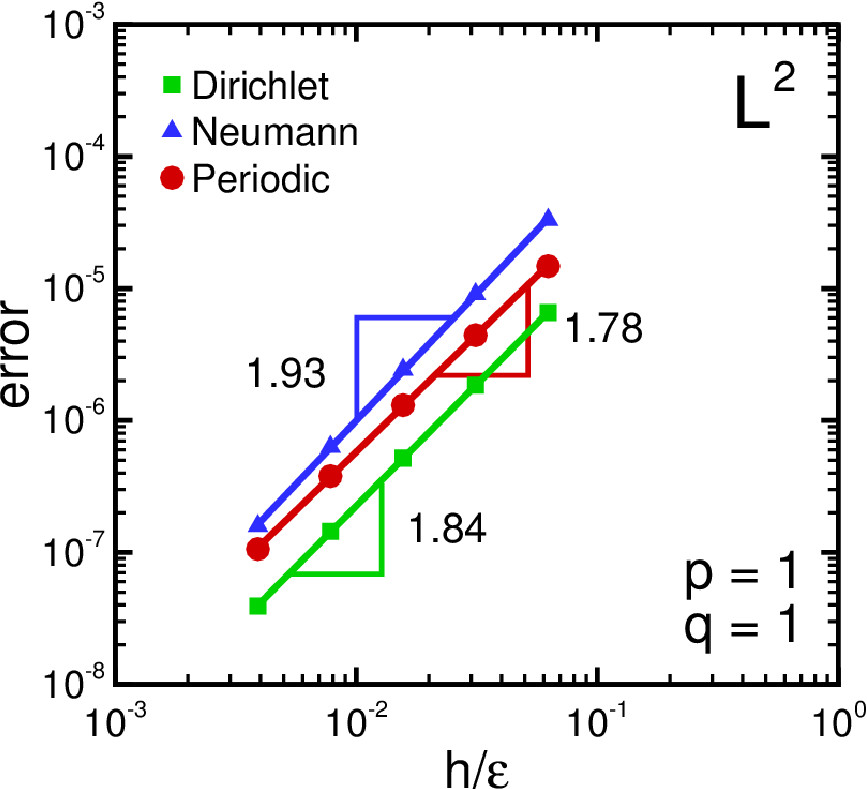} \hfill
	\includegraphics[width=0.32\linewidth]{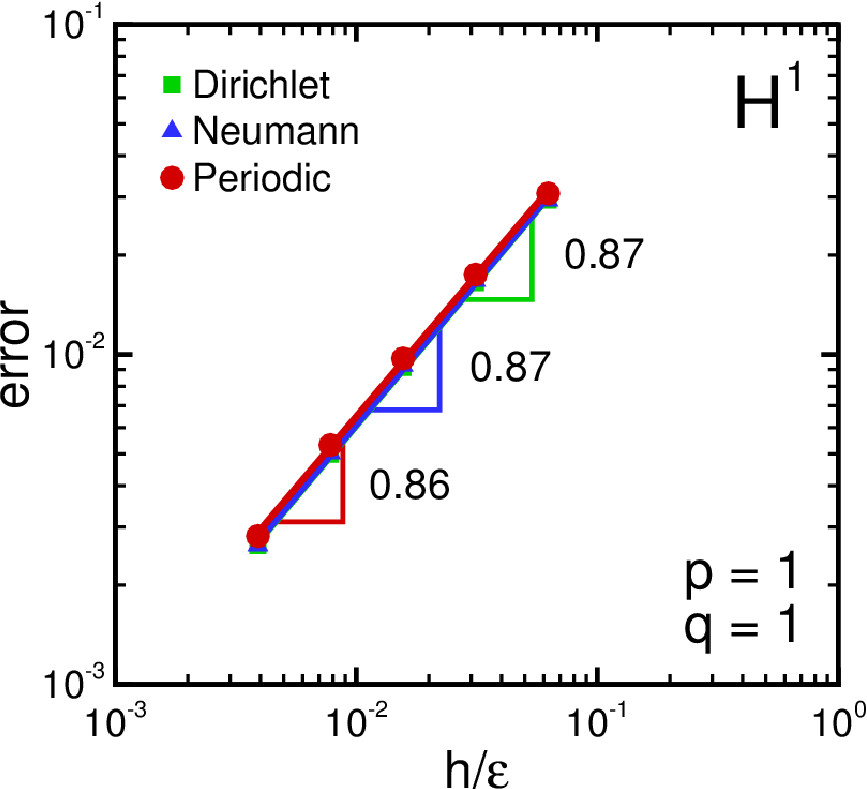} \hfill
	\includegraphics[width=0.32\linewidth]{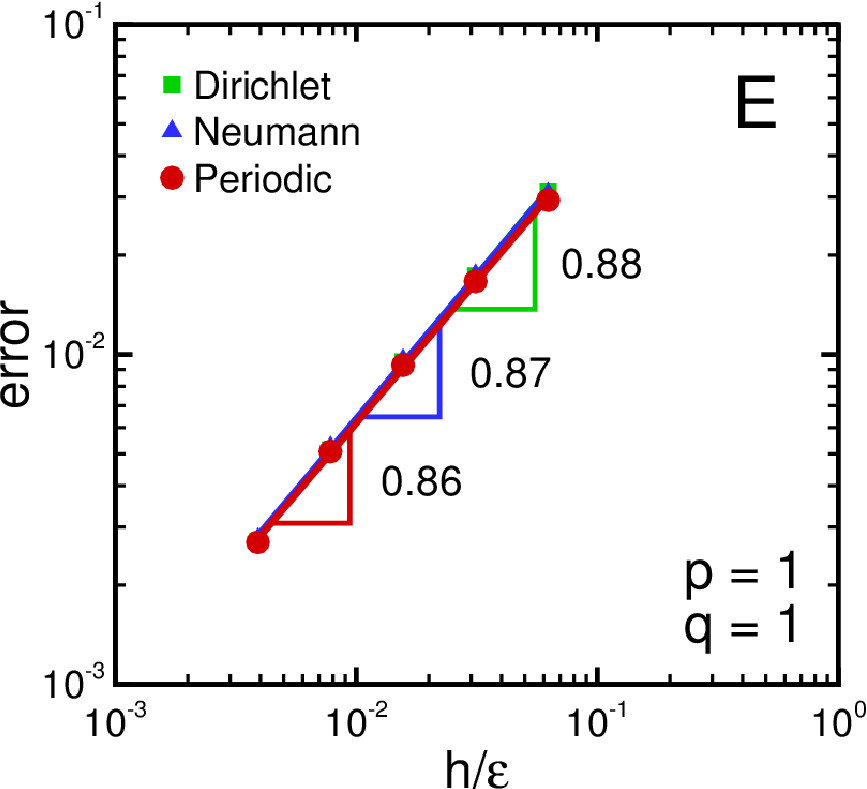} \\
    \includegraphics[width=0.32\linewidth]{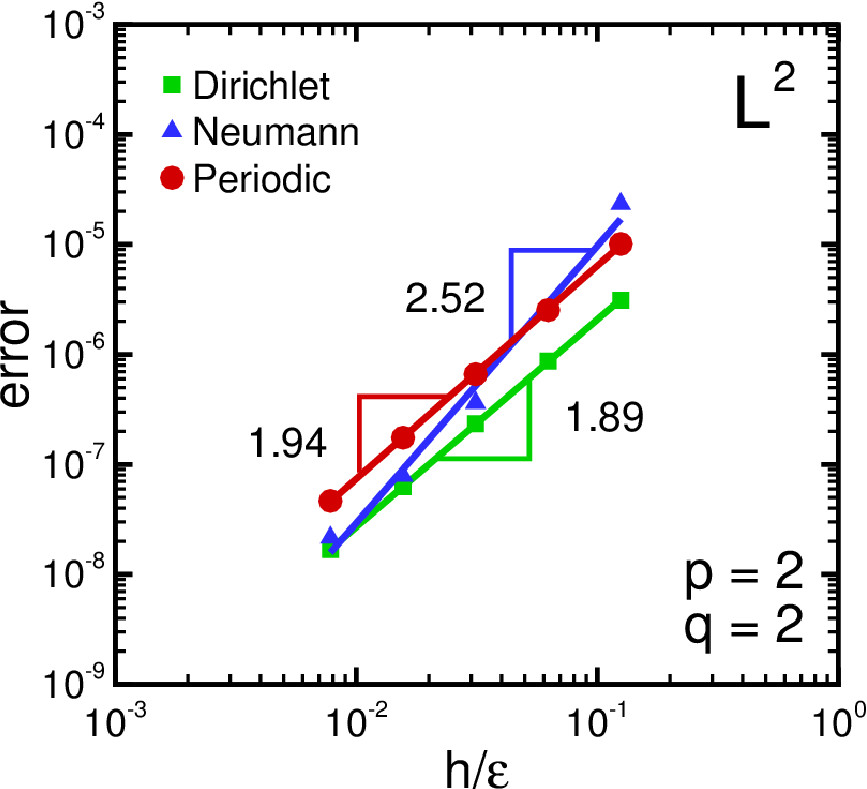} \hfill
	\includegraphics[width=0.32\linewidth]{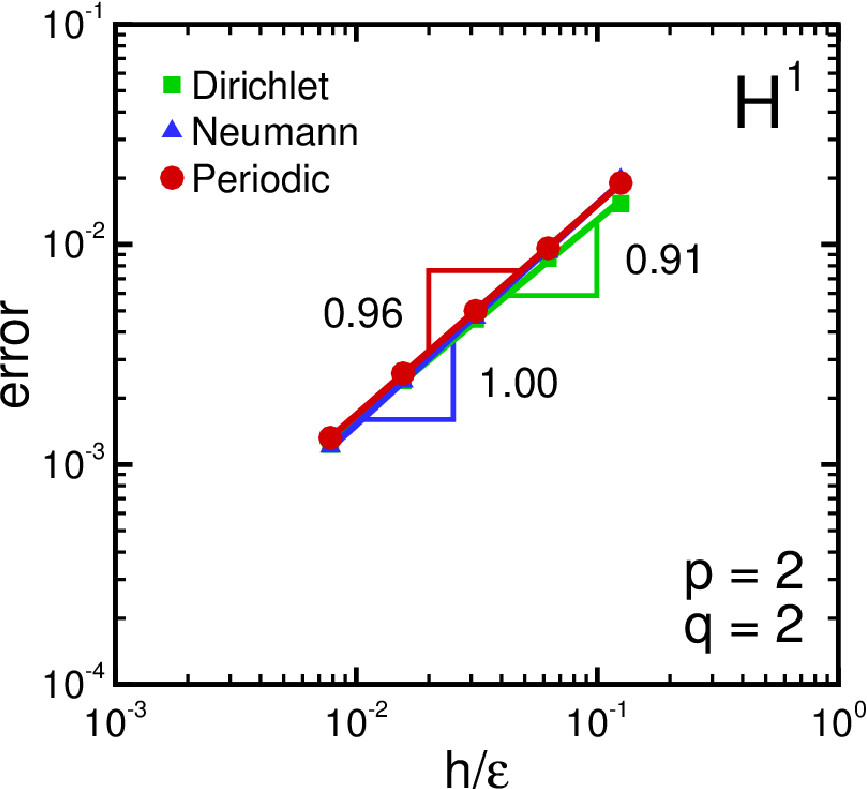} \hfill
	\includegraphics[width=0.32\linewidth]{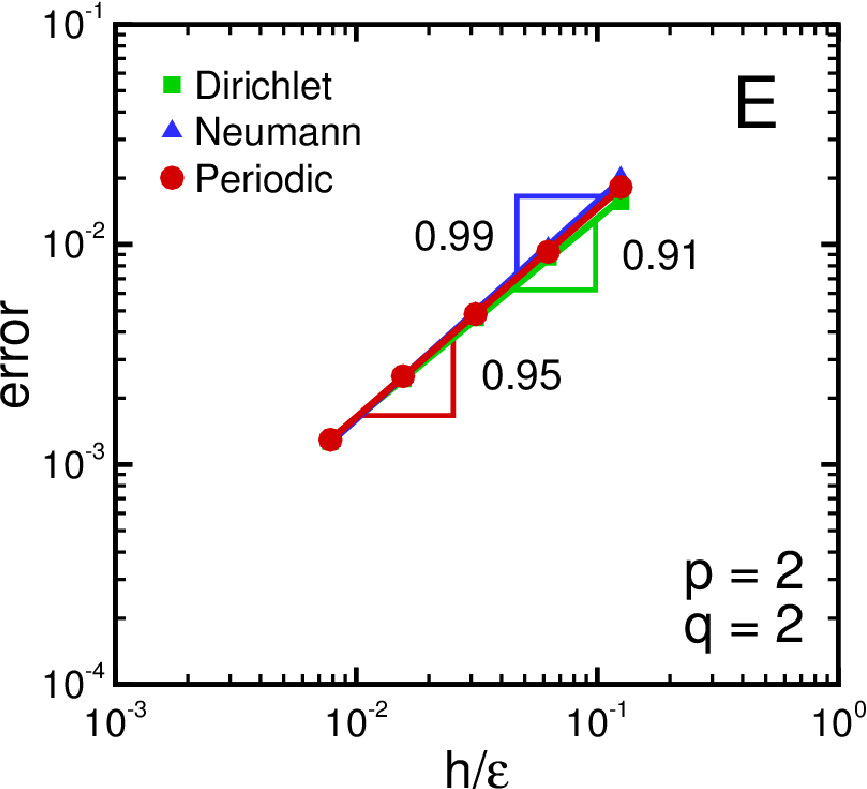}
	\caption{\textbf{Micro error convergence on the microscale for matrix-inclusion problem.} (first row) linear shape functions, (second row) quadratic shape functions, (from left to right:) $L^2$-, $H^1$- and energy-norm.}
	\label{fig:diagramm_matincl_q1_microdomain}
\end{Figure}

{\bf Micro error convergence on the microscale.} The results of an error calculation on one microdomain attached to the macroscopic quadrature point at [0.26, 0.26] is shown in the first row diagrams of Fig.~\ref{fig:diagramm_matincl_q1_microdomain} for linear shape functions. The convergence orders for different coupling conditions are in good agreement with each other. 

The optimal convergence orders of $q+1$ in the $L^2$-norm and $q$ in the $H^1$- as well as in the energy-norm is not reached due to the above mentioned reduced regularity. 

Figure~\ref{fig:diagramm_matincl_q1_microdomain} (second row diagrams) displays the errors on one microdomain for quadratic shape functions. Again the convergence orders of the different coupling conditions are in good agreement with each other except of for Neumann coupling in the $L^2$-norm. A closer look at these calculated errors reveals that for coarse discretizations the error is too large and for that reason converges faster than expected. If only the two finest discretizations are considered, convergence is in reasonable agreement with the other two coupling conditions.

%----------------------------------------------------------------------------------------------------------------------
\subsubsection{Comparison of methods for constant traction BC: semi-Dirichlet coupling versus perturbation technique}

As described in Sec.~\ref{subsec:Implementation-NeumannCoupling} two different techniques are considered and compared, which fulfill the constant traction BC and remove the rigid body motions from the RVE. The methods are the semi-Dirichlet coupling introduced by \cite{JaviliSaebSteinmann2017} and the perturbation technique going back to \cite{Miehe-Koch-2002}. For a comparison of the methods the above matrix-inclusion problem is considered but for visualization purposes the applied load and the microdomain size are increased compared to Sec.~\ref{subsubsec:mat-incl}. The methods are compared, first with respect to their accuracy, second with respect to kinematical implications of removing the rigid body motions from the RVE.  
   
\begin{Table}[htbp]
	\centering
\begin{tabular}{c c c c c c c}
	\hline
 $p$=$q$ & h/$\epsilon$ & 1/16 & 1/32 & 1/64 & 1/128 & 1/256 \\[1mm] 
         \hline 
    1    & semi-Dirichlet [in $10^{-6}$] & 28.3848 & 28.4078 & 28.4151 & 28.4173 & 28.4180 \\ 
         & perturbation [in $10^{-6}$]   & 28.3848 & 28.4078 & 28.4151 & 28.4173 & 28.4180 \\ 
	\hline 
    2    & semi-Dirichlet [in $10^{-6}$] & 28.5173 & 28.5265 & 28.5286 & 28.5291 & 28.5292 \\ 
         & perturbation [in $10^{-6}$]   & 28.5173 & 28.5265 & 28.5286 & 28.5291 & 28.5292 \\
	\hline 	
\end{tabular} 
\caption{\textbf{Constant traction BC: Comparison of semi-Dirichlet coupling and perturbation technique.} $L^2$-normed solution vector at different micro discretizations for linear and quadratic shape functions.}
\label{tab:comparison_semiDirichlet_vs_perturbation}
\end{Table}

First and foremost, both methods accurately fulfill the constant traction BC. Table \ref{tab:comparison_semiDirichlet_vs_perturbation} shows the $L^2$-norm of the solution vectors in the RVE for the two methods indicating that they yield the same results for various discretizations. For the perturbation technique, the parameters are chosen randomly with a maximum value of $10^{-5}$. If the perturbation parameter is varied in the range from $10^{-8}$ to $10^0$ for linear shape functions (and in the range of $10^0$ to $10^{-7}$ for quadratic shape functions), the results show only minor deviations for the large perturbation parameter. For the choice of $10^{-1}$ in the case of linear shape functions the $L^2$-norm for $h/\epsilon = 1/256$ amounts to 28.4179, for $10^0$ to 28.4050. 
 
\begin{Figure}[htbp]
	\centering
	\includegraphics[width=0.27\linewidth]{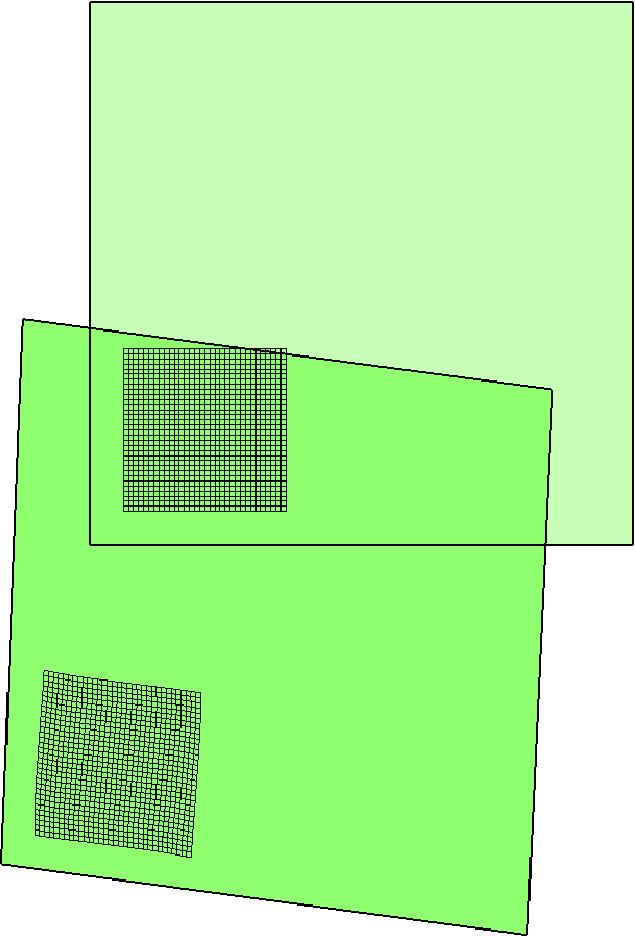} \hspace*{0.5cm}
	\includegraphics[width=0.27\linewidth]{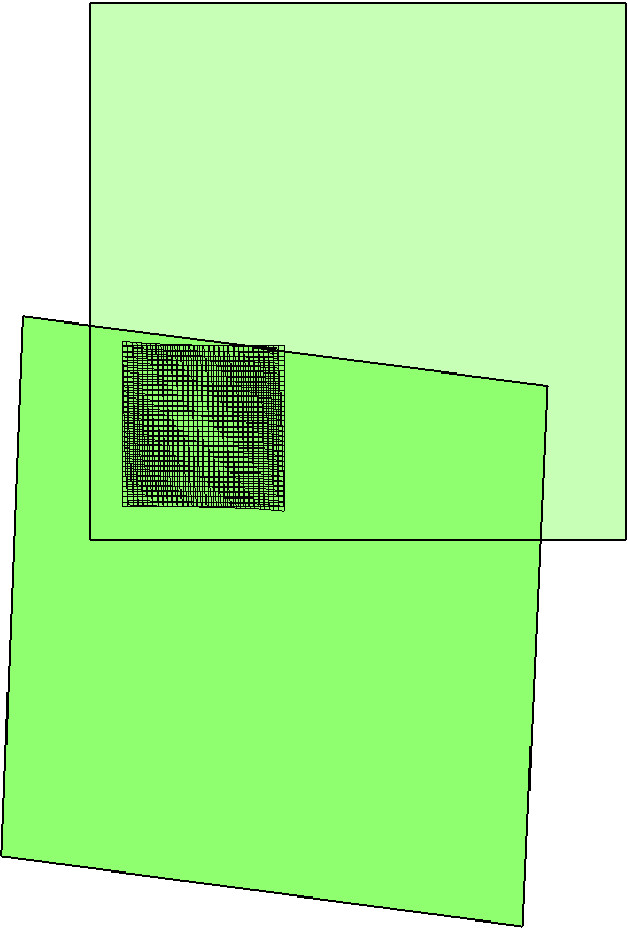} \hspace*{0.5cm}
	\includegraphics[width=0.27\linewidth]{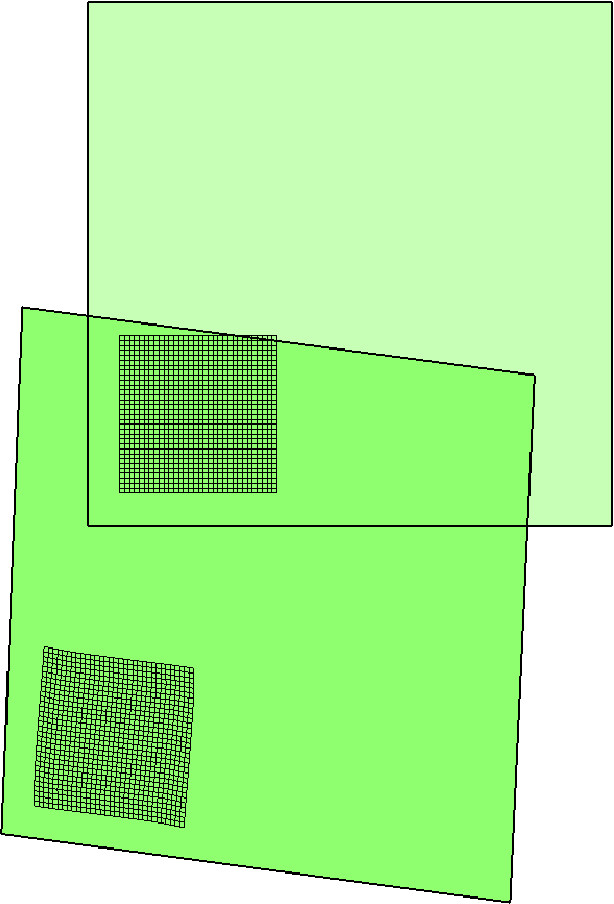}
	\caption{\textbf{Micro displacements for matrix-inclusion problem at constant tractions BC} Deformed and undeformed macroelement and microdomain for (left) semi-Dirichlet coupling and for the perturbation technique (center) without and (right) with  adding macroscopic rigid body motions.}
	\label{fig:MatIncl_microdisplacements_semi_dirichlet_vs_perturbation}
\end{Figure} 

Figure \ref{fig:MatIncl_microdisplacements_semi_dirichlet_vs_perturbation} (left) displays for the semi-Dirichlet coupling the macroelement and the RVE at the lower left quadrature point in the undeformed and deformed configurations. For that case the micro displacements obviously fit into the macroscopic displacement field. Figure~\ref{fig:MatIncl_microdisplacements_semi_dirichlet_vs_perturbation} (center) similarly displays the same macroelement and RVE for the perturbation technique. It is obvious that the calculated micro displacements lack the rigid body motions following from the macroscopic deformation. In order to add the missing kinematical embedding, the displacement of the corresponding macroscopic quadrature point and the rotation of the macroscopic element can be added to the calculated micro displacement field, which results in the deformed configuration of Fig.~\ref{fig:MatIncl_microdisplacements_semi_dirichlet_vs_perturbation} (right).

{\color{black}

An important aspect is the efficiency of the two methods. It is obvious that the semi Dirichlet coupling is more expensive due to the fact that the system of equations has to be solved more than once in contrast to the pertubation technique. 

\begin{Figure}[htbp]
	\centering
	\includegraphics[height=0.28\linewidth]{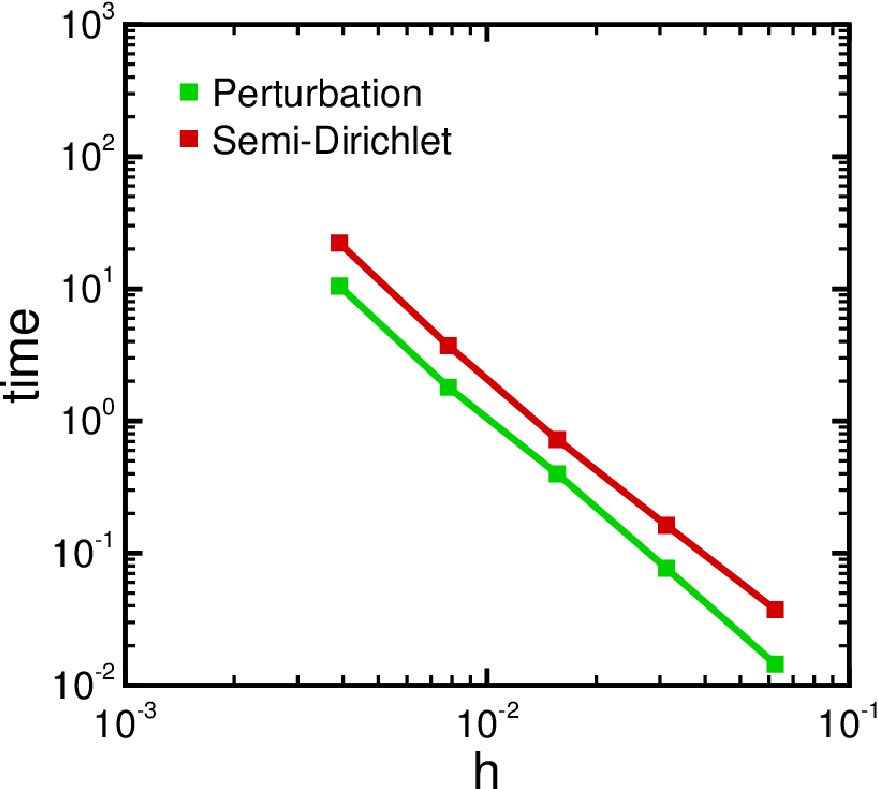} \hspace*{10mm}
	\includegraphics[height=0.28\linewidth]{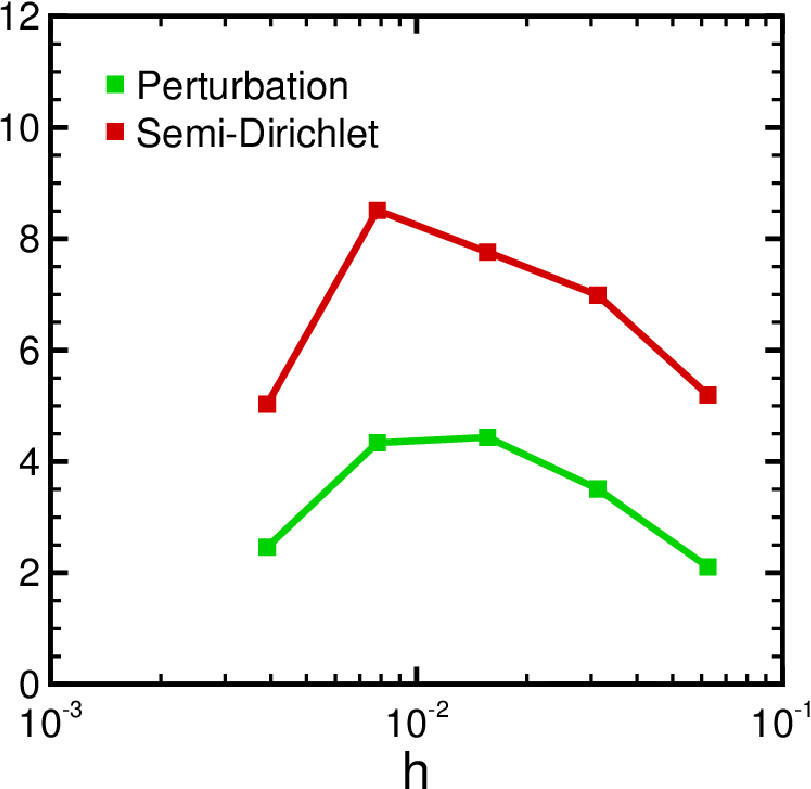} \caption{{\color{black}\textbf{Computational time for perturbation technique and semi-Dirichlet.} Absolute computational time (left) and relative percentaged computational time related to full computation time on the micro level.}} 
	\label{fig:perturbation_vs_semi-dirichlet}
\end{Figure} 

Figure \ref{fig:perturbation_vs_semi-dirichlet} displays the computational times of both methods, with absolute values in the left and, in the right, the percentage of the total computational time (including the stiffness matrix calculation) on the micro level. The results show that there is a difference between the two methods methods. 
} 

In conclusion, both techniques accurately fulfill the constant traction BC and yield the same microscopic stresses and strains. The perturbation technique is accurate for a wide range of the perturbation parameters. In this context it should also be mentioned that the semi-Dirichlet coupling method is robust with respect to the choice of the nodes in the RVE to which the additional Dirichlet constraints are applied. The only difference between both methods is the embedding of the rigid body motions following from the macroscopic displacements in the semi-Dirichlet technique. For the calculation of the microdisplacements which are used in the transformation matrix \eqref{eq:k-mac-element-6} it does not matter that the perturbation technique lacks the rigid body motions. If not only microscopic stresses and strains are of interest but equally the microscopic displacements including the macroscopic displacement state, either semi-Dirichlet coupling can be used or the perturbation technique along with rigid body motions enriched kinematics as detailed above. The lower numerical effort favors the iteration-free perturbation technique. 
 
\subsubsection{Chessboard microstructure}

\begin{Figure}[htbp]
	\centering
	\includegraphics[width=0.33\linewidth]{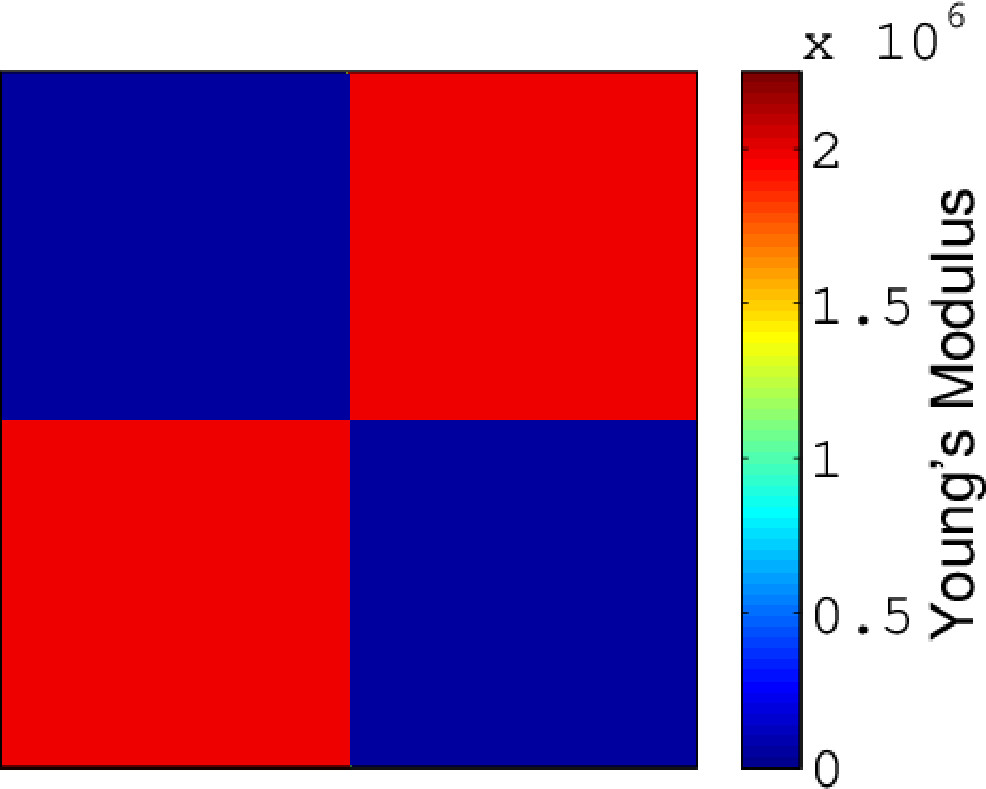}
	\caption{\textbf{Chessboard-type microstructure.} Distribution of Young's modulus on the micro domain.}
	\label{fig:chessboard}
\end{Figure}

While in the first example of the matrix-inclusion problem the material at the RVE boundary was homogeneous, we choose a chessboard-type microstructure, where the heterogeneity is expanded from the micro domain's interior to its boundaries, see  Fig. \ref{fig:chessboard}. The aim is to investigate the impact of micro-coupling conditions on the results for that case. 

The chessboard pattern of Young's modulus distribution exhibits two phases with $E_1 = 2000000$ $[F/L^2]$ and $E_2 = 40000$ $[F/L^2]$. The stiffness contrast of the phases is $E_1/E_2 = 50$, for the Poisson's ratio it holds $\nu=0.2$. 

\begin{Figure}[htbp]
	\centering
	\includegraphics[width=0.32\linewidth]{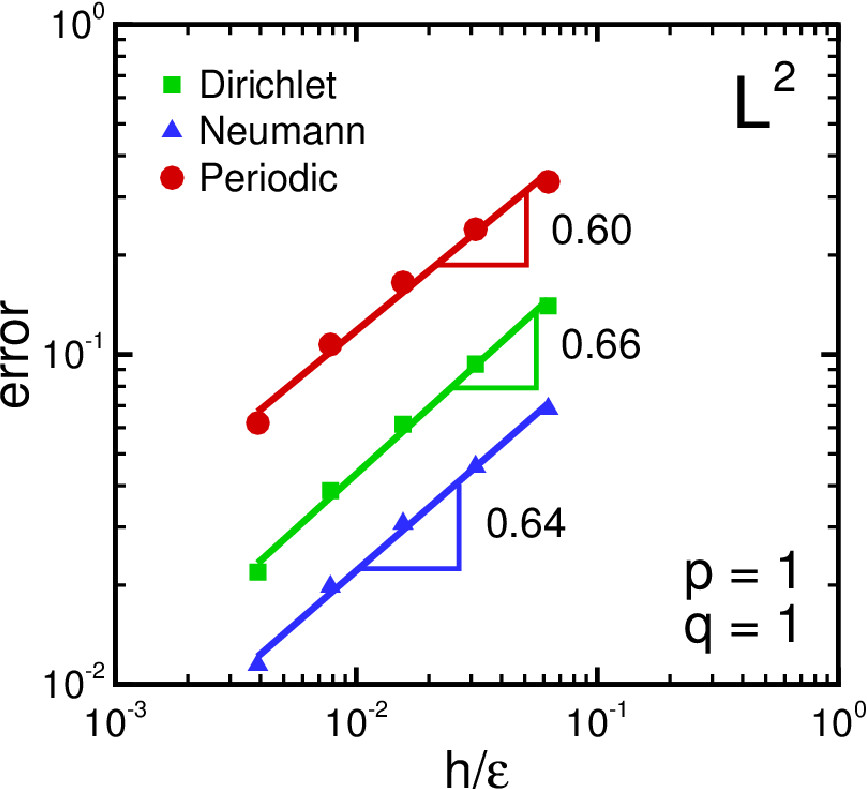} \hfill
	\includegraphics[width=0.32\linewidth]{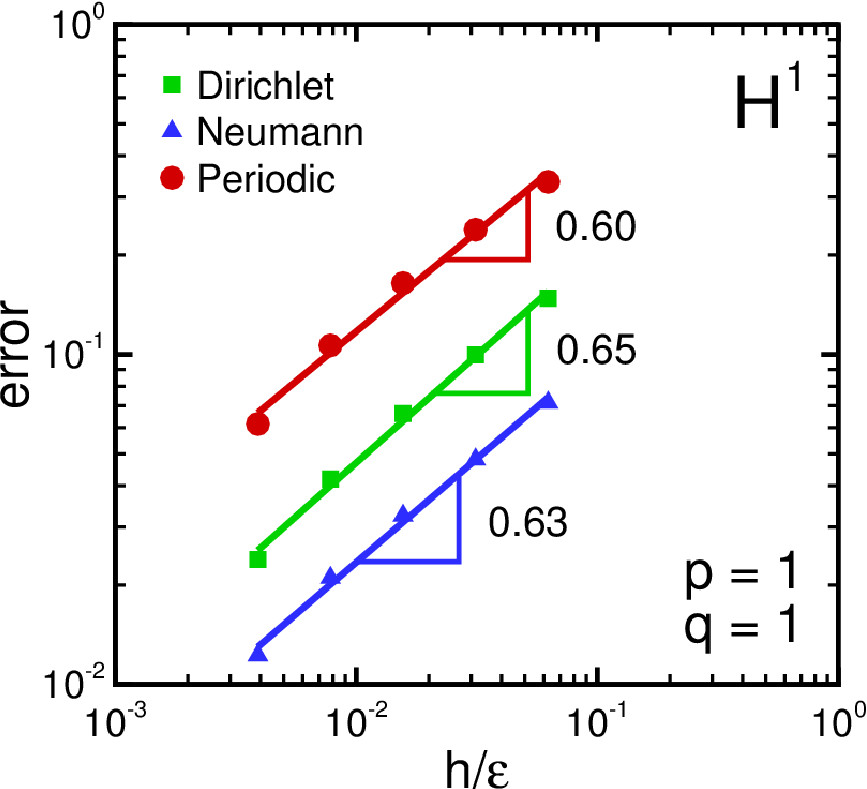} \hfill
	\includegraphics[width=0.32\linewidth]{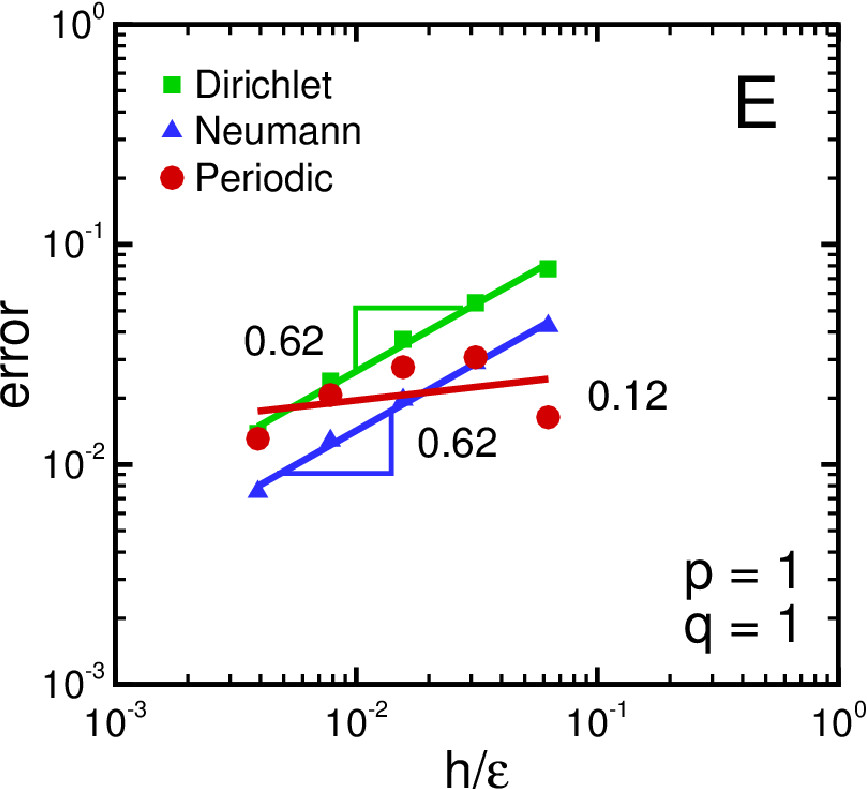} \\
	\includegraphics[width=0.32\linewidth]{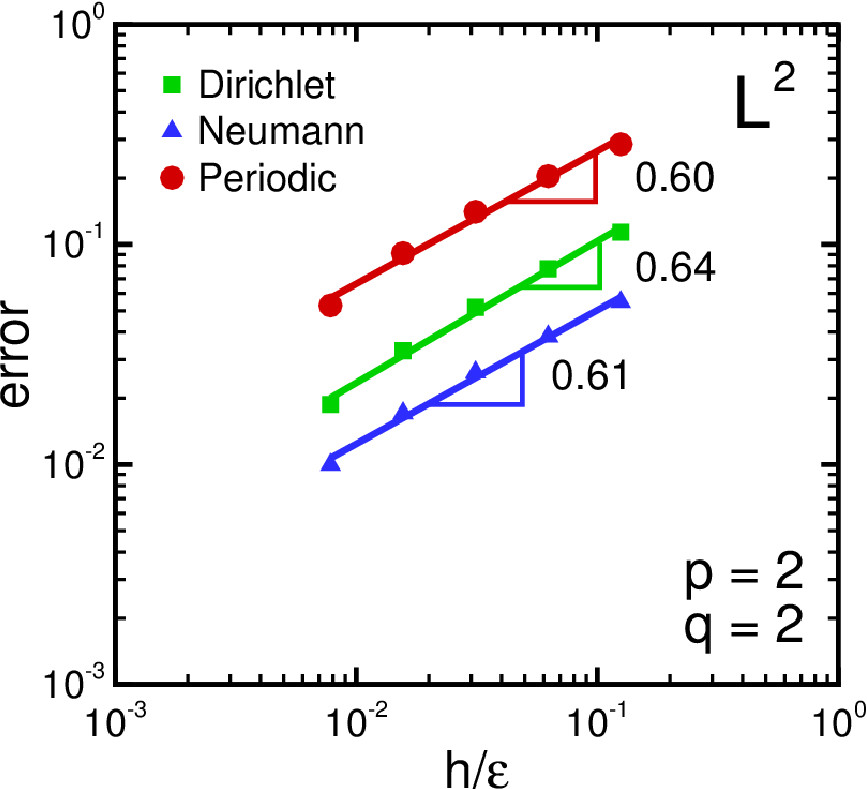} \hfill
	\includegraphics[width=0.32\linewidth]{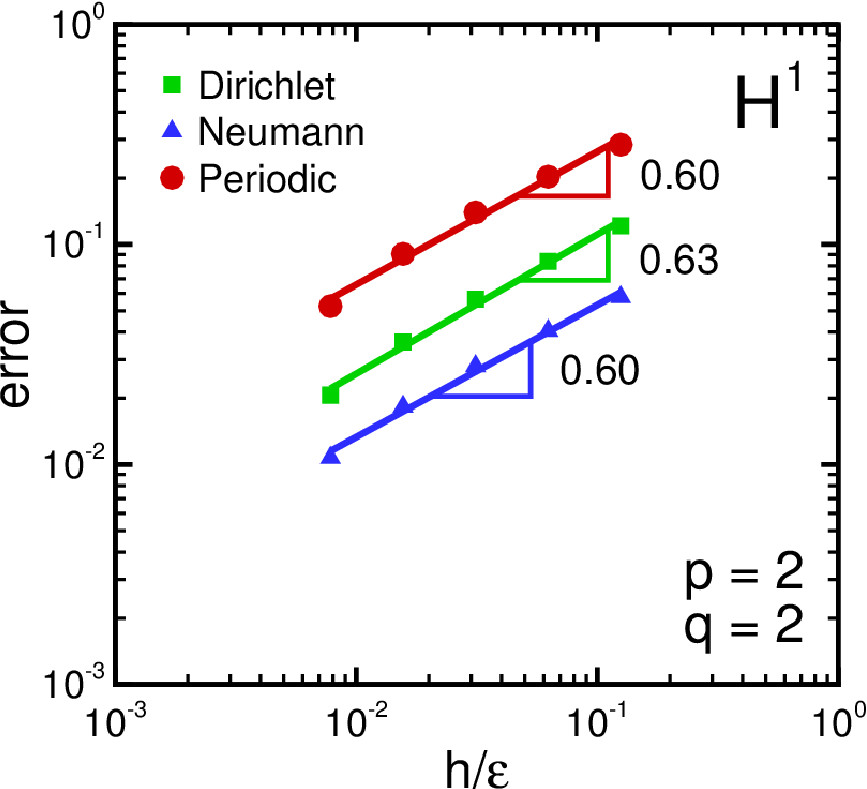} \hfill
	\includegraphics[width=0.32\linewidth]{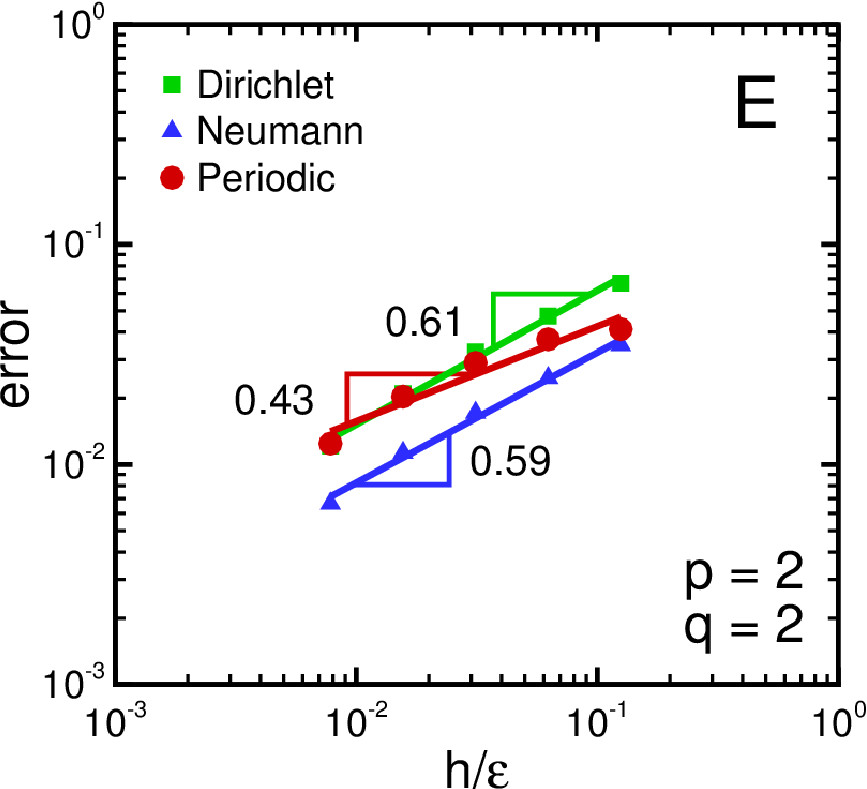}
	\caption{\textbf{Micro error convergence on the macroscale for chessboard stiffness pattern.} (first row) linear shape functions, $p$=$q$=$1$, (second row) quadratic shape functions, $p$=$q$=$2$, (from left to right:) $L^2$-, $H^1$-and energy-norm.}
	\label{fig:diagramm_chessboard_q1}
\end{Figure}

The results for linear shape functions are displayed in Fig. \ref{fig:diagramm_chessboard_q1} (first row). Again, the different coupling conditions agree well in the convergence order. The relative error however exhibits larger deviations between the coupling conditions. An exception is the energy-norm for PBC which leads to far worse results. Starting with the coarse discretizations on the right side of Fig. \ref{fig:diagramm_chessboard_q1} (first row, right) the error first increases with finer meshes and finally decreases again. The results of the calculations with rather fine micro meshes fit well into the results for Dirichlet and Neumann coupling, while the calculated errors for coarse meshes seem to be too small.

Again the optimal convergence order can not be reached in any of the norms due to the reduced regularity of the micro problem which is again based in the stiffness-jump at the interface between the two phases. The deviation from the optimal convergence order is even larger compared to the matrix-inclusion problem.

\textbf{Remark:} The convergence results for PBC in the energy-norm deserve a closer investigation; the analysis reveals that stresses in the macroscopic quadrature points do not exhibit sufficient accuracy. The entries of the homogenized elasticity tensor $\mathbb{A}^0$, which is used to calculate macroscopic stresses, converge with orders in the range from 0.62 ($\mathbb{A}_{12}^0$) to 0.72 ($\mathbb{A}_{11}^0$) which is in the range of the convergence orders of the $L^2$- and $H^1$-norm. The investigation of the homogenized elasticity tensor showed that there is a major absolute error in the single entries. The error of the coarsest discretization is about 50\% of the numerical values of the reference solution in all entries. For Neumann coupling in contrast, the errors are in the range from 7--12\%, and for Dirichlet coupling in the range of 8--21\%. These findings suggest that the questionable results for the error in the energy-norm for PBC is caused by the major absolute error of the homogenized tensor.  
 
Figure \ref{fig:diagramm_chessboard_q1} (second row) shows the results for quadratic shape functions. The results do not differ significantly from the results for linear shape functions. Again the error in the energy-norm seems to be too small for periodic coupling conditions and rather coarse discretizations.

The optimal convergence order of $2q = 4$ is clearly missed and the use of quadratic shape functions does not improve the convergence order.

\begin{Figure}[htbp]
	\centering
	\includegraphics[width=0.32\linewidth]{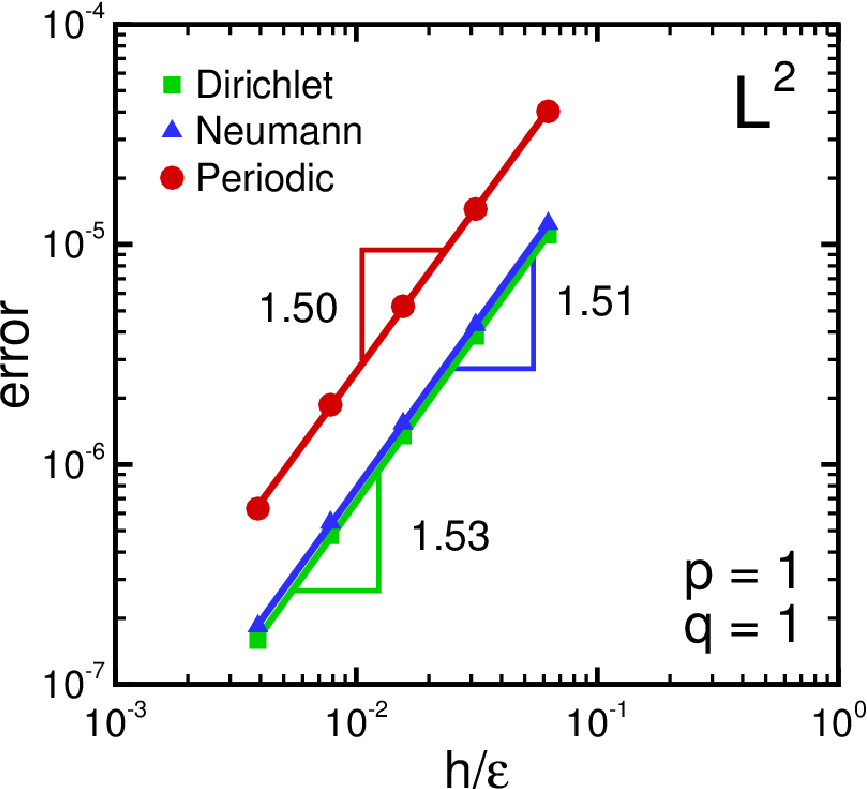} \hfill
	\includegraphics[width=0.32\linewidth]{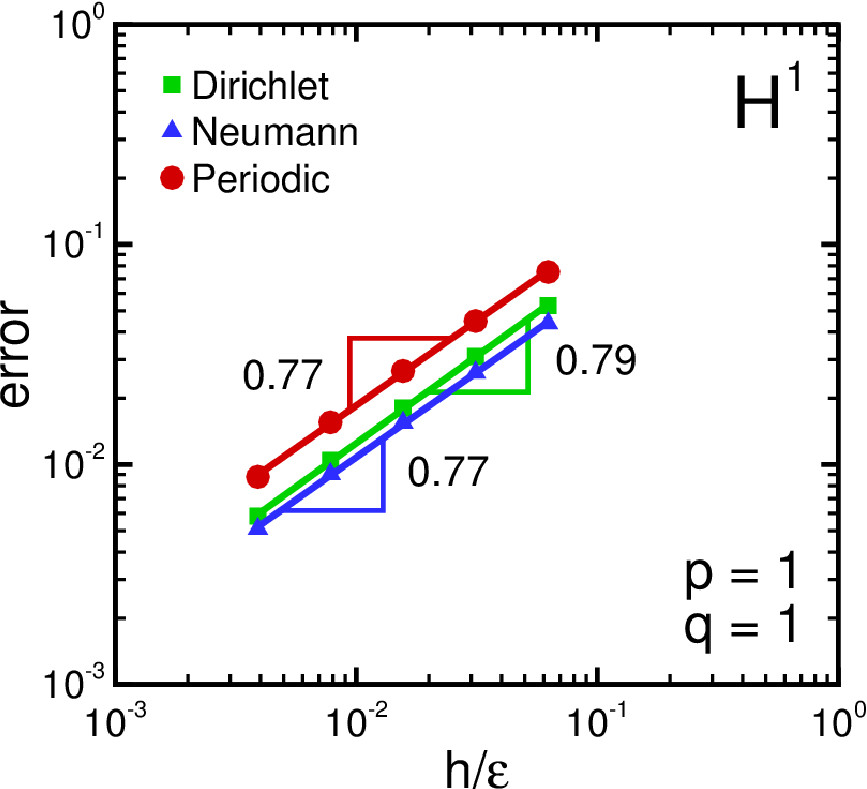} \hfill
	\includegraphics[width=0.32\linewidth]{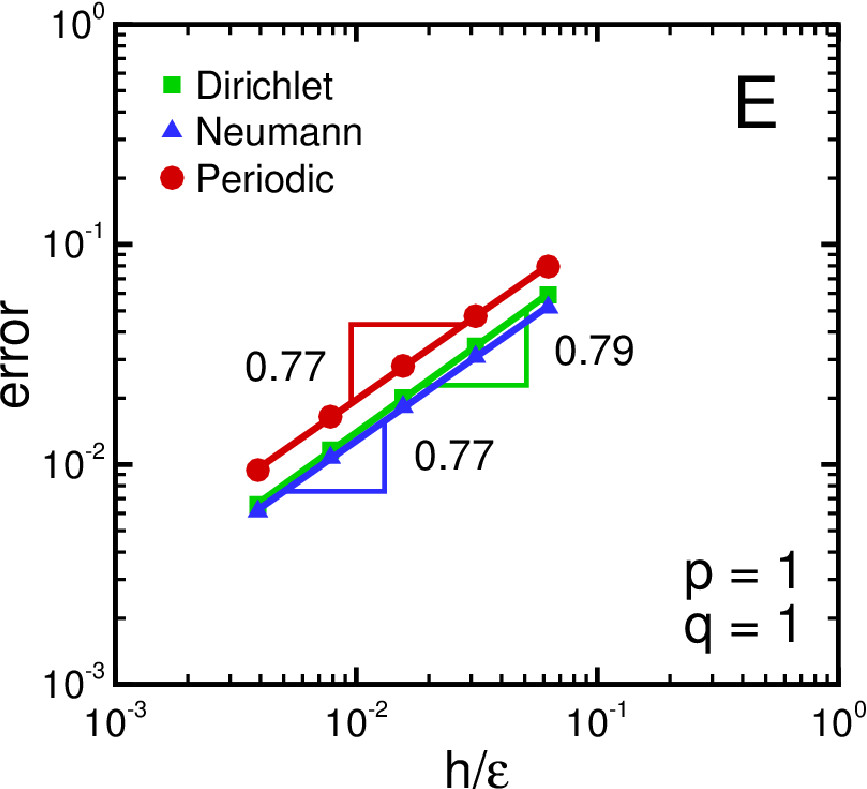} \\
	\includegraphics[width=0.32\linewidth]{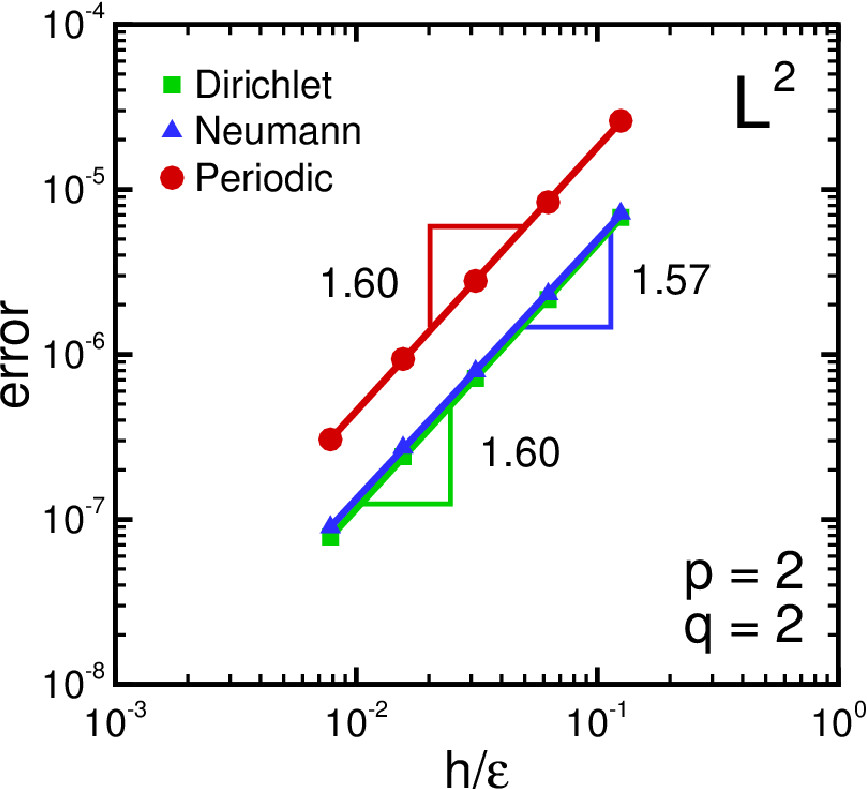} \hfill
	\includegraphics[width=0.32\linewidth]{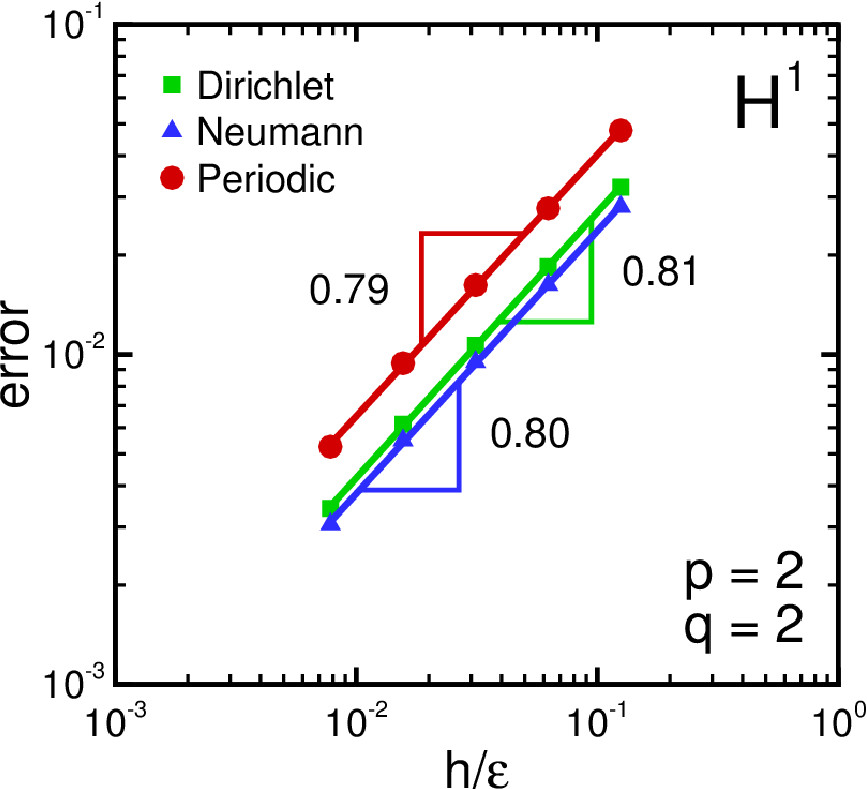} \hfill
	\includegraphics[width=0.32\linewidth]{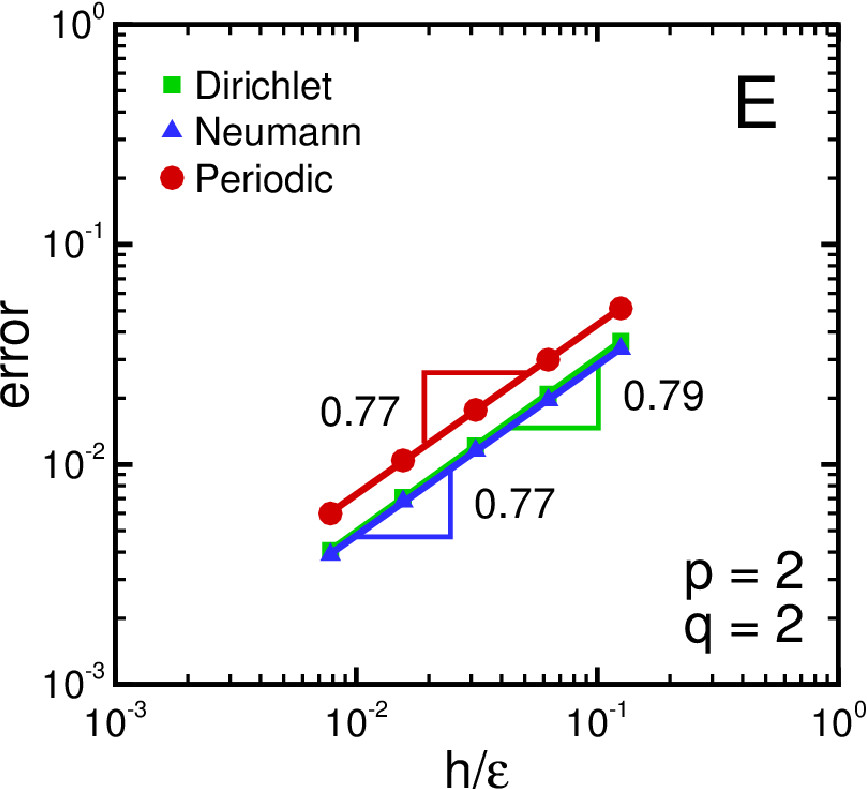}
	\caption{\textbf{Micro error convergence on the microscale for chessboard stiffness pattern.} (first row) linear shape functions, $p$=$q$=$1$, (second row) quadratic shape functions, $p$=$q$=$2$, (from left to right:) $L^2$-, $H^1$-and energy-norm.}
	\label{fig:diagramm_chessboard_q1_microdomain}
\end{Figure}

{\bf Micro error convergence on the microscale.} The results of an error calculation on the microdomain related to the macroscopic quadrature point at [0.26, 0.26] with linear shape functions can be found in Fig. \ref{fig:diagramm_chessboard_q1_microdomain} (first row). The convergence orders for different coupling conditions exhibit good agreement. 

Figure \ref{fig:diagramm_chessboard_q1_microdomain} (second row) indicates that for quadratic shape functions the convergence orders are in good agreement for the different coupling conditions. The convergence order however is not improved for quadratic shape functions compared to the linear case, which is due to the low regularity.

\subsubsection{Sine wave distribution}

The low regularity of the micro BVP in the first two examples is the reason why convergence for quadratic shape functions shows a strong deviation from the nominal order. Aiming at the full convergence order of $2q$ for the micro error a sine wave-type Young's modulus distribution is chosen, which is expected to exhibit high regularity for its smooth stiffness distribution. Therein, the minimum Young's modulus is $E_{min} = 40000$ $[F/L^2]$, the maximum is $E_{max} = 50000$ $[F/L^2]$.

\begin{Figure}[htbp]
	\centering
	\includegraphics[width=0.34\linewidth]{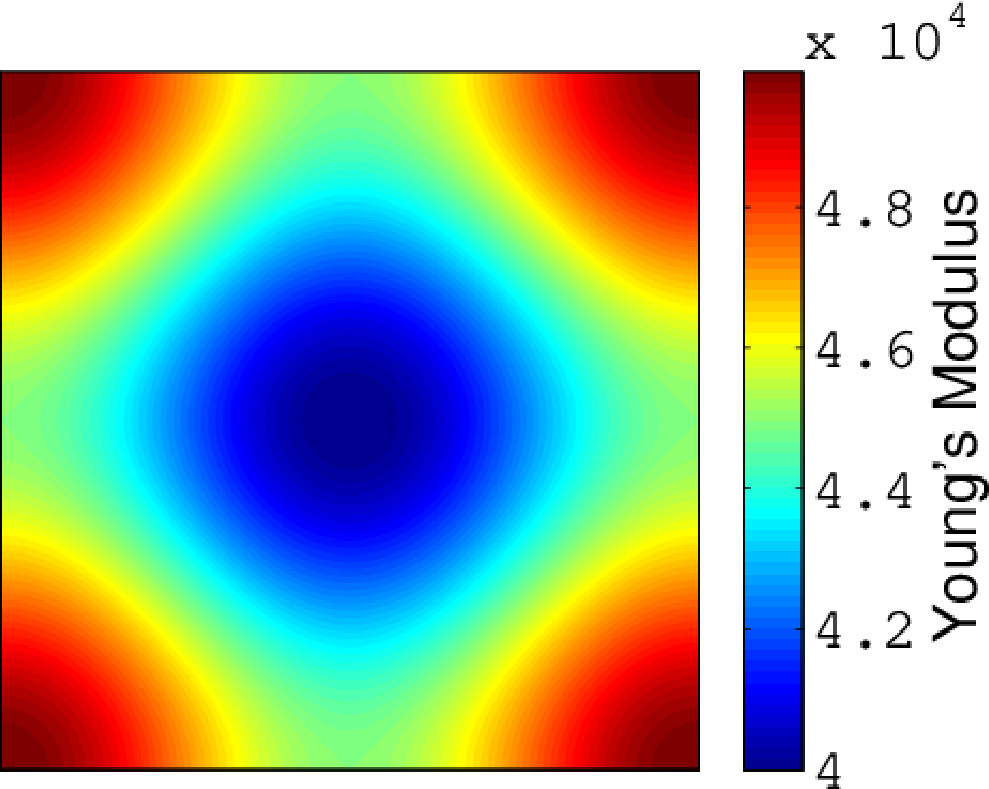} \hspace*{16mm}
	\includegraphics[width=0.34\linewidth]{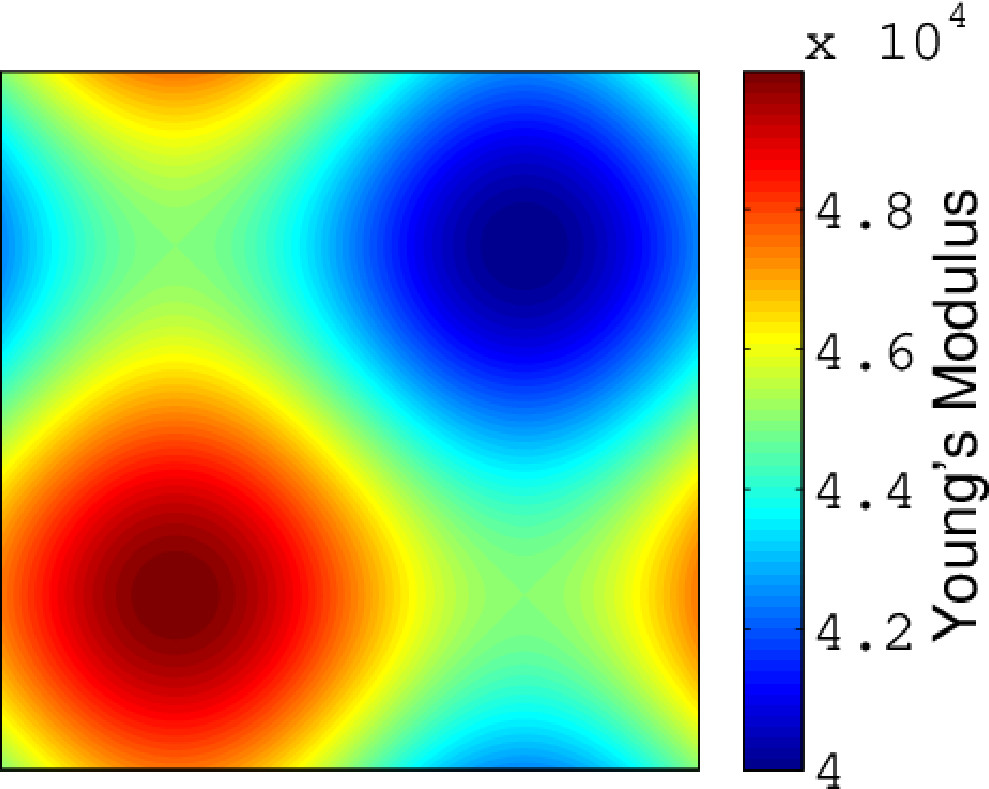}  
	\caption{\textbf{Sine wave distribution problem.} Distribution of Young's modulus on (left) a unit cell with cubic symmetry and (right) an alternative unit cell.}
	\label{fig:sinewave}
\end{Figure}

The Young's modulus distribution on the micro domain is depicted in Fig.~\ref{fig:sinewave}. The unit cell in the left reflects the cubic symmetry of the periodic structure; an alternative definition (among many others) of the unit cell is displayed on the right of Fig.~\ref{fig:sinewave}. While the stiffness results for Neumann and Dirichlet coupling depend on the choice of the unit cell, for PBC stiffness is invariant with respect to that choice. In the following we use the unit cell in the right of Fig.~\ref{fig:sinewave}. 

\begin{Figure}[htbp]
	\centering
	\includegraphics[width=0.32\linewidth]{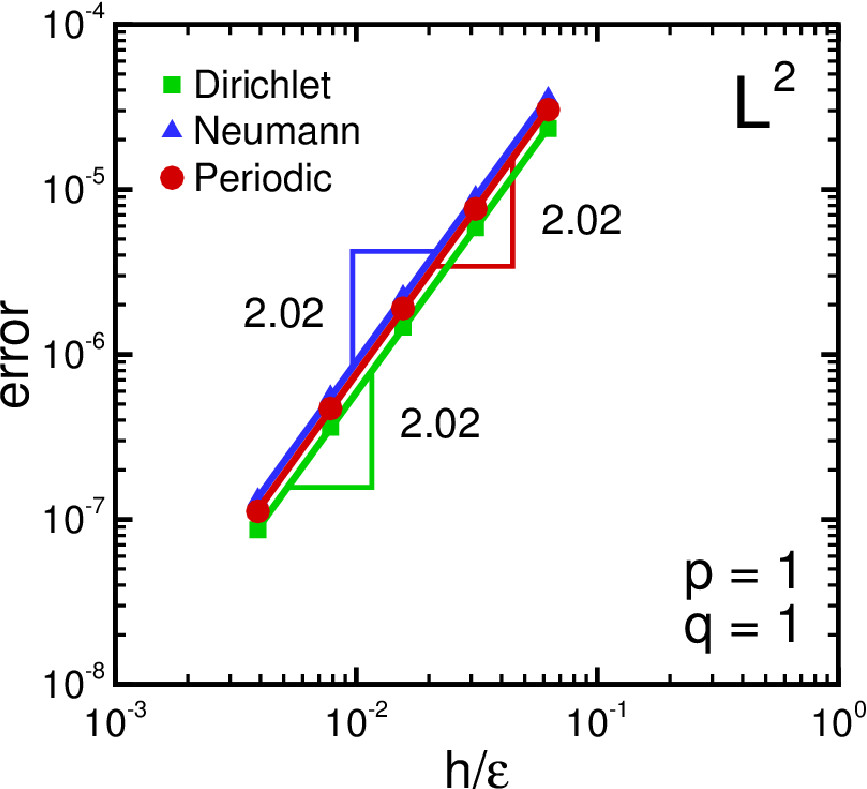} \hfill
	\includegraphics[width=0.32\linewidth]{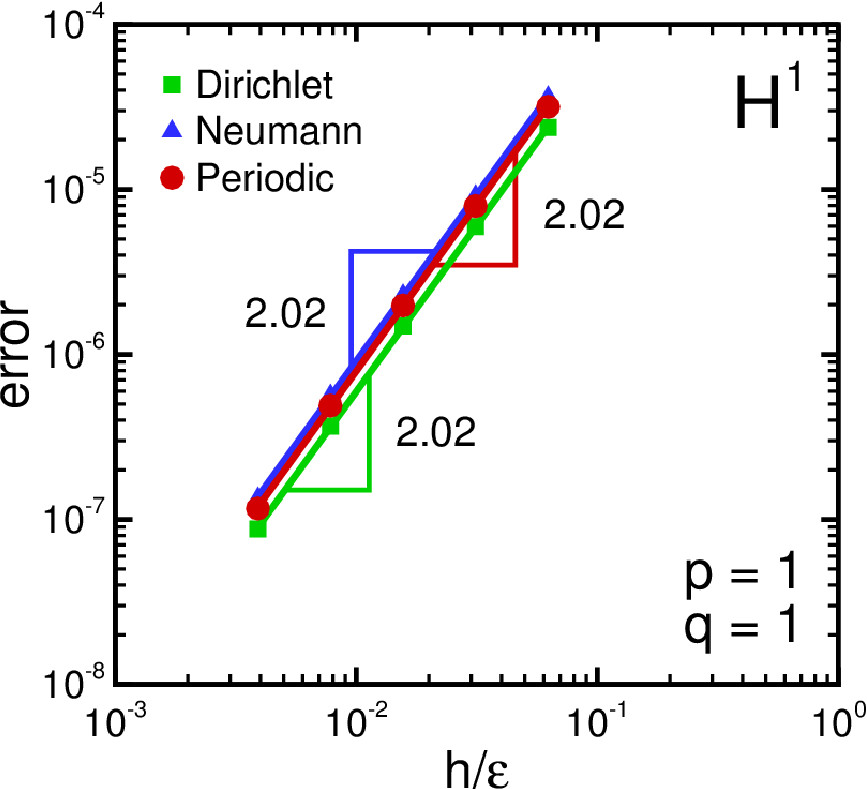} \hfill
	\includegraphics[width=0.32\linewidth]{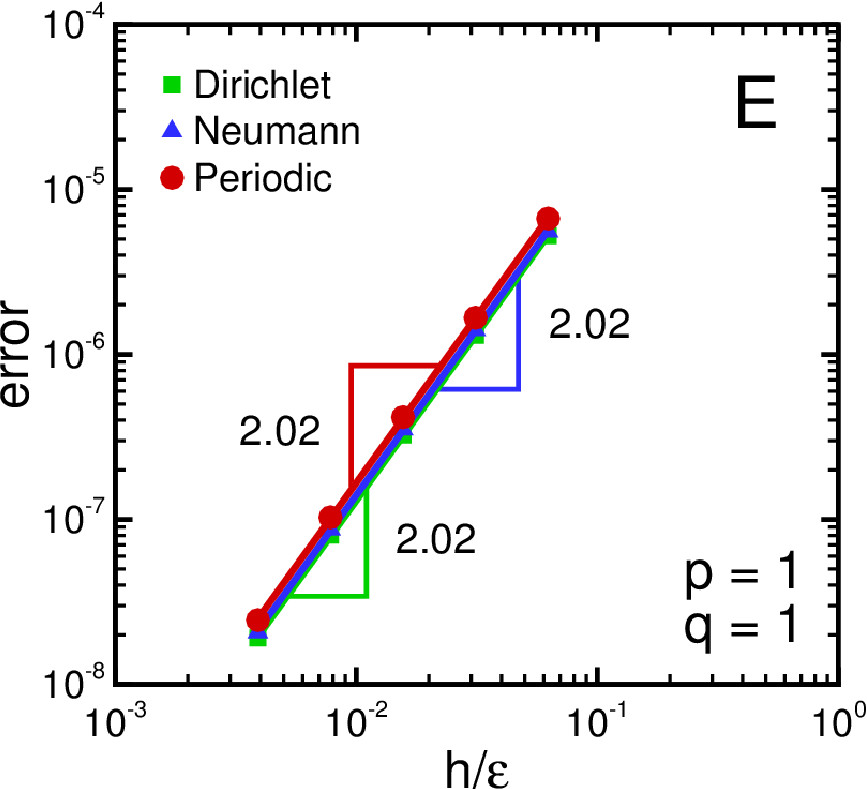} \\
	\includegraphics[width=0.32\linewidth]{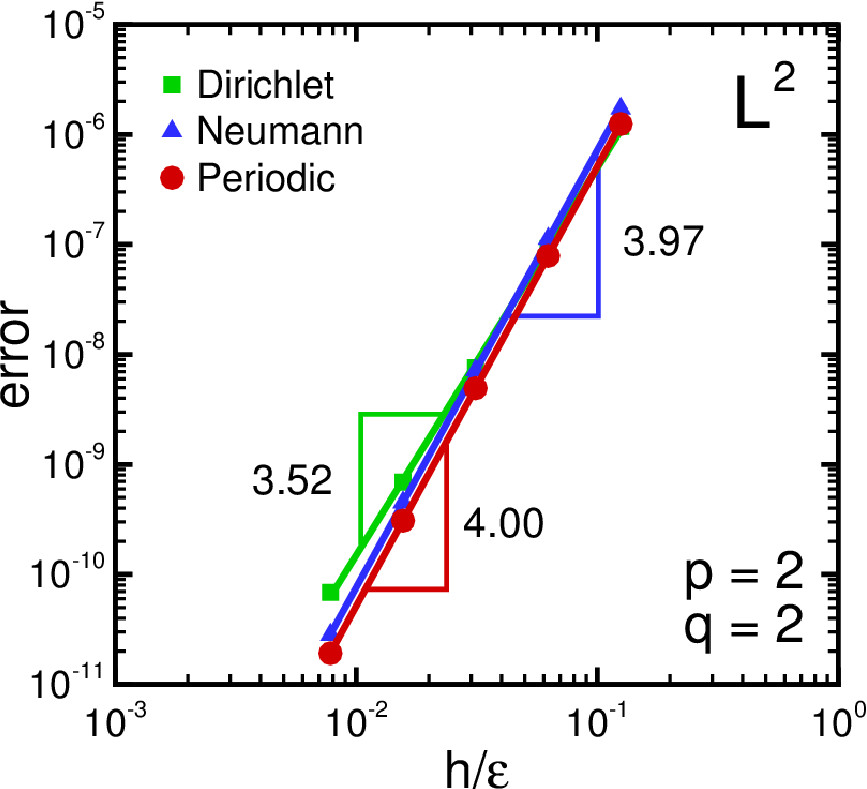} \hfill
	\includegraphics[width=0.32\linewidth]{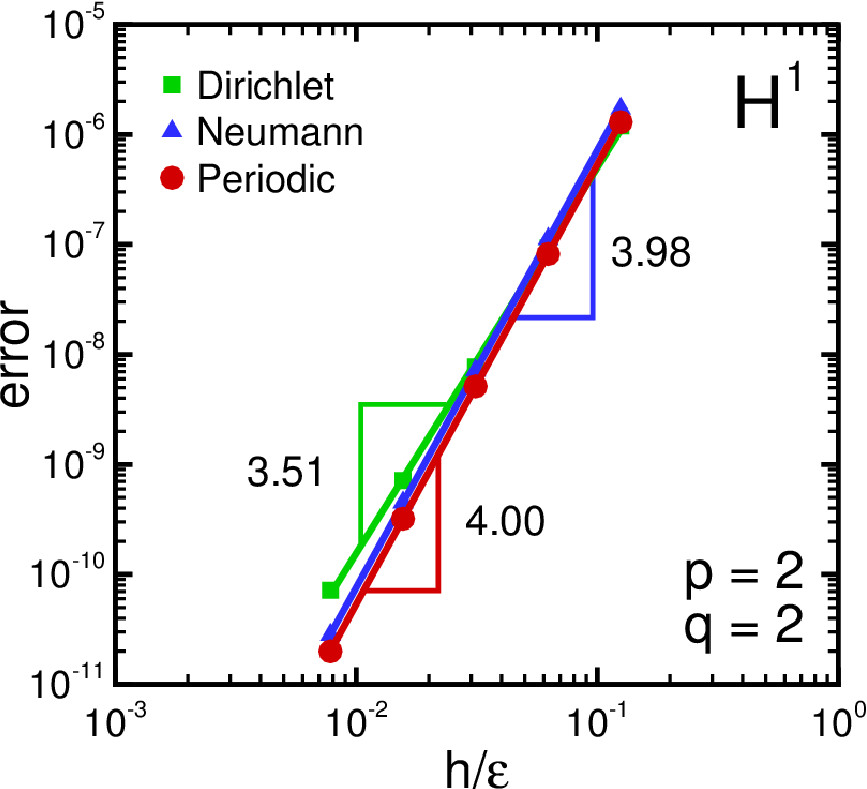} \hfill
	\includegraphics[width=0.32\linewidth]{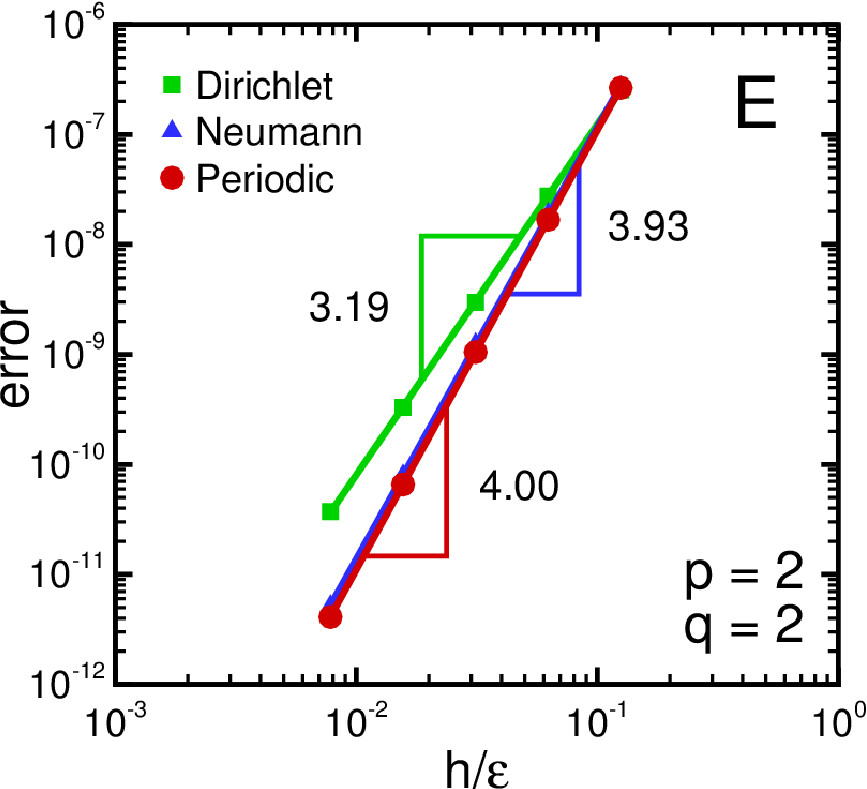}
	\caption{\textbf{Micro error convergence on the macroscale for a microstructure with sine wave stiffness distribution.} (first row) linear shape functions, $p$=$q$=$1$, (second row) quadratic shape functions, $p$=$q$=$2$, (from left to right:) $L^2$-, $H^1$-and energy-norm.}
	\label{fig:diagramm_sinewave_q1}
\end{Figure}

Figure \ref{fig:diagramm_sinewave_q1} (first row) shows the results for linear shape functions. In all norms and for all coupling conditions the convergence order of the calculated errors is 2.02. The numerical values of the relative errors also show only minor deviations. The sine wave distribution enables full regularity of the solution as indicated by the full theoretical convergence order in all norms.
 
The results for quadratic micro shape functions are displayed in the second row of Fig.~\ref{fig:diagramm_sinewave_q1}. For periodic and Neumann coupling conditions the optimal convergence order of $2q = 4$ is virtually achieved in all norms, while for Dirichlet coupling a reduced order is observed. The numerical values of the relative errors of Neumann and periodic coupling are in good agreement, while the values for Dirichlet coupling exhibit good agreement with the estimates only for coarse discretizations, for finer discretizations they worsen most notably in the energy-norm.

In conclusion, the regularity of the micro BVP enables full convergence order and --opposed to the first two examples-- a higher convergence order for quadratic shape functions than for linear shape functions. 

\begin{Figure}[htbp]
	\centering
	\includegraphics[width=0.32\linewidth]{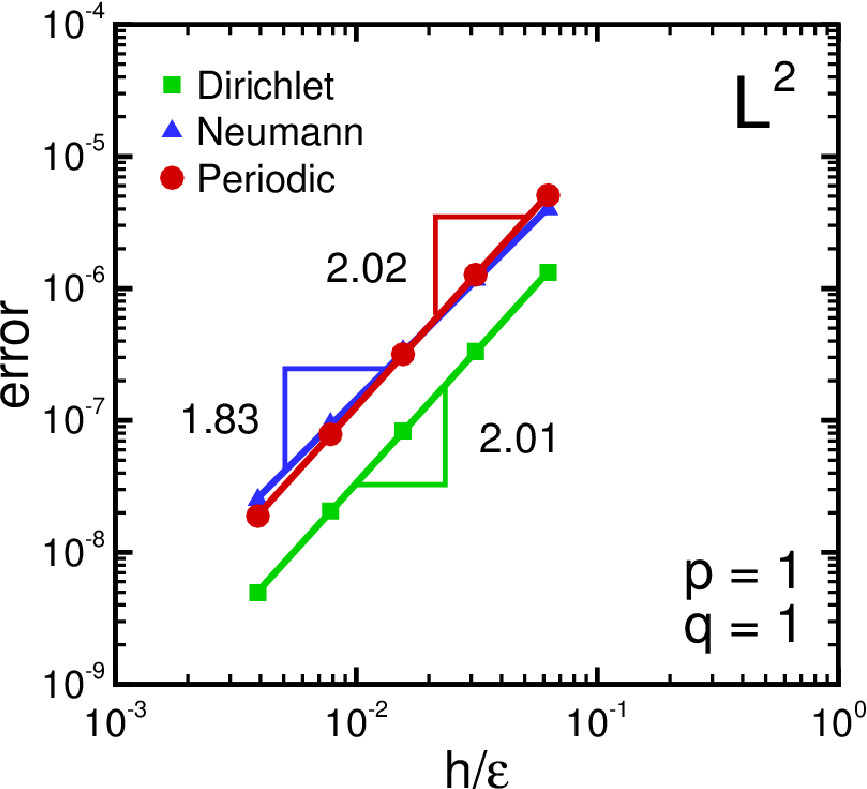} \hfill
	\includegraphics[width=0.32\linewidth]{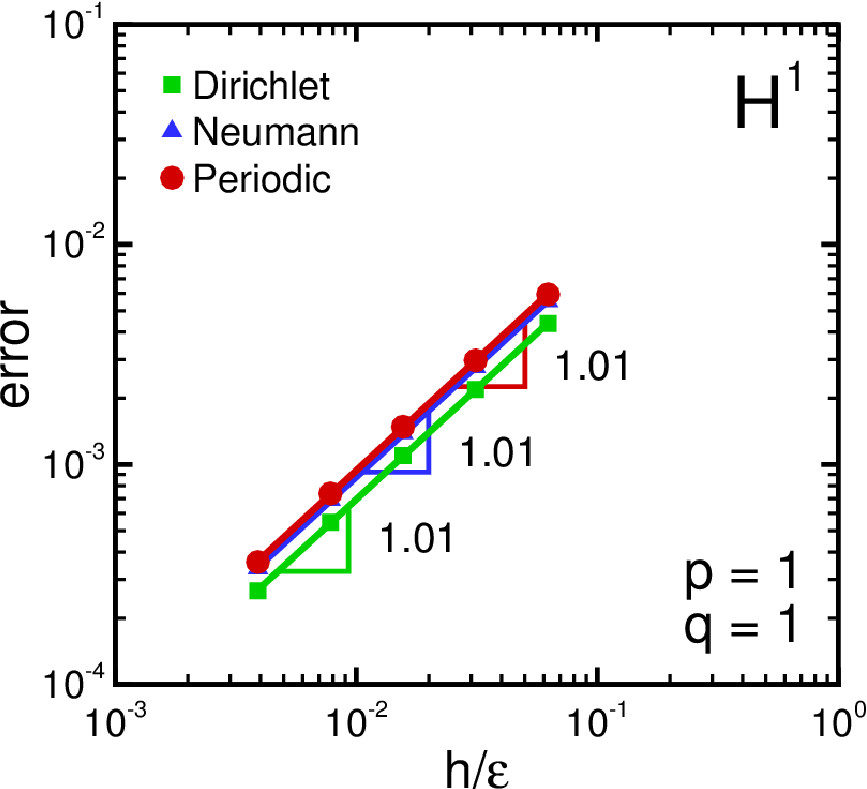} \hfill
	\includegraphics[width=0.32\linewidth]{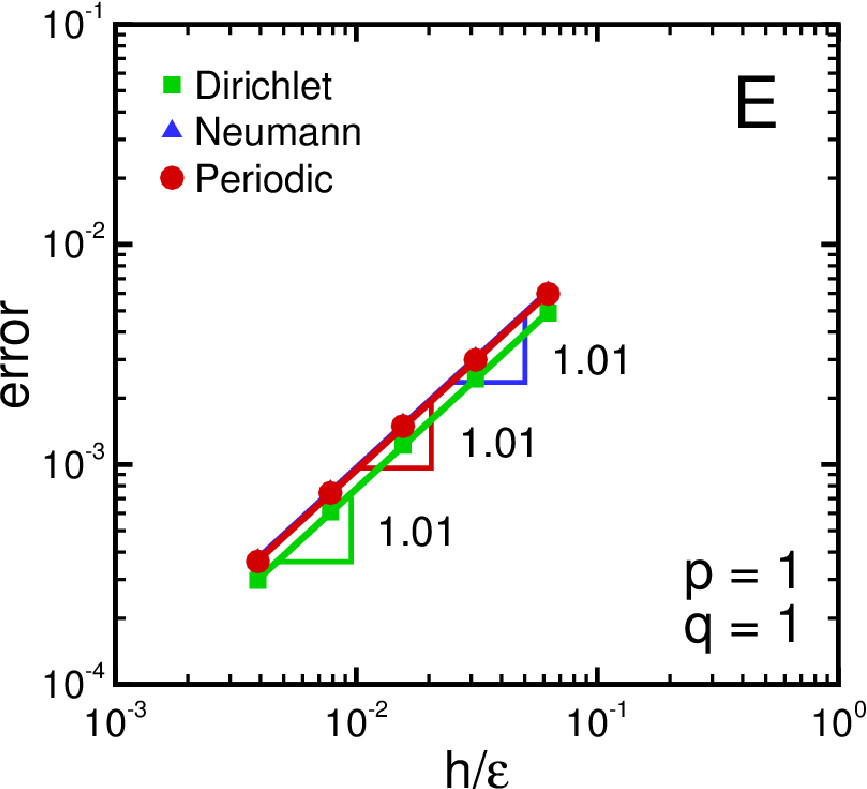} \\
	\includegraphics[width=0.32\linewidth]{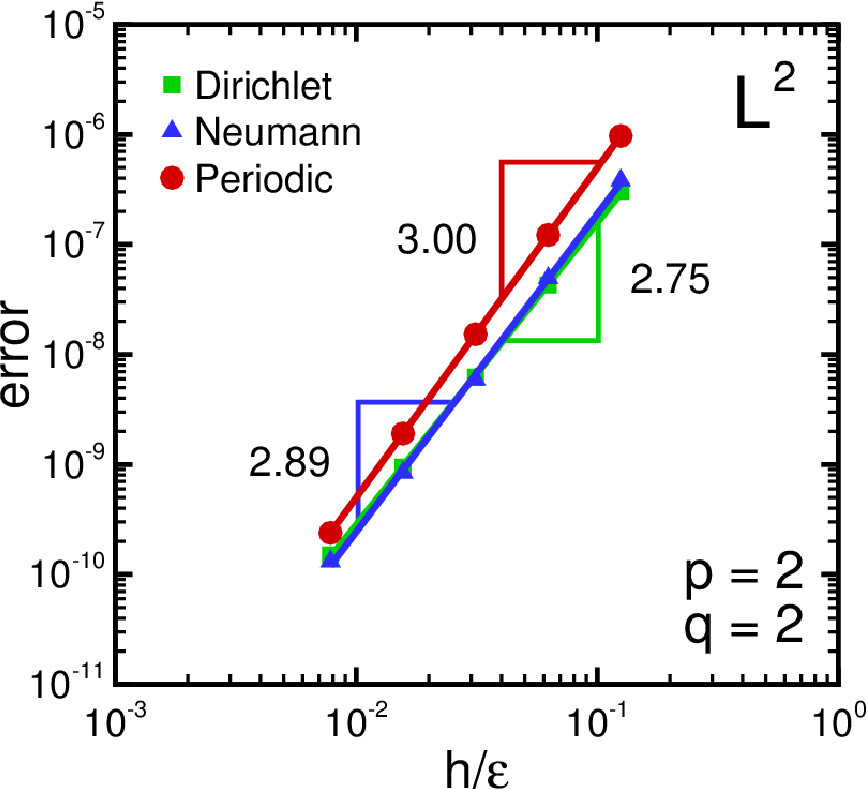} \hfill
	\includegraphics[width=0.32\linewidth]{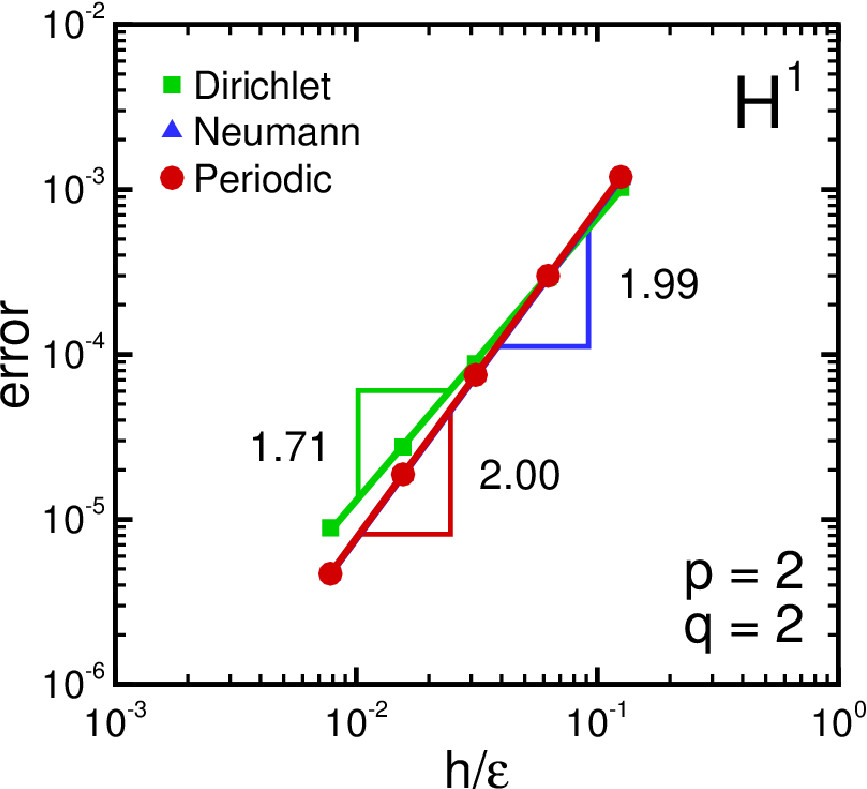} \hfill
	\includegraphics[width=0.32\linewidth]{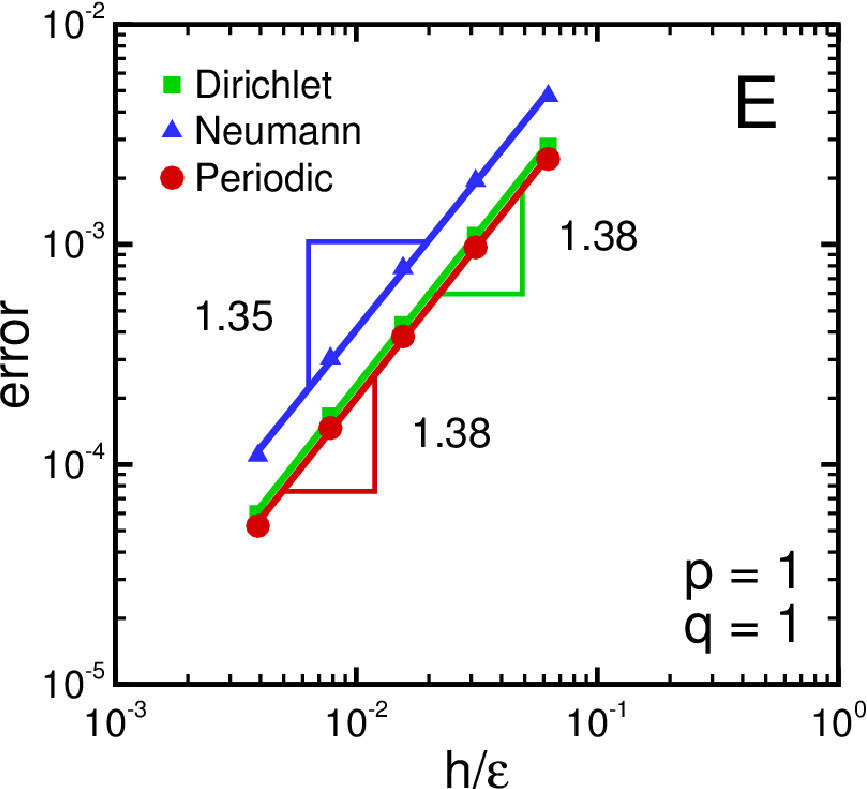}
        \caption{\textbf{Micro error convergence on the microscale for a microstructure with sine wave stiffness distribution.} (first row) linear shape functions $p$=$q$=$1$, (second row) quadratic shape functions $p$=$q$=$2$, (from left to right:) $L^2$-, $H^1$-, and energy-norm.}  
	\label{fig:diagramm_sinewave_q1_microdomain}
\end{Figure}

{\bf Micro error convergence on the microscale.} The error calculation on the microdomain related to the macroscopic quadrature point at [0.26, 0.26] for $q=1$ yields the results that are displayed in Fig. \ref{fig:diagramm_sinewave_q1_microdomain} (first row). The convergence orders of $q+1$ in the $L^2$-norm and $q$ in the $H^1$- and energy-norm are achieved for Dirichlet and periodic coupling, for Neumann coupling there are some minor deviations in the $L^2$-norm. 

The results of the error calculation on the same microdomain for quadratic shape functions are displayed in Fig. \ref{fig:diagramm_sinewave_q1_microdomain} (second row). For periodic coupling the optimal convergence order is achieved in all norms, for Neumann coupling the optimal convergence order is restricted to the $H^1$- and the energy-norm. Dirichlet coupling however, shows again minor reductions in all three norms, which is consistent with the reduced convergence order of the micro error on the macroscale.

\textbf{Remark:} The measured convergence orders being almost in perfect agreement with the a priori estimates could suggest that the observed regularity is due to the low stiffness contrast (1:1.25).  Additional analyses employing an increased stiffness contrast of up to 1:25 yield the same convergence orders and thereby rebut this hypothesis. Instead it is the smoothness of Young's modulus distribution that enables the regularity in terms of full convergence orders.

\subsection{Macro convergence analysis}

After the assessment of micro errors both on the micro as well as on the macroscale, the macro error convergence is investigated in the following. 

\subsubsection{Square cantilever}
\label{subsubsec:SquareCantilever}
 
%\begin{Figure}[htbp]
%	\centering
%	\psfrag{l}{1}
%	\psfrag{t}{0.1}
%	\includegraphics[height=5.0cm]{model_standard_beam}
%	\caption{Model of a square cantilever.}
%	\label{fig:model_standard_beam}
%\end{Figure}

\begin{Figure}[htbp]
	\centering
	\includegraphics[height=5.0cm]{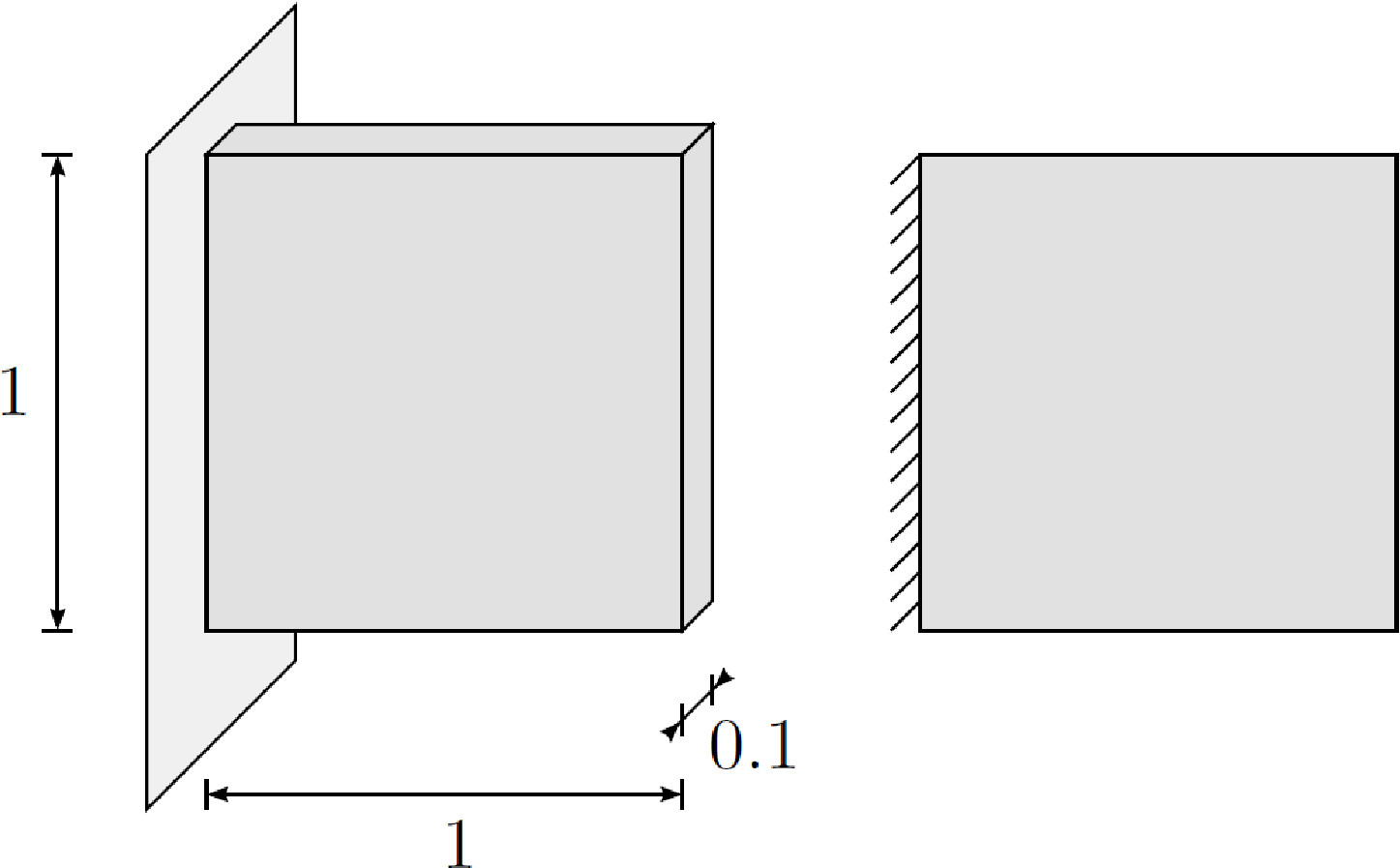}
	\caption{Model of a square cantilever.}
	\label{fig:model_standard_beam}
\end{Figure}

In the first numerical example we consider the square cantilever of Fig.~\ref{fig:model_standard_beam}, which is subject to volume forces of $\bm f = [0, -10]^T$ $[F/L^2]$. The microstructure is the already introduced sine wave-type distribution of Young's modulus.

\begin{Figure}[htbp]
	\centering
	\includegraphics[width=0.32\linewidth]{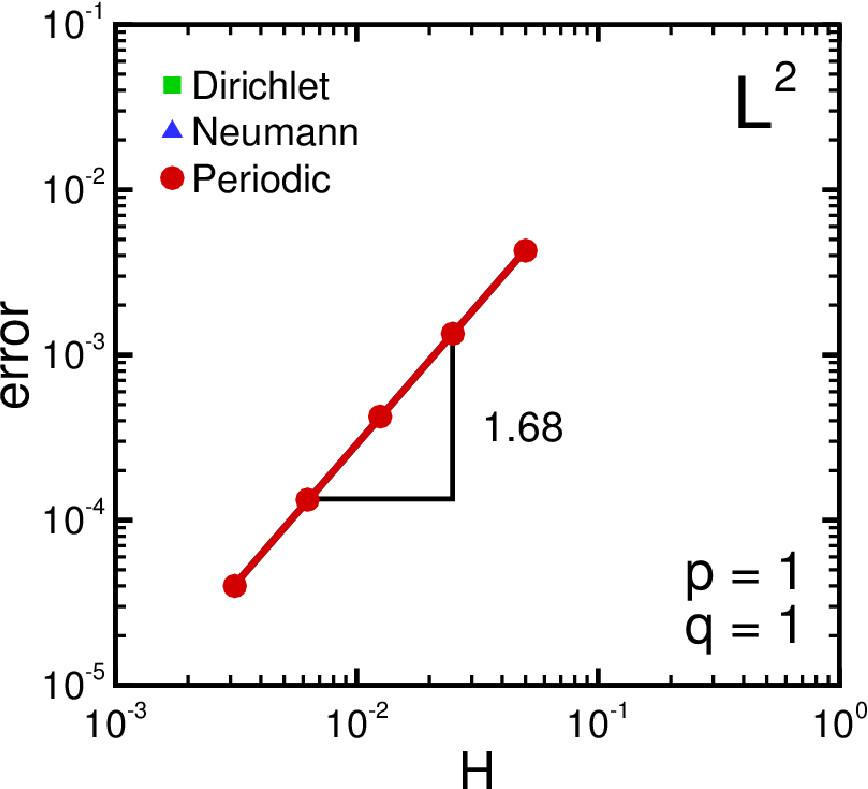} \hfill 
	\includegraphics[width=0.32\linewidth]{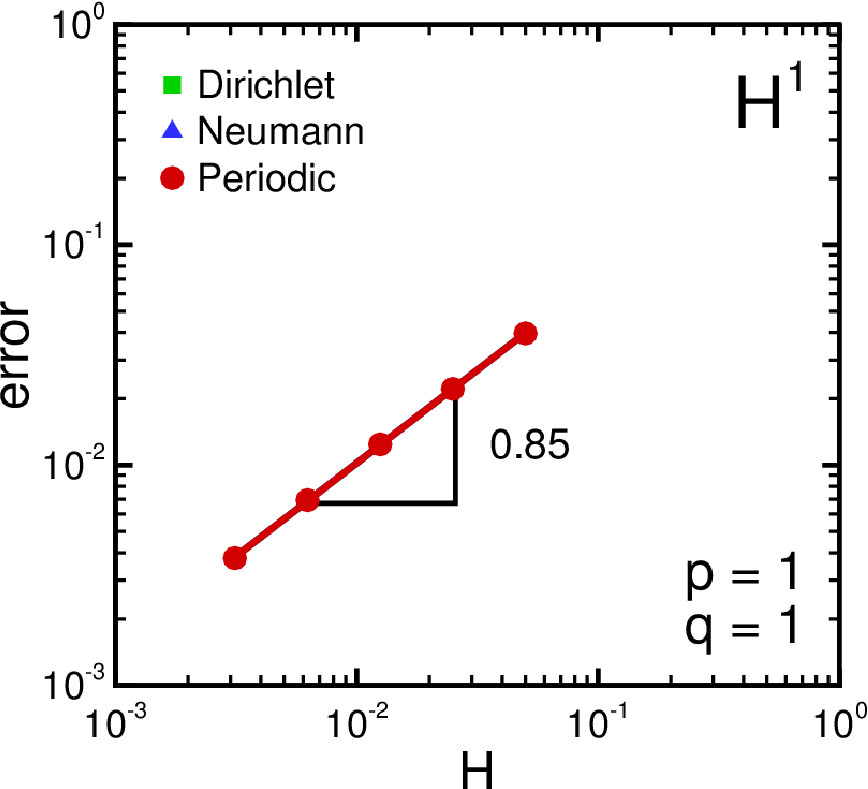} \hfill 
	\includegraphics[width=0.32\linewidth]{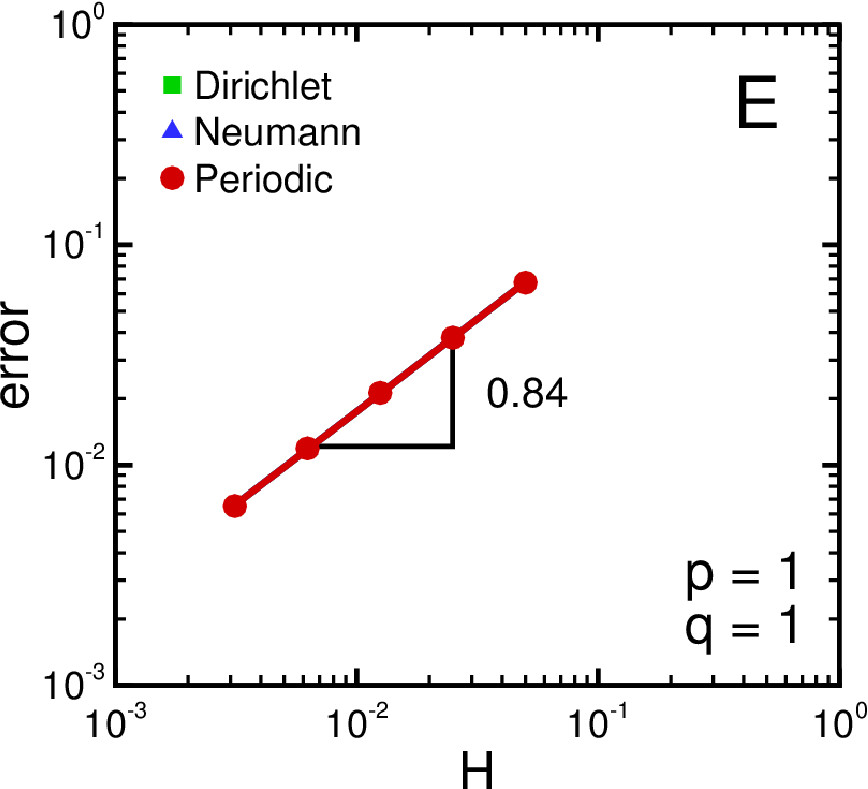} \\
	\includegraphics[width=0.32\linewidth]{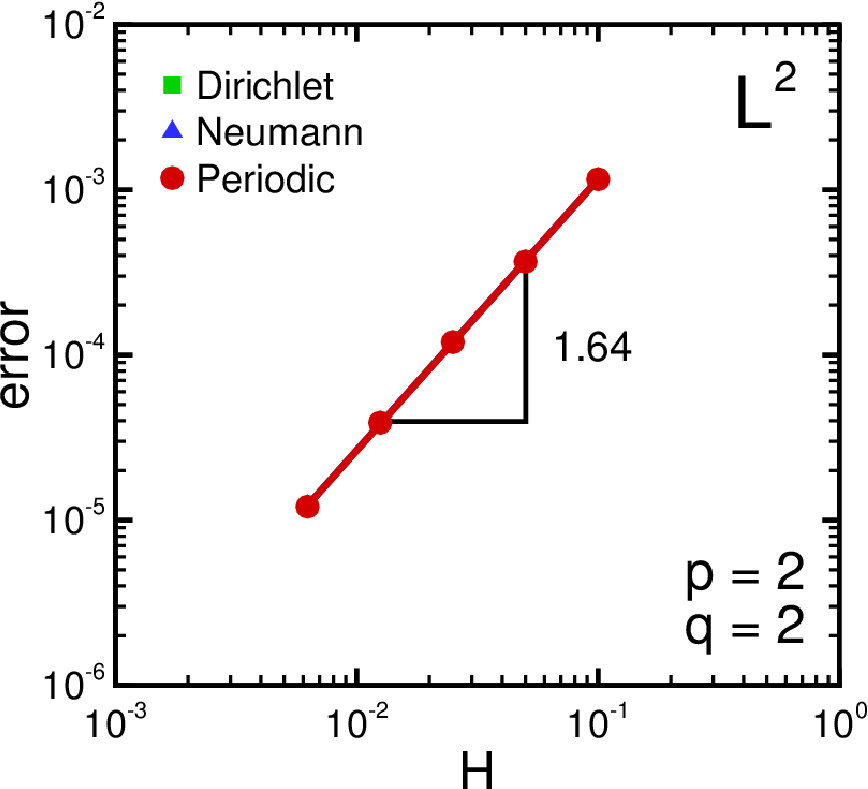} \hfill 
	\includegraphics[width=0.32\linewidth]{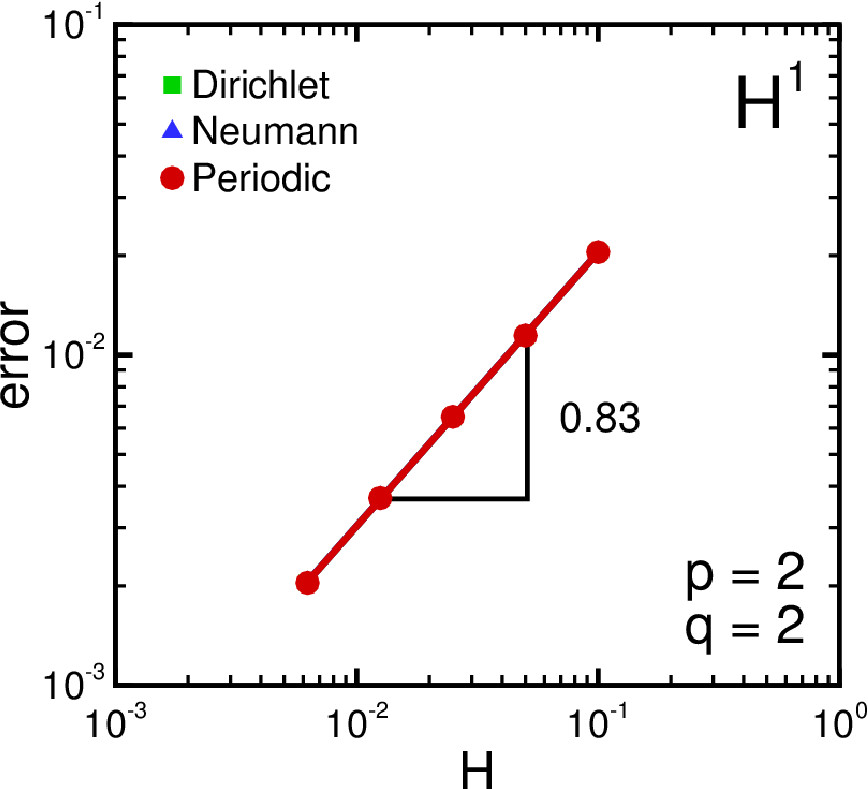} \hfill 
	\includegraphics[width=0.32\linewidth]{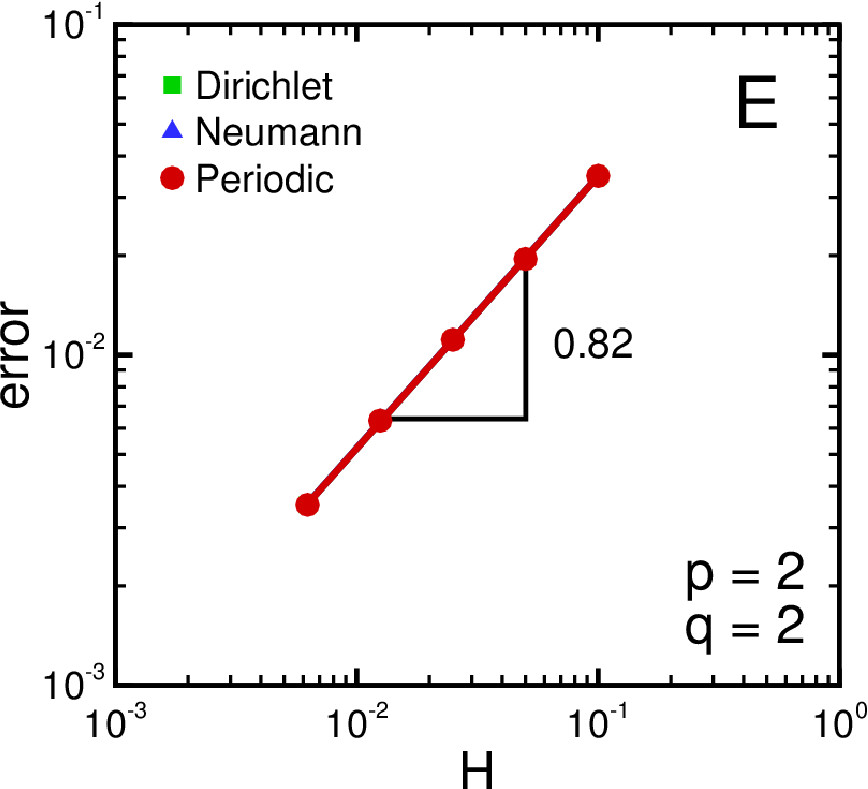}
	\caption{\textbf{Macroconvergence for square cantilever under volume force.} (first row) linear shape functions $p$=$q$=$1$, (second row)
	quadratic shape functions $p$=$q$=$2$, (from left to right) $L^2$-, $H^1$-, and energy-norm.}
	\label{fig:diagramm_cantilever_beam_p1}
\end{Figure}

The results of the convergence analysis for linear shape functions is shown in the first row of Fig.~\ref{fig:diagramm_cantilever_beam_p1}. The optimal convergence orders, of $p+1$ in the $L^2$-norm and of $p$ in the $H^1$- and energy-norm are not achieved, but the deviations are small. The micro-macro coupling condition has virtually no influence on the macro error convergence. The results employing quadratic shape functions as displayed in Fig.~\ref{fig:diagramm_cantilever_beam_p1} (second row) exhibit only very minor deviations compared to linear shape functions; first, the results for all coupling conditions coincide in each norm. Second, the convergence orders are 1.64 in the $L^2$-norm, 0.83 in the $H^1$-norm and 0.82 in the energy-norm, respectively. In conclusion, the optimal convergence orders are not obtained, the problem of order reduction is not cured at all by increasing the polynomial order of the shape functions.

The order reduction is caused by corner singularities at the clamped end of the square cantilever.

\begin{Figure}[htbp]
	\centering
	\includegraphics[height=4.0cm]{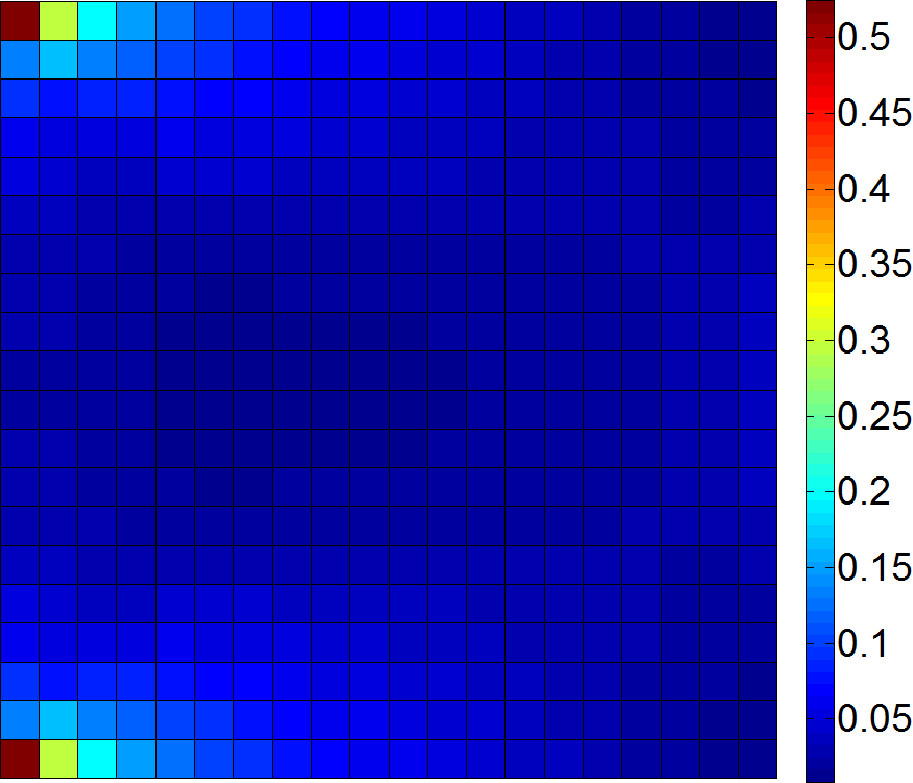} \hspace*{16mm} 
	\includegraphics[height=4.0cm]{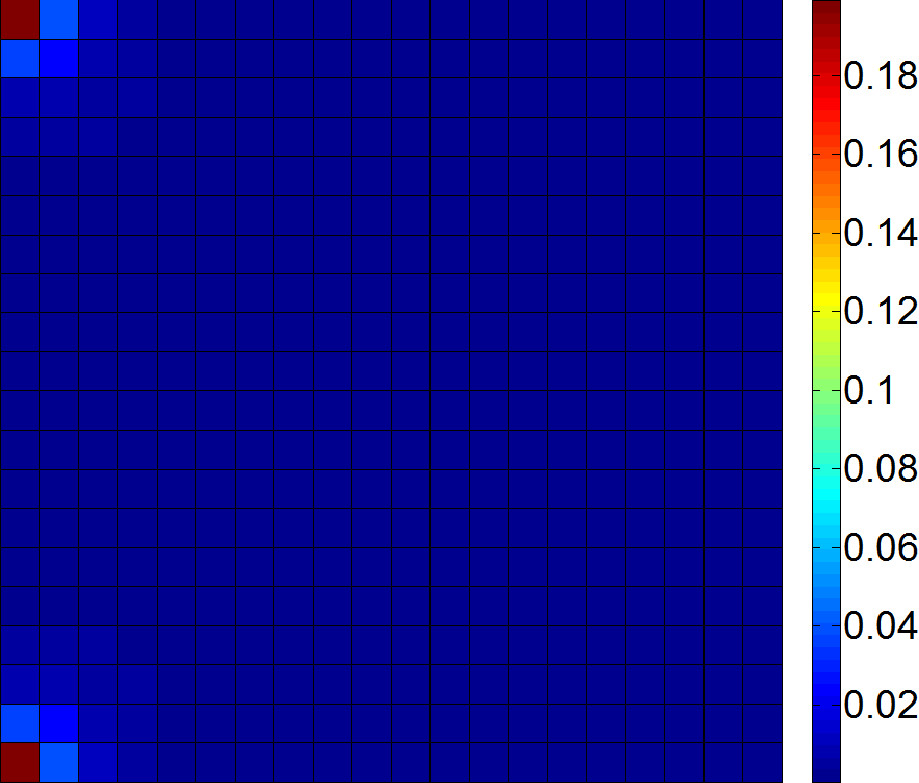}
	\caption{\textbf{Square cantilever under volume force.} Relative elementwise energy-error distribution on the macrodomain for (left) $p$=$q$=$1$ and (right) $p$=$q$=$2$.}
	\label{fig:error_distribution_cantilever_beam}
\end{Figure}

The relative elementwise error on the macrodomain is shown in Fig. \ref{fig:error_distribution_cantilever_beam}. The error is computed by the ratio of the error in each element and the average energy per element. For the visualization of the error distribution the errors were calculated in the single elements of a coarse macro mesh. For both polynomial orders of shape functions the maximum absolute error is located in the corners of the cantilever's bearing. The high relative error of more than 50\% for linear and more than 18\% for quadratic shape functions indicates that the total energy in these elements is considerably higher than the average energy per element to which it is related here.

\subsubsection{Tapered cantilever}
\label{subsubsec:TaperedCantilever}

%\begin{Figure}[htbp]
%	\centering
%	\psfrag{h}{2}
%	\psfrag{l}{1}
%	\psfrag{t}{0.1}
%	\psfrag{w}{$\alpha$}
%	\includegraphics[height=6.0cm]{model_tapered_beam}
%	\caption{Model of a tapered cantilever.}
%	\label{fig:model_tapered_beam}
%\end{Figure}

\begin{Figure}[htbp]
	\centering
	    \includegraphics[height=6.0cm]{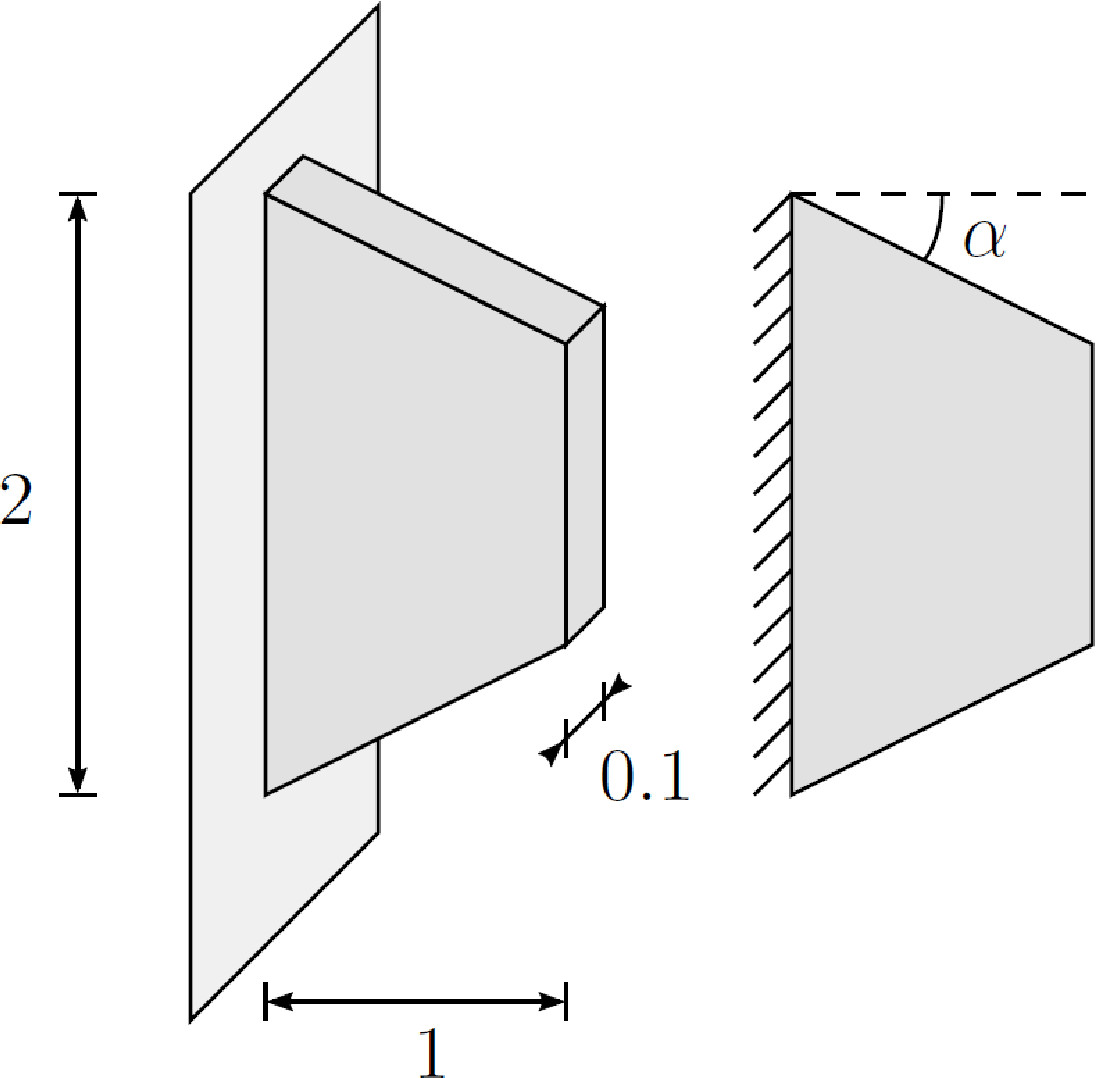}
	\caption{Model of a tapered cantilever.}
	\label{fig:model_tapered_beam}
\end{Figure}

If the angle $\alpha$ is chosen sufficiently large ($\alpha>28.4^{\circ}$), then the design of a tapered cantilever as in Fig.~\ref{fig:model_tapered_beam} avoids the singularities of the square cantilever plate, which was proven by analytical means in \cite{Roessle-2000}. Here we choose $\alpha=30.4^{\circ}$. Again, a volume load\footnote{According to \cite{Roessle-2000} stress-free boundaries on $\partial \mathcal{B}_{N}$ along with volume loads are an additional condition for the regularity of the BVP.} of $\bm f = [0, -10]^T$ $[F/L^2]$ is applied to the cantilever and the sine wave-type microstructure with a Young's modulus contrast of 1.25 is chosen.

\begin{Figure}[htbp]
	\centering
	\includegraphics[width=0.32\linewidth]{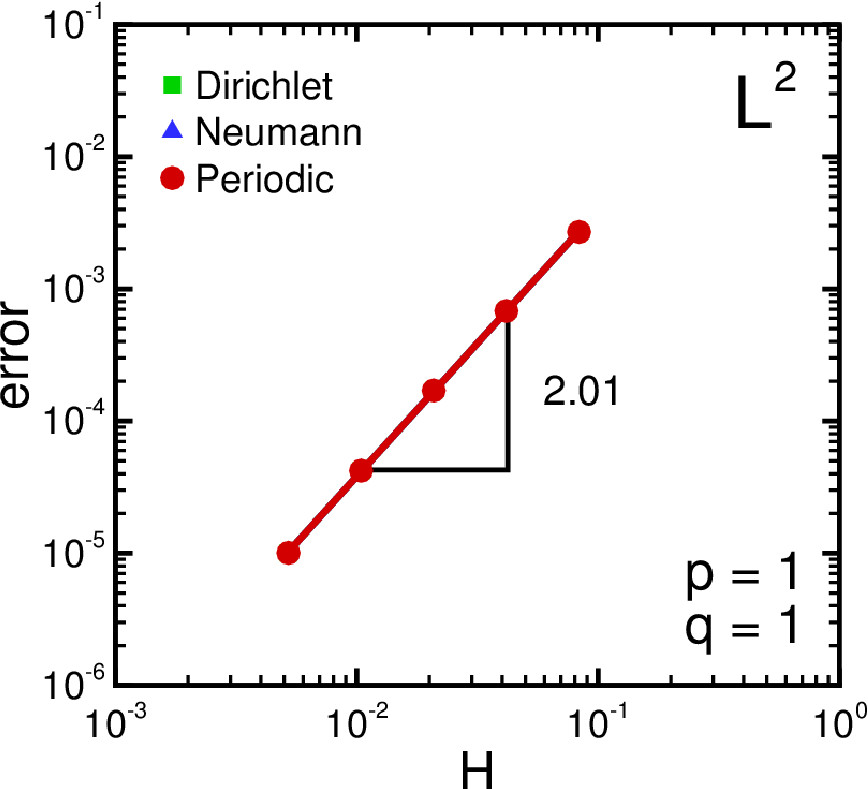} \hfill 
	\includegraphics[width=0.32\linewidth]{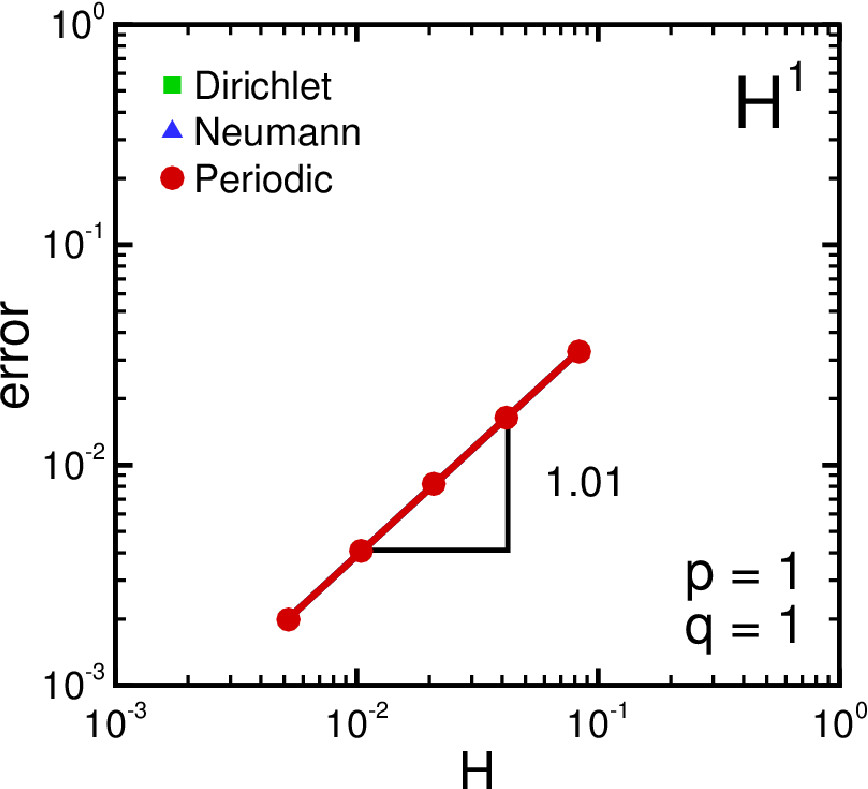} \hfill 
	\includegraphics[width=0.32\linewidth]{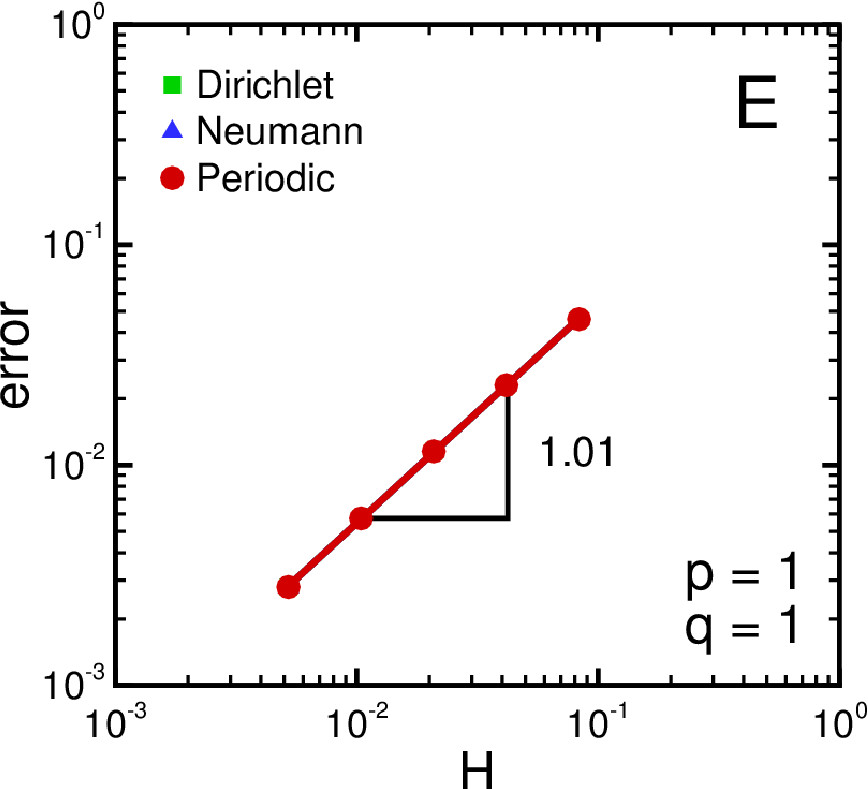} \\
	\includegraphics[width=0.32\linewidth]{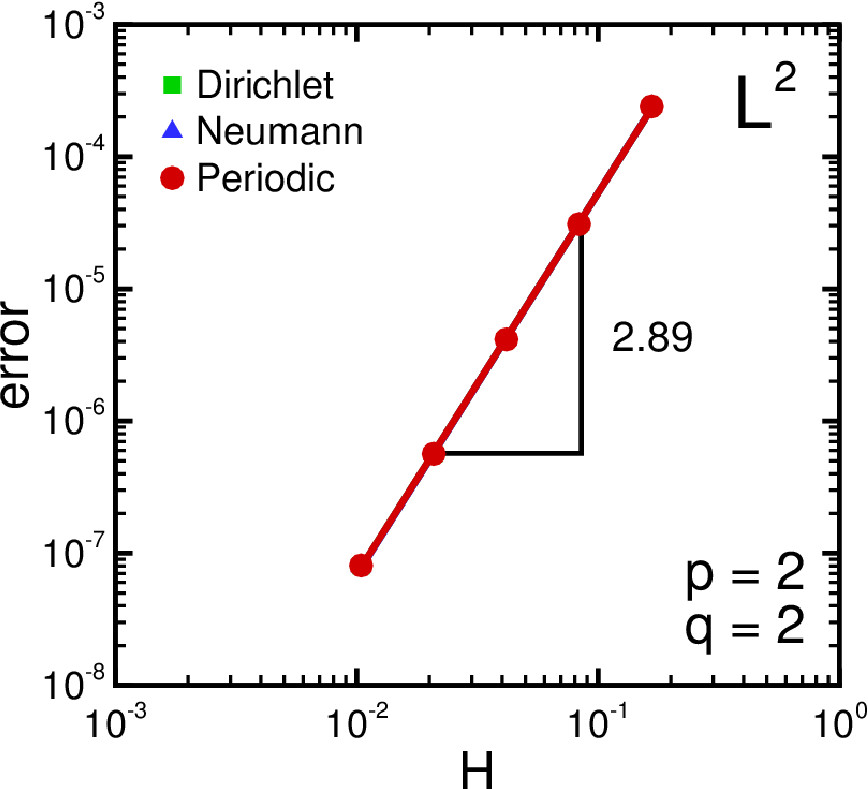} \hfill 
	\includegraphics[width=0.32\linewidth]{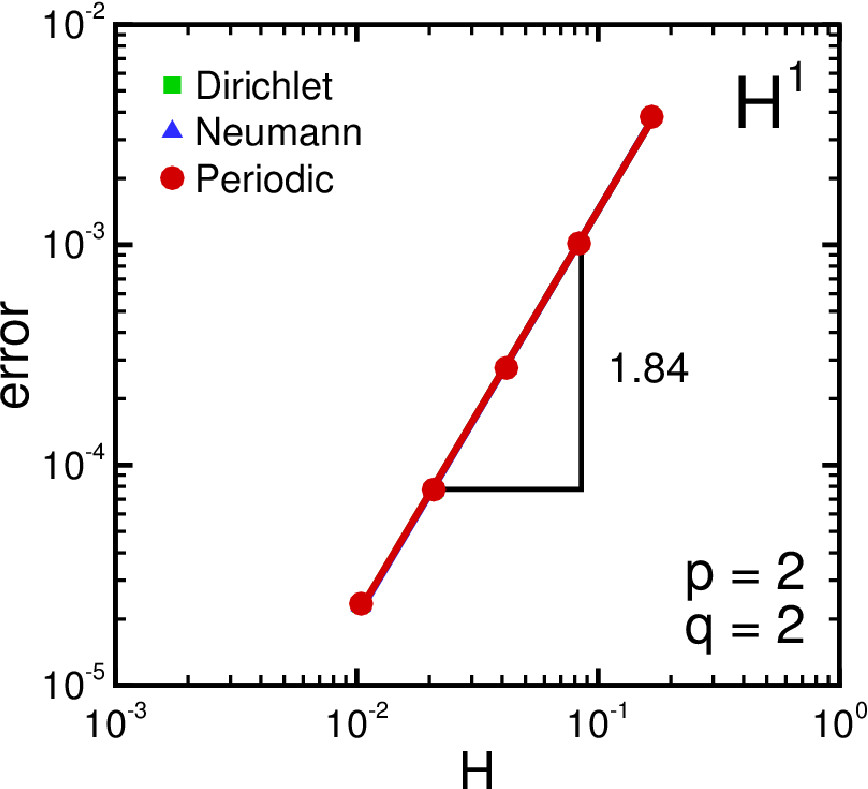} \hfill 
	\includegraphics[width=0.32\linewidth]{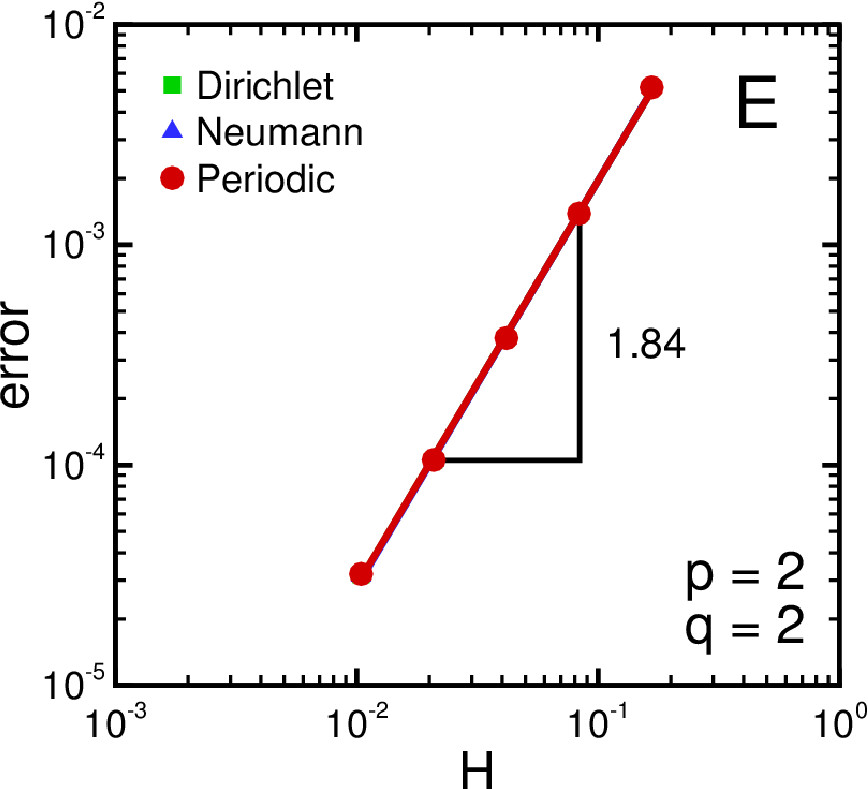}
	\caption{\textbf{Tapered cantilever under volume force.} Macroconvergence keeping the microdiscretization fixed, (first row) linear  shape functions $p$=$q$=$1$, (second row) quadratic shape functions $p$=$q$=$2$, (from left to right) $L^2$-, $H^1$-, and energy-norm.}
	\label{fig:diagramm_tapered_beam_p1}
\end{Figure}

The convergence orders for linear shape functions are displayed in the diagrams of the first row in Fig. \ref{fig:diagramm_tapered_beam_p1}. In all norms the optimal convergence order is achieved. Again, the results of the different coupling conditions coincide.

The results for quadratic shape functions in the second row of Fig. \ref{fig:diagramm_tapered_beam_p1} exhibit minor deviations from the full, nominal convergence orders (2.89 instead of 3 in the $L^2$-norm, and 1.84 instead of 2 in the $H^1$- and energy-norm). However, compared to the considerable order reduction for the square plate, the present deviations are small.
 
\begin{Figure}[htbp]
	\centering
	\includegraphics[height=6.0cm]{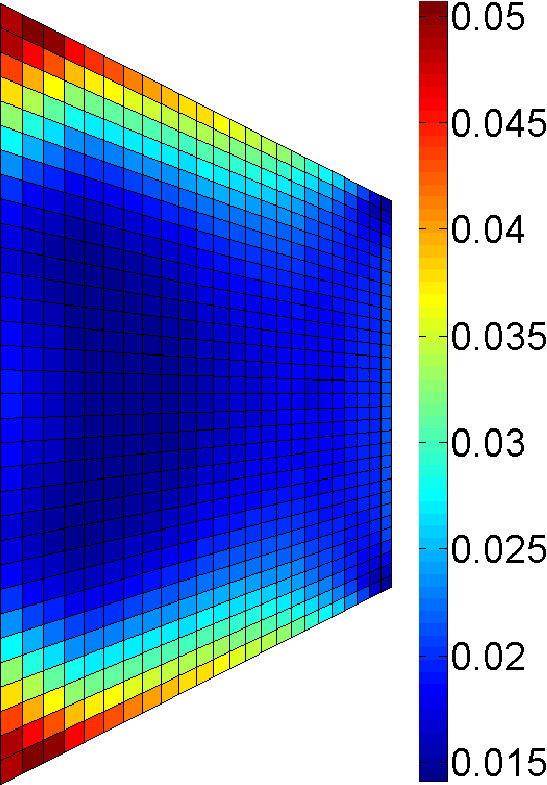} \hspace*{16mm}
	\includegraphics[height=6.0cm]{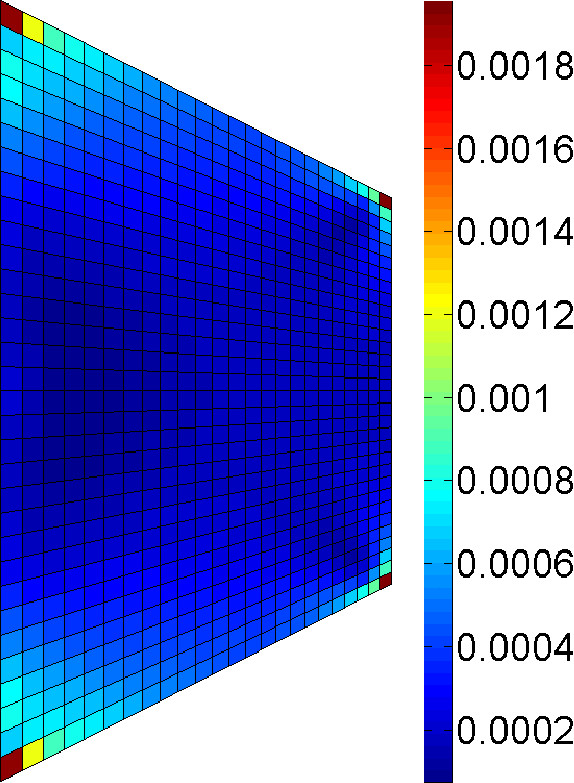}
	\caption{\textbf{Tapered cantilever under volume force.} Relative elementwise energy-error distribution on the macrodomain for (left) $p$=$q$=$1$ and (right) $p$=$q$=$2$.}
	\label{fig:error_distribution_tapered_beam}
\end{Figure}

Figure \ref{fig:error_distribution_tapered_beam} displays the relative elementwise error distribution on the macrodomain of the tapered cantilever for linear and for quadratic shape functions. The relative error in the energy-norm is computed as in the previous example. The error distribution reveals that there is no longer a singularity in the lower and upper left corner due to the bearing of the plate. This leads to significantly lower relative errors in these areas.

\subsection{Optimal uniform micro-macro refinement strategy}
\label{subsec:Optimal-mic-mac-refinement-strategy}

To investigate the optimal uniform micro-macro refinement strategy for linear and quadratic shape functions, the tapered cantilever is chosen as macro problem and the sine wave distribution is chosen as micro problem for their excellent regularity. The ratio of maximum to minimum Young's modulus is increased to 2.5.

\subsubsection{Linear shape functions}

For linear shape functions both on the micro and on the macro level the error in the $L^2$-norm is expected to converge in the order of $p+1=2$ on the macro and in the order of $2q=2$ on the micro level. Since micro and macro error converge in the same order, the micro mesh has to be refined in the same order as the macro mesh in order to achieve the optimal convergence order, $N_{mic}=(N_{mac})^{p+1/2q}=N_{mac}$. 

The error in the $H^1$-/energy-norm converges in the order of $p=1$ on the macro level and in the order of $2q=2$ on the micro level. Here the micro error converges in a higher order which means that the micro mesh does not have to be refined ''in the same order'' as the macro mesh, see \cite{JeckerAbdulle2016} Tab. 1 on p.5 and in the present work Tab.~\ref{tab:Best-mic-mac-RefinementStrategies}, $N_{mic}=(N_{mac})^{p/2q}= (N_{mac})^{1/2}$. 

\begin{Figure}[htbp]
	\centering
	\includegraphics[width=0.33\linewidth]{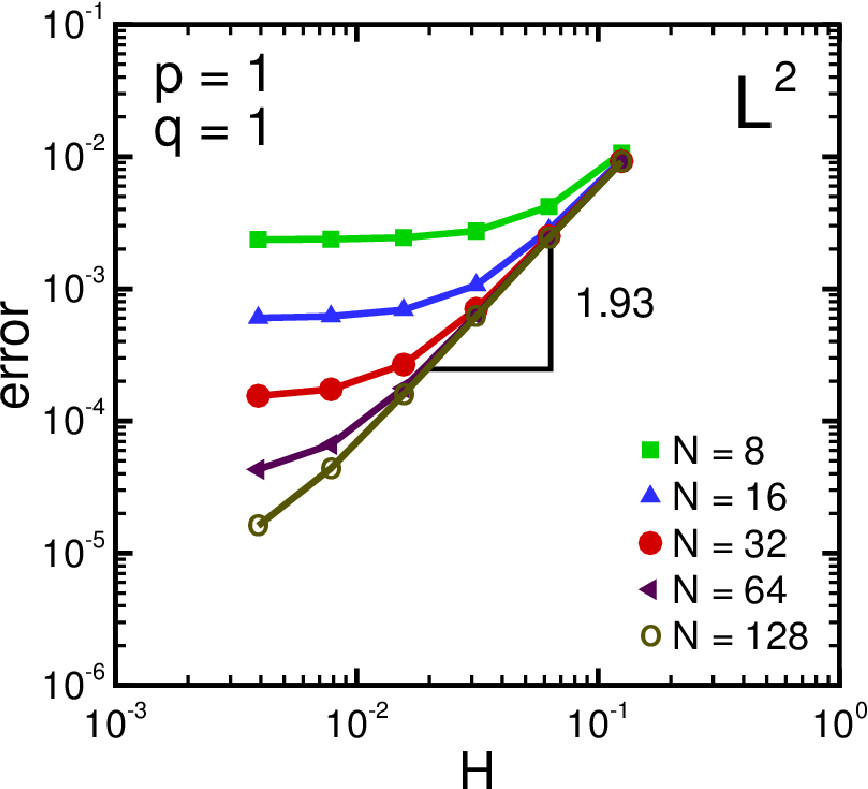} \hspace*{12mm}
	\includegraphics[width=0.33\linewidth]{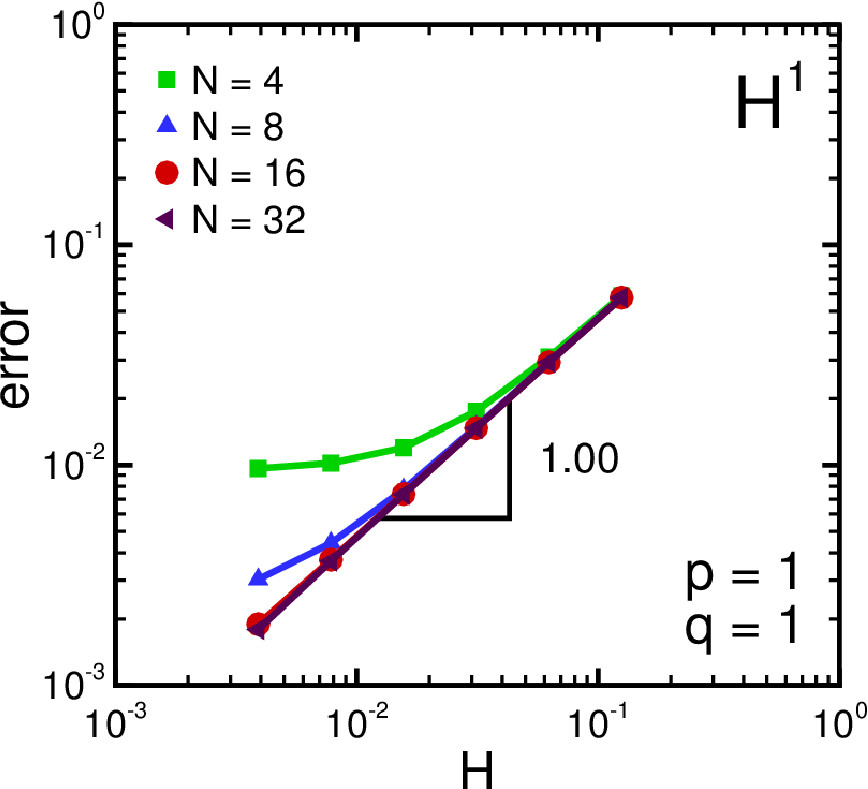}
	\caption{\textbf{Optimal uniform micro-macro refinement strategy.} (left) Error in the $L^2$-norm and (right) in the $H^1$-norm for linear shape functions, $p$=$q$=$1$ with {\color{black}$N=\epsilon/h$}.}
	\label{fig:refinement_strategy_p1q1}
\end{Figure}

Figure~\ref{fig:refinement_strategy_p1q1} shows the FE-HMM errors for different micro- and macrodiscretizations. For each line the microdiscretization is kept fixed {\color{black} where $N$ is the number of elements per edge on the micro domain,  $N=\epsilon/h$}; each marker in the diagrams denotes one macrodiscretization. If the micro mesh is not refined in the same order as the macro mesh, the error of the FE-HMM solution diverges from the line of optimal convergence in the $L^2$-norm. In the $H^1$-norm the micro mesh does not need not to be refined in the same order as the macro mesh so that the point where the single lines of fixed microdiscretizations diverge from the line of optimal convergence is shifted to finer macrodiscretizations.

\subsubsection{Quadratic shape functions}

The use of quadratic shape functions leads to optimal convergence orders of $p+1=3$ on the macro level and $2q=4$ on the micro level in the $L^2$-norm, $N_{mic}=(N_{mac})^{p+1/2q}=(N_{mac})^{3/4}$. In the $H^1$-norm we have $p=2$ on the macro level and again $2q=4$ on the micro level, $N_{mic}=(N_{mac})^{p/2q}= (N_{mac})^{1/2}$. 

\begin{Figure}[htbp]
	\centering
	\includegraphics[width=0.33\linewidth]{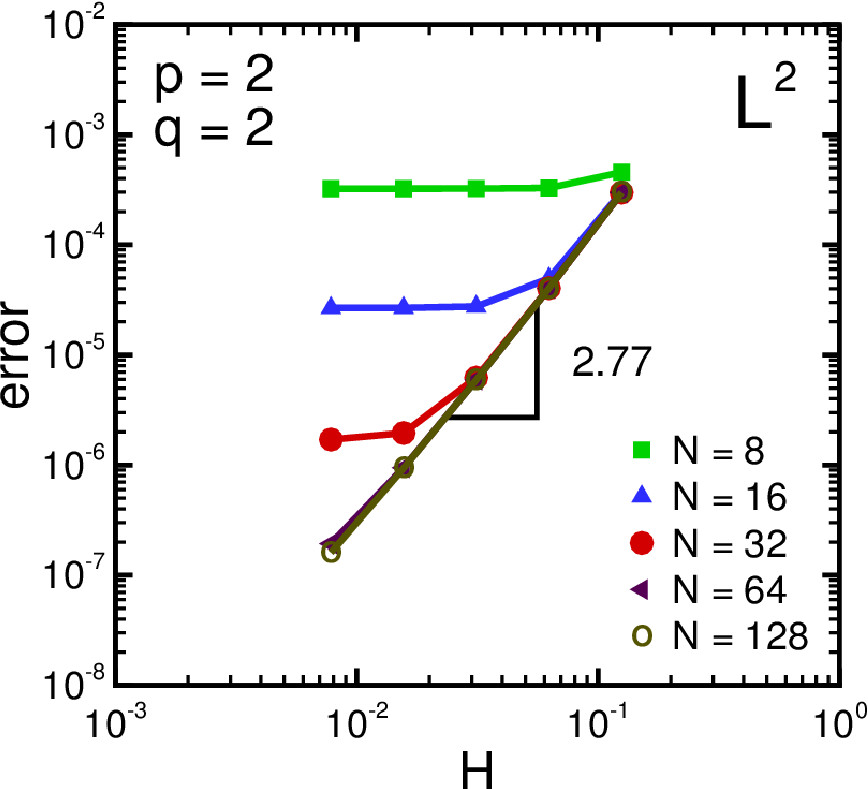} \hspace*{12mm}
	\includegraphics[width=0.33\linewidth]{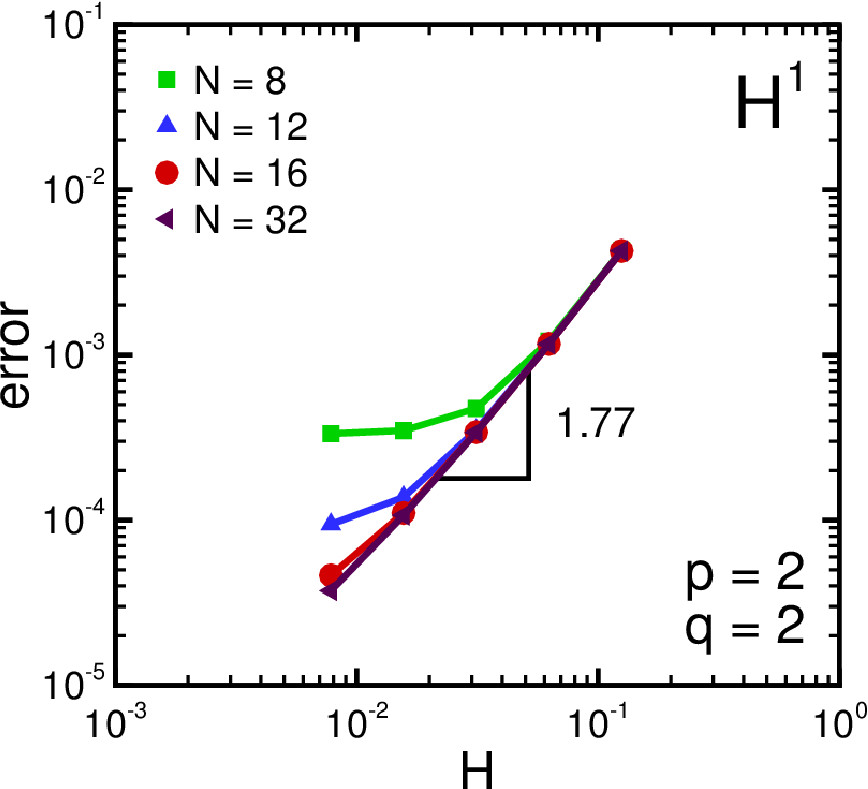}
	\caption{\textbf{Optimal micro macro refinement strategy.} Error in $L^2$-norm (left) and $H^1$-norm (right) with quadratic shape functions, $p$=$q$=$2$ with {\color{black}$N=\epsilon/h$}.}
	\label{fig:refinement_strategy_p2q2} 
\end{Figure}

The results of the optimal refinement strategy analysis is shown in Fig.~\ref{fig:refinement_strategy_p2q2} for quadratic shape functions. Analogue to the linear case the micro mesh has to be refined in a reduced order for optimal convergence in the $H^1$-norm compared to the $L^2$-norm.

\subsection{Accuracy of error estimation}
\label{subsec:error-estimation}

{\color{black} The accuracy of error estimation is compared with true error computation for the tapered cantilever subject to body forces along with a sine wave type Young's modulus distribution on the microscale and with PBC. Error estimation is based on the superconvergent patch recovery and on a simple averaging of elementwise stresses and strains. 
The results displayed in the diagrams of Fig.~\ref{fig:error_estimator} indicate that the estimated errors are in good agreement with the calculated errors for both linear as well as quadratic shape functions. 

\begin{Figure}[htbp]
	\centering
	\includegraphics[width=0.32\linewidth]{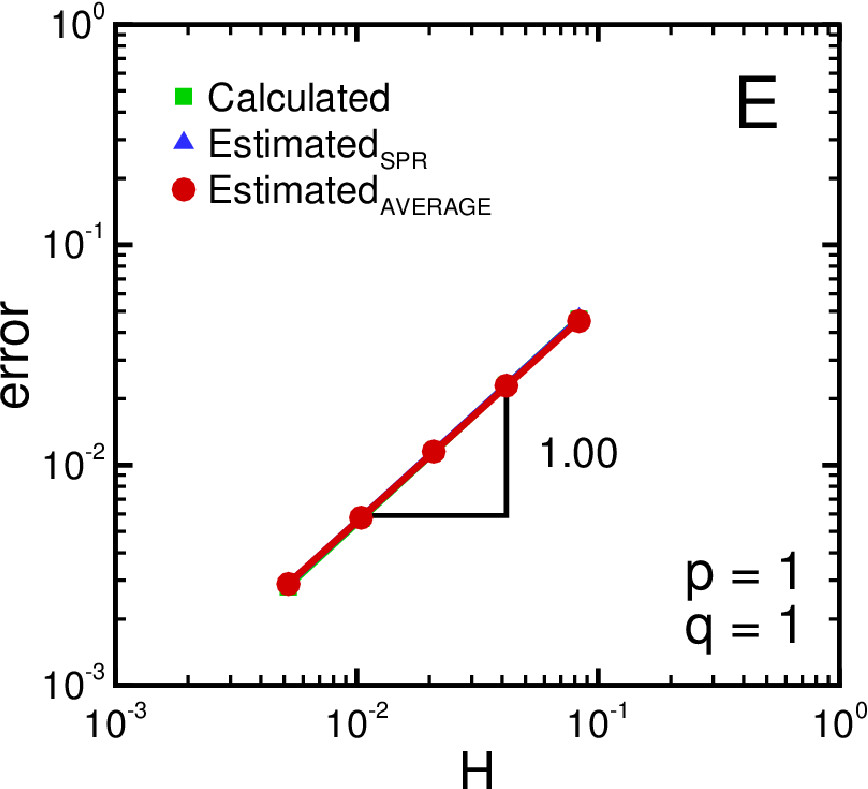} \hfill
	\includegraphics[width=0.32\linewidth]{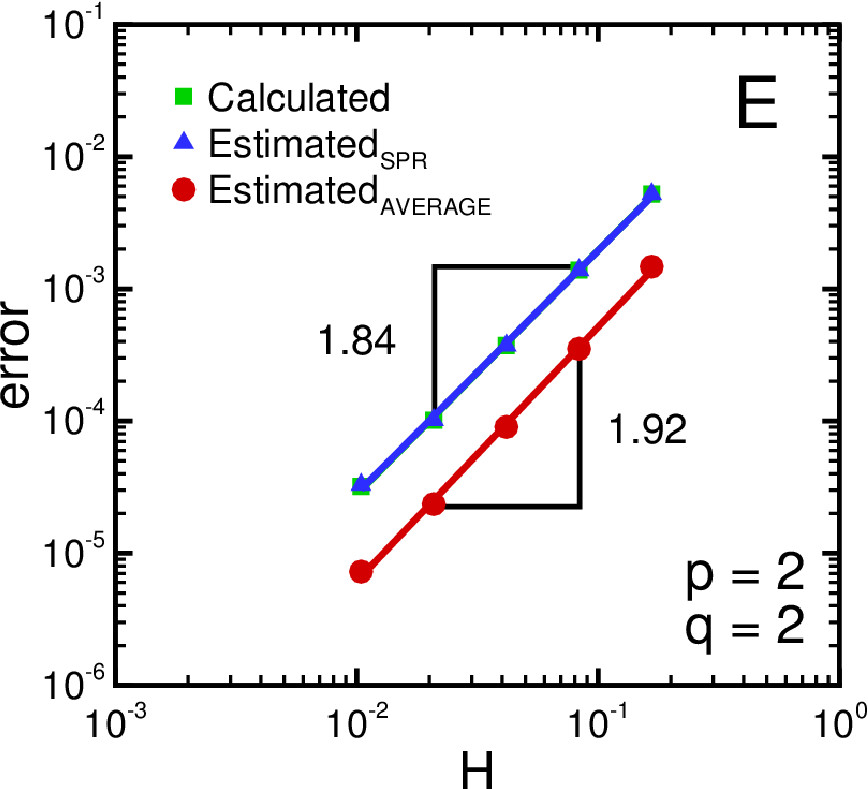} \hfill
	\includegraphics[width=0.32\linewidth]{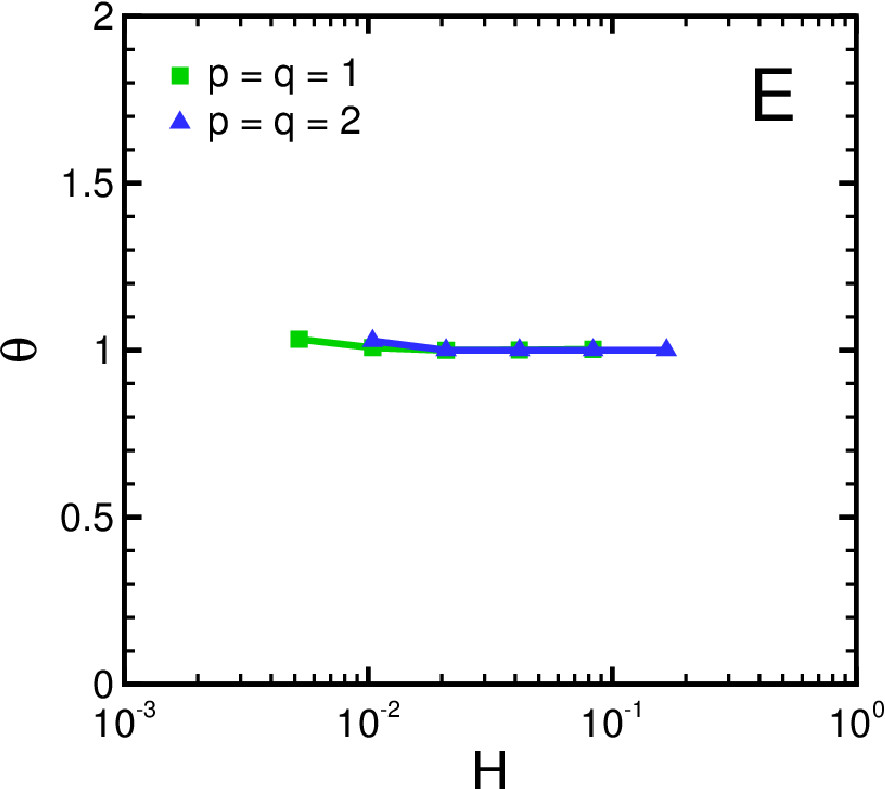}
	\caption{\color{black}\textbf{Accuracy of macro error estimation.} Calculated and estimated error in the energy-norm for (left) $p$=$q$=$1$, for (center) $p$=$q$=$2$, and (right) the corresponding effectivity index $\theta$.}
	\label{fig:error_estimator}
\end{Figure}

Notice that if the error estimator is restricted to the macro error part, the computation of the true total error 
along the lines of an optimal uniform micro-macro refinement strategy as verified in Sec.\ref{subsec:Optimal-mic-mac-refinement-strategy} cannot be replaced by suchlike error estimation. 
\\
We mention in passing that in the present example for the considered discretizations the macro error is orders of magnitude larger than the micro error, a result that will be quantitatively analyzed for the example in the consecutive Subsec.~\ref{subsec:error-deomposition}
}
 
{\color{black}

\subsection{Efficient decomposition of the true errors}
\label{subsec:error-deomposition}

The efficient decomposition of the computed error into its macro and micro parts according to Tab.~\ref{tab:rationale-for-error-decomposition-i} in Subsec.~\ref{subsec:error-computation} shall be demonstrated for the tapered beam with the sine wave type Young's modulus distribution along with PBC. Here, in contrast to Sec.~\ref{subsec:error-estimation}, a lineload is applied to the free end of structure.  
	
Steps 1.) and 2.) in Tab.~\ref{tab:rationale-for-error-decomposition-i} are carried out along with the corresponding reference solutions ($H\rightarrow 0$ by $1536 \times 2048$ elements, and $h\rightarrow 0$ by $1024 \times 1024$ elements). The microscopic error then directly follows as the difference between total error and macroscopic error. 

\begin{Figure}[htbp]
	\centering                             
	\includegraphics[width=0.32\linewidth]{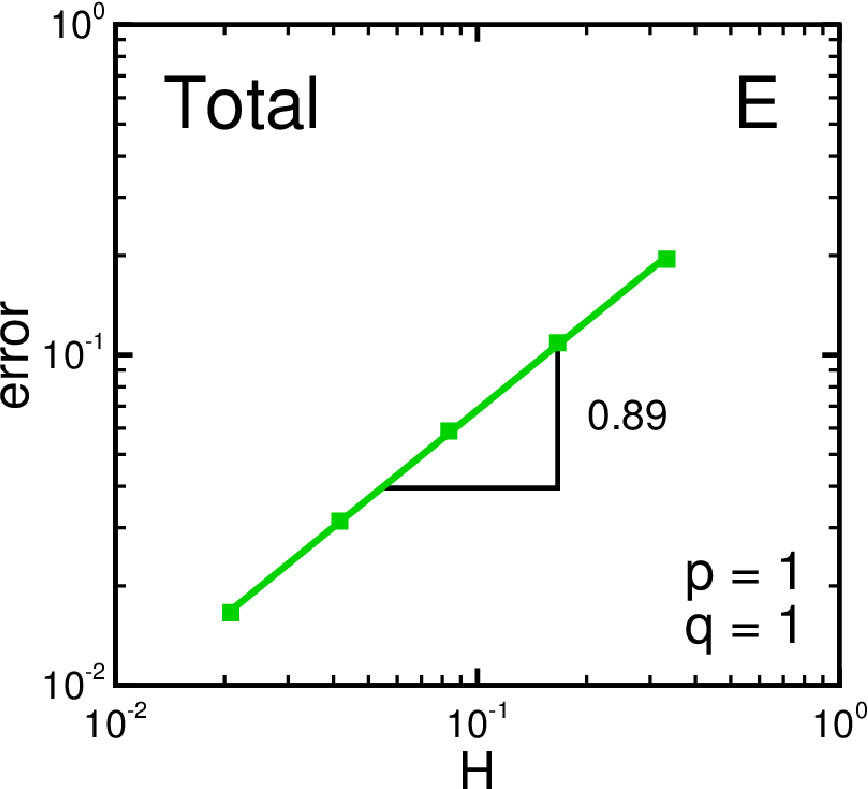} \hfill
	\includegraphics[width=0.32\linewidth]{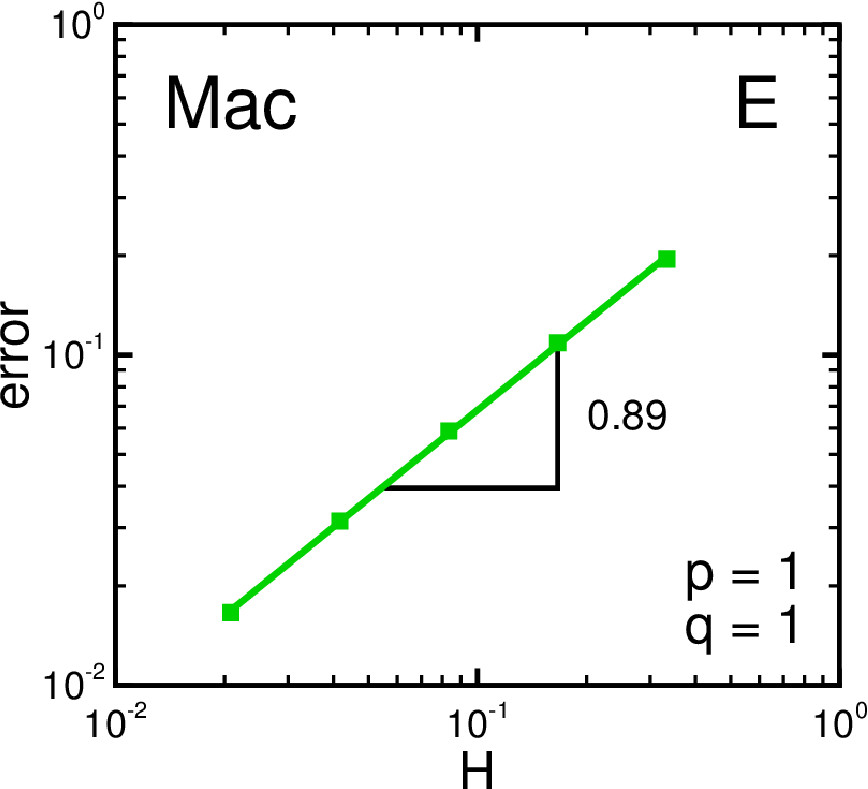} \hfill
	\includegraphics[width=0.32\linewidth]{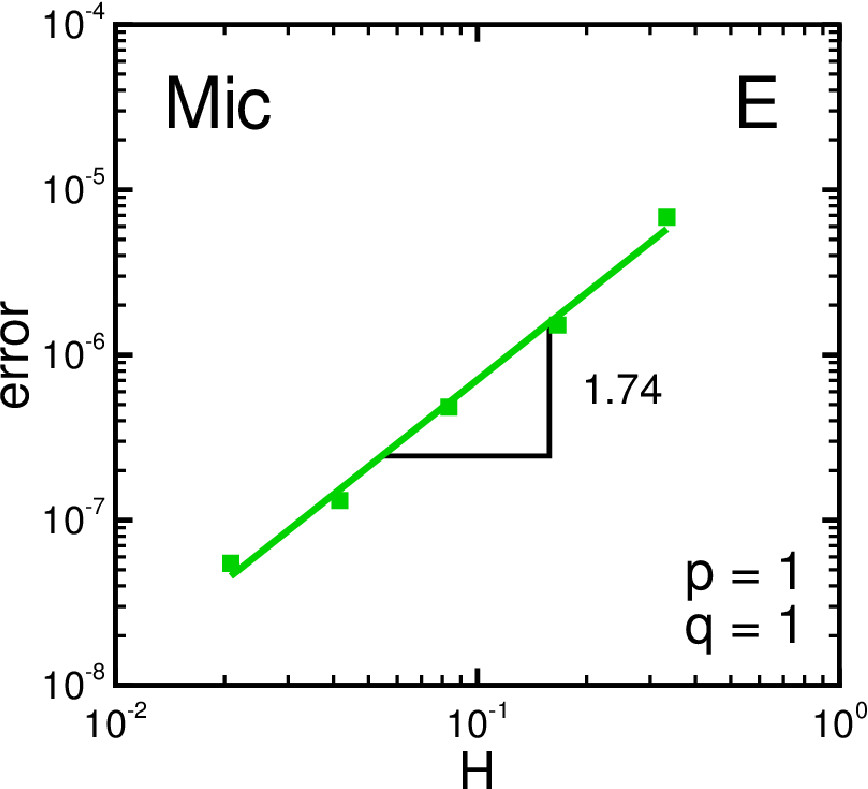}
	\caption{\color{black}\textbf{Tapered cantilever under line load.} Total error (left), macroscopic error (center) and microscopic error (right).}
	\label{fig:error_decomposition_tapered_beam}
\end{Figure}

Figure \ref{fig:error_decomposition_tapered_beam} shows the magnitude and convergence of the three different errors. The considerably larger values of the macro error compared to the micro error indicate the stronger influence of the macro discretization on the accuracy compared with the influence of the micro discretization. As a consequence, the total error and the macroscopic error are close together. The convergence rates are somewhat below the theoretical values which is true for each of the three errors. 

Figure \ref{fig:error_decomposition_tapered_beam_2} displays the distributions of the total error and the micro error. Recall that suchlike error decomposition clearly cannot be carried out by the error estimator operating on the macroscale, compare Remark 2 in Sec.~\ref{subsubsec:error_estimator}.

\begin{Figure}[htbp]
	\centering
	\includegraphics[height=5.7cm]{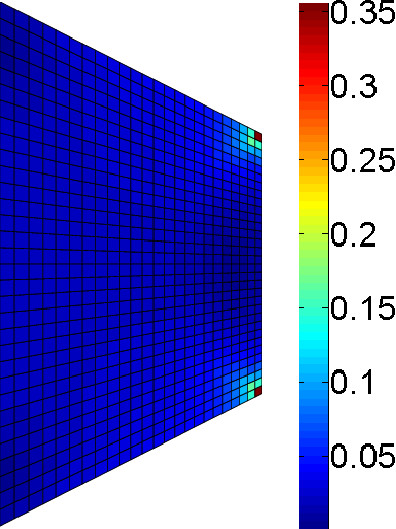} \hspace*{20mm}
	\includegraphics[height=6.3cm]{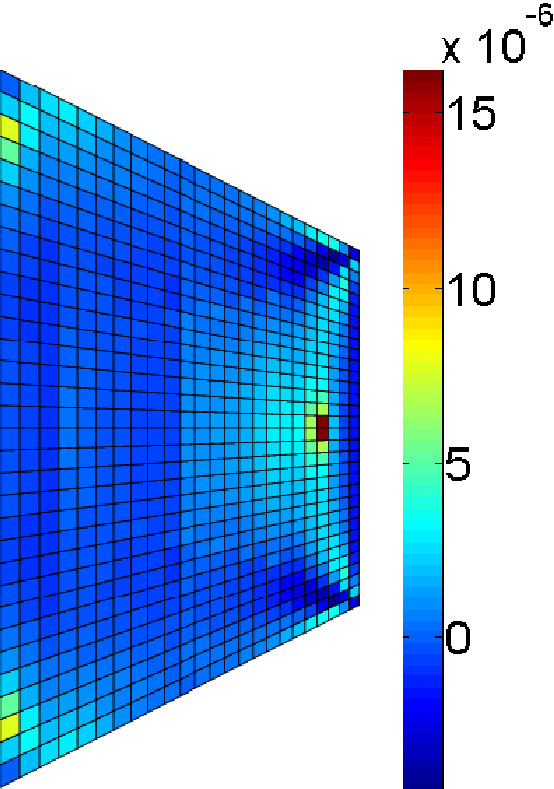} 
	\caption{\textbf{{\color{black}Tapered cantilever under line load.}} {\color{black}Relative elementwise energy-error distribution on the macrodomain for $p$=$q$=$1$ with (left) the total discretization error and (right) the micro error.}}
	\label{fig:error_decomposition_tapered_beam_2}
\end{Figure}
}

\subsection{Modeling error}

Next, the modeling error for Dirichlet coupling along with a noninteger ratio $\delta/\epsilon$ is examined, cf. \eqref{eq:Total-Error-estimate-L2}--\eqref{eq:Total-Error-estimate-Energy} along with \eqref{eq:ModelingError}. In the analysis the macro problem of the square plate under volume forces is used, the micro problem is the sine wave stiffness distribution. In order to investigate the convergence and show the modeling error, the macro discretization is continuously refined, while a very fine micro discretization ensures negligible micro errors. The calculations are run with $\delta/\epsilon = 1$ for PBC to indicate the optimal convergence of the macro problem without any modeling error. For Dirichlet coupling the cases {\color{black}$\delta/\epsilon \in \{1.0, 1.1, 5/3, 2.0\}$ are considered somewhat increasing the range of $\delta/\epsilon$} in  \cite{JeckerAbdulle2016}. The reference solution for error calculation is obtained for PBC and a very fine macro mesh.

\begin{Figure}[htbp]
	\centering
	\includegraphics[width=0.33\linewidth]{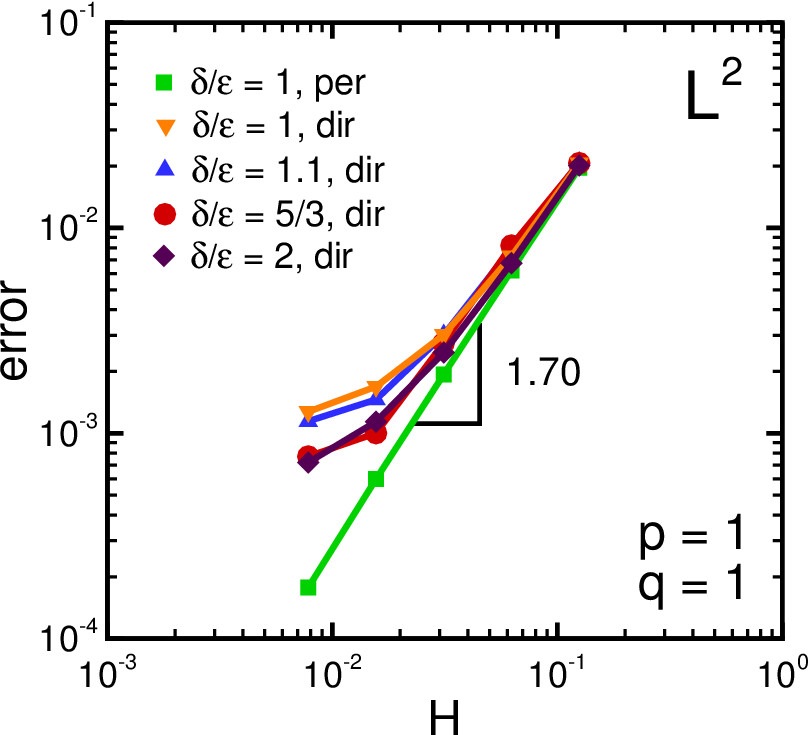}
	\caption{\textbf{Error convergence for different ratios of $\delta/\epsilon$.} Calculated errors for periodic coupling with $\delta/\epsilon = 1$ and for  Dirichlet coupling with {\color{black}$\delta/\epsilon \in \{1.0, 1.1, 5/3, 2.0\}$}, {\color{black}$\epsilon=0.005$}.}
	\label{fig:delta_epsilon}
\end{Figure}

Figure \ref{fig:delta_epsilon} depicts the results of the error calculations. For $\delta/\epsilon = 1$ and PBC the expected constant order is observed; it is the case already reported in Fig.~\ref{fig:diagramm_cantilever_beam_p1} for $p$=$q$=1 showing a somewhat reduced convergence order of 1.70. In contrast to the constant convergence for the reference case of $\delta/\epsilon = 1$ along with PBC the case of {\color{black} Dirichlet coupling along with various $\delta/\epsilon$ ratios} exhibit an offset between the corresponding curves and the reference solution. This offset indicates the modeling error, which becomes increasingly dominant compared to the macro error for macro mesh refinement.  
{\color{black} Notice that the modeling error increases for an increasing ratio $\epsilon/\delta$ for $\delta > \epsilon$ in agreement with \eqref{eq:ModelingError}}.
\section{Summary and conclusions}
\label{sec:summary}

The aim of the present work was the numerical analysis of energetically consistent micro-coupling conditions in the homogenization framework of a two-scale finite element method. The obtained results are valid for FE-HMM and FE$^2$ for the coincidence of the methods, which was shown in \cite{EidelFischer2018}. The main results shall be summarized.  
 
\begin{enumerate}
 \item The {\bf micro error {\color{black} convergence} analysis} for different micro problems, micro-macro coupling conditions and polynomial orders of shape functions has led to the following results.
 \begin{enumerate}
        \item We have clarified the distinction between the micro-error measured on the microscale with the theoretical convergence order of $q$+$1$ and $q$ in the $L^2$-norm and the $H^1$-/energy-norm, respectively, and the micro error as propagated to the macroscale with the -somewhat surprising- convergence order of $2q$ in all ($L^2$-, $H^1$-, energy-) norms. 
	\item For sufficiently regular micro problems the a priori error estimates of FE-HMM have been confirmed for each of the considered ($L^2$-, $H^1$-, energy-) norms.
	The micro-coupling conditions show no significant deviation from each other in the measured convergence order.
	\item The regularity of a micro BVP requires a microstructure with smooth distribution of the heterogeneous material parameters (here: Young's modulus). Then the contrast of maximum to minimum material parameters does not influence the convergence order. Vice versa, a stiffness-jump at interfaces in the RVE lowers the convergence order for linear shape functions, and quadratic shape functions do not cure the order reduction. 
   \end{enumerate}
\item The {\bf constant traction (Neumann) condition} 
   \begin{enumerate} 
        \item Two methods for Neumann conditions have been compared, the recently introduced semi-Dirichlet ansatz \cite{JaviliSaebSteinmann2017} with the mass-type perturbation technique for regularization \cite{Miehe-Koch-2002}. The two methods coincide in the goal but differ in the methodic procedure to remove rigid body motions from the RVE and the corresponding singularity of the stiffness matrix.
	\item Both methods are accurate in fulfilling the condition of constant traction. The approach of Miehe \& Koch turns out to be remarkably insensitive to the particular choice of the perturbation parameter in a wide parameter range. It is most simple to implement and fast. The Semi-Dirichlet ansatz carries out explicit static condensation of rigid body motions by additional Dirichlet conditions of the RVE. Since this approach requires an iterative solution, it is more expensive than the perturbation technique. 
    \end{enumerate}   
\item The {\bf macro error {\color{black} convergence} analysis} for two different macro problems applying different micro-macro coupling conditions and linear as well as quadratic shape functions has led to the following results.
   \begin{enumerate}
	\item For fully regular macro BVPs the error estimates have been confirmed in all norms. Singularities spoil the convergence such that the order is below the theoretical estimate for linear shape functions, which implies no improvement in the order for quadratic shape functions. 
	\item Two macroproblems underpin the aforementioned statements; a square cantilever which does not achieve the full order due to singularities in the corners of the bearing. At these points and their direct neighborhood the error are maximal. A tapered cantilever plate avoids by its geometry these singularities and enables therefore full order in agreement with the a priori error estimates, which is true for linear shape functions and for quadratic shape functions along with only minor deviations.  
	\item The choice of the coupling conditions on the RVE does not affect the macro convergence order nor the quantitative macro error. 
   \end{enumerate}
\item {\bf Error estimation and optimal mesh refinements} \\
     {\color{black} The recovery-type, superconvergent error estimator of \cite{SPR} was implemented on the macroscale; it exhibits the following properties:}
   \begin{enumerate} 
	\item {\color{black} The error estimator is accurate as indicated by an efficiency index close to unity; the estimated error almost equals the computed discretization error.}  
	\item {\color{black} The optimal uniform micro-macro refinement strategies directly following from the a priori error estimates were confirmed for linear and quadratic shape functions. These refinement strategies are of considerable practical value since they enable the optimal convergence of the total error while keeping the numerical effort minimal.} 
	\item {\color{black} Although the error estimator on the macroscale exhibits different estimates at (macro-, micro-) discretizations $(H, h)$ and $(H, h/2)$, it merely measures macro discretization errors. The difference in the error estimations for various $h$ is on the macroscale not a discretization error but a modeling error, since indeed the same type of constitutive law is used on the microscale but for different material parameters, which depend on the microdiscretization $h$.} 
    \end{enumerate}
\item {\bf Analysis of the modeling error.} The modeling error for Dirichlet coupling along with various $\delta/\epsilon$ ratios was identified and made measurable by uniform macro mesh refinements along with fine micro meshes. The simulation results underpin the a priori estimate {\color{black} in that the modeling error increases for an increasing $\epsilon/\delta$ ratio with $\delta>\epsilon$. Moreover, for $H \rightarrow 0$ and $h \rightarrow 0$ the modeling error persists as a discretization-independent residual, again consistent with the estimate.} 
\end{enumerate}
 
\bigskip

{\bf Acknowledgements.} Bernhard Eidel acknowledges support by the Deutsche Forschungsgemeinschaft (DFG) within the Heisenberg program (grant no. EI 453/2-1). Simulations were performed with computing resources granted by RWTH Aachen University under project ID prep0005. 
 
\bigskip

{\bf Declaration of Interest.} None. 

\vfill
\newpage

\begin{appendix}

\addcontentsline{toc}{section}{Appendix}
\renewcommand{\thesubsection}{\Alph{section}.\arabic{subsection}}
\renewcommand{\theequation}{\Alph{section}.\arabic{equation}}
\renewcommand{\thefigure}{\Alph{section}.\arabic{figure}}
\renewcommand{\thetable}{\Alph{section}.\arabic{table}}
\newcommand {\ssectapp}{
                        \setcounter{equation}{0}
                        \setcounter{figure}{0}
                        \setcounter{table}{0}
		                \subsection
                        }

\setcounter{equation}{0}

%\appendix
 
\sect{\color{black} Appendix}
\label{sect:Miscellaneous}

%---------------------------------------------------------------------------------------------------------
\subsection{Derivation of the micro-to-macro stiffness transformation matrix}
\label{subsec:TrafoMatrix}

The derivation of the macro element stiffness matrix part $\bm k^{e,mac}_{IJ}$ in \eqref{eq:k-mac-element-5} shall be detailed:
 \begin{eqnarray}          
\bm k^{e,mac}_{IJ} 
&=& B^e_H\left[\bm N_I^H, \bm N_J^H\right]  \nonumber\\
&=& \sum_{l=1}^{N_{qp}} \dfrac{\omega_{K_l}}{|K_{\delta}|} 
           \int_{K_{\delta}} \left(\bm L \bm u^{h(I,x_i)}_{K_{\delta}} \right)^T \mathbb{A}^{\epsilon} \, \bm L \bm u^{h(J,x_i)}_{K_{\delta}} \, dV 
           \nonumber \\
&=& \sum_{l=1}^{N_{qp}} \dfrac{\omega_{K_l}}{|K_{\delta}|} 
           \int_{K_{\delta}} \left(\bm L \sum_{m=1}^{M_{mic}} \bm N^h_{m, K_{\delta}} \, \bm d^{h(I,x_i)}_{m} \right)^T \mathbb{A}^{\epsilon} 
                             \, \bm L \sum_{n=1}^{M_{mic}} \bm N^h_{n,K_{\delta}} \, \bm d^{h(J,x_i)}_{n} \, dV 
           \nonumber \\
&=& \sum_{l=1}^{N_{qp}} \dfrac{\omega_{K_l}}{|K_{\delta}|} 
           \sum_{T \in \mathcal{T}_h} \, \int_{T} \, \left(\bm L \sum_{m=1}^{n_{node}} \bm N^h_{m,K_{\delta}} \, \bm d^{h(I,x_i)}_{m} \right)^T \mathbb{A}^{\epsilon} 
                             \, \bm L \sum_{n=1}^{n_{node}} \bm N^h_{n,K_{\delta}} \, \bm d^{h(J,x_i)}_{n} \, dV 
           \nonumber \\          
&=& \sum_{l=1}^{N_{qp}} \dfrac{\omega_{K_l}}{|K_{\delta}|} 
           \sum_{T \in \mathcal{T}_h} \, \int_{T} \, \left(\bm L \sum_{m=1}^{n_{node}} \bm N^h_{m,K_{\delta}} \, \bm d^{h(I,x_i)}_{m} \right)^T \mathbb{A}^{\epsilon} 
                             \, \bm L \sum_{n=1}^{n_{node}} \bm N^h_{n,K_{\delta}} \, \bm d^{h(J,x_i)}_{n} \, dV 
           \nonumber \\           
&=& \sum_{l=1}^{N_{qp}} \dfrac{\omega_{K_l}}{|K_{\delta}|} 
           \left(\bm d^{h(I,x_i)}\right)^T 
           \sum_{T \in \mathcal{T}_h} \, \int_{T} \, \sum_{m=1}^{n_{node}} \sum_{n=1}^{n_{node}}
                             \left(\bm L \bm N^h_{m, K_{\delta}} \right)^T \mathbb{A}^{\epsilon} 
                             \, \bm L \bm N^h_{n, K_{\delta}} \, dV \, \bm d^{h(J,x_i)}_{n} \,  
           \nonumber \\ 
&=& \label{eq:k-mac-element-2} 
           \sum_{l=1}^{N_{qp}} \dfrac{\omega_{K_l}}{|K_{\delta}|} 
           \left(\bm d^{h(I,x_i)}\right)^T 
           \, \sum_{T \in \mathcal{T}_h} \,
           \Big(
           \sum_{m=1}^{n_{node}} \sum_{n=1}^{n_{node}} 
                 \, \int_{T} \bm B_m^{e,T} \mathbb{A}^{\epsilon} \, \bm B^e_n \, dV 
           \bm d^{h(J,x_i)}_n  \Big) \\
           &=& \label{eq:k-mac-element-3}
                      \sum_{l=1}^{N_{qp}} \dfrac{\omega_{K_l}}{|K_{\delta}|} 
           \sum_{T \in \mathcal{T}_h} \,
           \left(
           \sum_{m=1}^{n_{node}} \sum_{n=1}^{n_{node}} 
                \left(\bm d^{h(I,x_i)}_m \right)^T \, \bm k^{e,mic}_{mn} \bm d^{h(J,x_i)}_n
            \right) \nonumber \\
          &=&  \label{eq:k-mac-element-4}
           \sum_{l=1}^{N_{qp}} \dfrac{\omega_{K_l}}{|K_{\delta}|} 
           \left(\bm d^{h(I)} \right)^{T} \, 
           \sum_{T \in \mathcal{T}_h}   
           \bm k_{K_{\delta}}^{e,mic} \,  \bm d^{h(J)} 
            \nonumber \\      
          &=& \label{eq:k-mac-element-5app}
           \sum_{l=1}^{N_{qp}} \dfrac{\omega_{K_l}}{|K_{\delta}|} 
           \, \left( \bm d^{h(I)} \right)^T \, \bm K^{mic}_{K_{\delta}} \,  \bm d^{h(J)}  \, , \nonumber 
\end{eqnarray} 
where $\bm d^{h(I)}=\left(\,\bm d^{h(I,x_1)} | \bm d^{h(I,x_2)} | \bm d^{h(I,x_3)} \, \right)$ for $n_{dim}=3$. The assembly of $\bm k^{e,mac}_{IJ}$ results in  
\begin{eqnarray}
 \bm k^{e,mac}_{K} &=& \sum_{l=1}^{N_{qp}} \dfrac{\omega_{K_l}}{|K_l|} 
               \, \, \bm T^{T}_{K_l}\, \bm K^{mic}_{K_l} \, \bm T_{K_l}     \\
 \mbox{with} \quad \bm T_{K_l} &=& \bigg[ \Big[ \big[ \bm d^{h(I,x_i)} \big]_{i=1,\ldots,n_{dim}} \Big]_{I=1, \ldots, N_{node}} \bigg]   \,.        
\end{eqnarray}

\vfill
\newpage

\subsection{Definition of norms}
\label{subsec:Definition-norms}

The norms used in the present work are defined according to 
\begin{eqnarray}
\mbox{$L^2$-norm:} \quad  || \bm u ||_{L^2(\Omega)} &:=& \sqrt{\int_{\Omega} \bm u : \bm u \, dV } \, ,
\label{eq:L2-norm} \\
\mbox{$H^1$-norm:} \quad  || \bm u ||_{H^1(\Omega)} &:=&  
\sqrt{ \left( \sum_{i,j=1}^d \int_{\Omega} \left( \dfrac{\partial u_i}{\partial x_j} \right)^2 \, dV
	+  \sum_{i=1}^d \int_{\Omega} \left(u_i\right)^2 \, dV \right) }\, ,
\label{eq:Hilbert-norm} \\
\mbox{energy-norm:} \quad  ||\,\bm u\,||_{A(\Omega)} &:=& \sqrt{\int_{\Omega} \mathbb{A} \, \bm \varepsilon(\bm u) : \bm \varepsilon(\bm u) \, dV} 
= \sqrt{\bm d^T \, \bm K \, \bm d }\, .
\label{eq:energy-norm}
\end{eqnarray}

\subsection{{\color{black}Direct implementation of Dirichlet and periodic coupling conditions}}
\label{subsec:DirectSolutionCouplingConditions}

{\color{black}
	
	When Dirichlet or periodic boundary conditions are implemented in a direct manner without using Lagrange multipliers the micro system of equations reads as
	
	\begin{equation}
	\label{eq: sys_direct}
		\bm K^{mic}_{K_l} \bm d^{h(I,x_i)} = \bm 0
	\end{equation}
	
	for a macroscopic unit displacement state $(I,x_i)$, a notation that is dropped in the following for notational convenience along with the 
	subscript $K_l$ indicating quadrature point $l$ of macro element $K$. %, so that the system of equations is
	
	For Dirichlet coupling conditions the micro displacements $\bm d^h$ are known since they directly follow from the macroscopic displacement field. For that 
	reason we decompose the microscopic displacement vector
	
	\begin{equation}
		\bm d^h = [\bm d^h_D, \ \bm d^h_F]^T
	\end{equation} 
	
	into the known displacements $\bm d^h_D$ and an unknown part $\bm d^h_F$. Then the system of equations follows the form
	
	\begin{equation}
	\left[ \begin{array}{cc}
		\bm K^{mic}_{DD} & \bm K^{mic}_{FD} \\
		\bm K^{mic}_{DF} & \bm K^{mic}_{FF} \\
	\end{array} \right]
	\left[ \begin{array}{c}
		\bm d^{h}_D \\
		\bm d^{h}_F \\
	\end{array}\right]
	=
	\left[ \begin{array}{c}
		\bm 0   \\
		\bm 0    \\
	\end{array}\right] \, .
	\end{equation}	
	
	It follows that
	
	\begin{equation}
		\bm K^{mic}_{DF} \bm d_D^h + \bm K^{mic}_{FF} \bm d_F^h = \bm 0 
	\end{equation}
	
	which yields the unknown micro displacements according to 
	
	\begin{equation}
		\bm d_F^h = - \left( \bm K^{mic}_{FF} \right)^{-1} \left(\bm K^{mic}_{FD} \bm d_D^h \right) \, .
	\end{equation}
	
	For periodic coupling conditions the direct implementation accounts for the fact that not the micro displacements but the fluctuations between macroscopic and microscopic 
	displacement fields are periodic on opposite boundaries. 
	
	First we introduce the fluctuations on the microlevel as the difference between macroscopic and microscopic displacements
	
	\begin{equation}
		\tilde{\bm d}^h = \bm d^H - \bm d^h \, .
	\end{equation}
	
	Since the system of equations from \eqref{eq: sys_direct} only contains the microscopic displacements, we decompose them into the homogeneous deformation part following from
	the macroscopic displacements $\bm d^H$ and their fluctuations $\tilde{\bm d}^h$
	
	\begin{equation}
	\bm d^h = \bm d^H + \tilde{\bm d}^h \, ,
	\end{equation}
	
	which is inserted into \eqref{eq: sys_direct} and yields the solution 
	
	\begin{equation}
	       \label{eq:app:Solving-for-fluctuations}
		\tilde{\bm d}^h = \left(\bm K^{mic}\right)^{-1} \left( - \bm K^{mic} \bm d^H \right) \, .
	\end{equation}
	
	Equation \eqref{eq:app:Solving-for-fluctuations} is solved accounting for periodic boundaries. Rigid body motions are removed from the system 
	simply by fixing one arbitrary node in each direction of space. A convenient choice for periodic structures is to set the displacements 
	of a node in the center of an RVE to zero.
 
}

\end{appendix}

\bibliographystyle{plainnat}
%\bibliographystyle{plaindin}
%\bibliographystyle{plaindin_shortname2}
%\bibliographystyle{elsarticle-num-names}
% ------- bib-datei --------------
%\bibliography{js_master_2009}
%\begin{appendix}
%\include{appendixA}
%\end{appendix}
\end{document}